\documentclass[a4paper, reqno, 10pt]{amsart}
\usepackage{amssymb}
\usepackage{amsmath}
\usepackage{extarrows}
\usepackage{bm}
\usepackage[centering]{geometry}
\usepackage{enumitem}
\usepackage{longtable}
\usepackage{multirow}
\usepackage{array}
\usepackage{etex}
\usepackage{pictex}
\usepackage{graphics}
\usepackage{amscd}
\usepackage[all]{xy}

\numberwithin{equation}{section}
\theoremstyle{definition}
\newtheorem{defn}{Definition}[section]
\newtheorem{rem}[defn]{Remark}

\newtheorem{exm}[defn]{Example}

\theoremstyle{plain}
\newtheorem{cor}[defn]{Corollary}
\newtheorem{thm}[defn]{Theorem}
\newtheorem{lem}[defn]{Lemma}
\newtheorem{deffact}[defn]{Definition - Fact}
\newtheorem{deflem}[defn]{Definition - Lemma}

\newtheorem{prop}[defn]{Proposition}
\newtheorem{fact}[defn]{Fact}
\newtheorem{defn-thm}[defn]{Definition-Theorem}

\newtheorem*{Freyd-Mitchell}{Freyd-Mitchell embedding Theorem}
\def\Fib{\operatorname{Fib}}
\def\CoFib{\operatorname{CoFib}}
\def\Weq{\operatorname{Weq}}
\def\Coker{\operatorname{Coker}}
\def\hoc{\operatorname{hoCoker}}
\def\Hom{\operatorname{Hom}}
\def\Ext{\operatorname{Ext}}
\def\Ker{\operatorname{Ker}}
\def\hok{\operatorname{hoKer}}
\def\Id{\operatorname{Id}}
\def\TFib{\operatorname{TFib}}
\def\TCoFib{\operatorname{TCoFib}}
\def\Ho{\operatorname{Ho}}
\def\Im{\operatorname{Im}}
\def\tri{\operatorname{\Delta}}
\def\X{\mathcal{X}}
\def\Y{\mathcal{Y}}
\def\A{\mathcal{A}}
\def\T{\mathcal{T}}
\newcommand{\xra}{\xlongrightarrow}
\title[Model structures on triangulated categories with proper class of triangles]{Model structures on triangulated categories \\ with proper class of triangles}
\author[Jian Cui, Pu Zhang] {Jian Cui, Pu Zhang$^*$ \\ \\  School of Mathematical Sciences \\
Shanghai Jiao Tong University,  \ Shanghai 200240, \ China }
\thanks{$^*$ Corresponding author}
\thanks{provinceanying$\symbol{64}$sjtu.edu.cn \ \ \ \ pzhang$\symbol{64}$sjtu.edu.cn}
\thanks{This work is supported by National Natural Science Foundation of China under Grant 12131015, and by Natural Science Foundation of Shanghai under Grant 23ZR1435100.}

\begin{document}
\maketitle
\begin{abstract} In contrast with the Hovey correspondence of abelian model structures from two compatible complete cotorsion pairs,
Beligiannis and Reiten give a construction of model structures on abelian categories from one hereditary complete cotorsion pair.
The aim of this paper is to extend this result to triangulated categories together with a proper class $\xi$ of triangles.
There indeed exist non-trivial
proper classes of triangles, and a proper class of triangles is not closed under rotations, in general. This is quite different from the class of
all triangles. Thus one needs to develop a theory of triangles in $\xi$ and hereditary complete cotorsion pairs in a triangulated category $\T$  with respect to $\xi$.
The Beligiannis - Reiten correspondence between weakly $\xi$-projective model structures on $\T$ and hereditary complete cotorsion pairs $(\X, \Y)$ with respect to $\xi$
such that the core $\omega = \X \cap \Y$ is contravariantly finite in $\T$ is also obtained. To study the homotopy category of
a model structure on a triangulated category, the condition in Quillen's
Fundamental theorem of model categories needs to be weakened,
by replacing the existence of pull-backs and push-outs by homotopy cartesian squares.

\vskip5pt


\end{abstract}

\section {\bf Introduction}

\subsection{Model structures and complete cotorsion pairs} \ Connections  between model structures and cotorsion pairs have been initiated by M. Hovey [Hov2] and A. Beligiannis and I. Reiten [BR].
There is an one-one correspondence between abelian model structures and the Hovey triples ([Hov2]). A Hovey triple gives two compatible complete cotorsion pairs; and
if the two compatible complete cotorsion pairs are hereditary, then the  Hovey triple is uniquely determined by the two cotorsion pairs ([G2]).
This Hovey correspondence is extended by J. Gillespie [G1] (see also [G3]) as an one-one correspondence
between exact model structures and the Hovey triples, on a weakly idempotent complete exact category. See also J. \v{S}t'ov\'{\i}\v{c}ek \cite {S}.
For a triangulated category with a proper class $\xi$ of triangles, X. Y. Yang [Y] considers
cotorsion pairs  with respect to $\xi$, and obtains the $\xi$-triangulated version of the Hovey correspondence.
H. Nakaoka and Y. Palu [NP] introduce an extriangulated category, which is a common generalization of an exact category and a triangulated category,
and generalize the Hovey correspondence to this setting up; moreover, remarkably, the homotopy category is proved to be triangulated, even if the two compatible complete cotorsion pairs are not necessarily hereditary. See [NP, Theorem 6.20], and [G3, Theorem 6.34].

\vskip5pt

On the other hand, there is an one-one correspondence between weakly projective model structures (or called the $\omega$-model structures) on abelian category $\A$ and hereditary complete cotorsion pairs
with the core $\omega$ contravariantly finite in $\A$ ([BR, VIII, 4.2, 4.13]). This $\omega$-model structure is fibrant (i.e., each object is fibrant)
and involves only one hereditary complete cotorsion pair.
This Beligiannis - Reiten correspondence is extended in [CLZ] to a weakly idempotent complete exact category.
A weakly projective model structure is not an exact model structure, in general; and the two approaches give the same result
if and only if $\A$ has enough projective objects and the model structure is projective (i.e., it is exact and each object is fibrant).
The Beligiannis - Reiten correspondence is also generalized to weakly idempotent complete extriangulated categories by J. S. Hu, D. D. Zhang and P. Y. Zhou [HZZ].

\vskip5pt

\subsection{A proper class $\xi$ of triangles in a triangulated category} \ To develop a relative homological algebra in a triangulated category,
A. Beligiannis [B] introduced the notion of a proper class $\xi$ of triangles in a triangulated category $(\T, [1], \tri)$.
A proper class of triangles enjoys pleasant properties: it is closed under the base changes and the cobase changes and it is saturated (see Definition \ref{properclass}).
The class $\tri$ of all the (distinguished) triangles and the class $\tri_0$ of splitting triangles are proper classes of triangles.
There indeed exist ``non-trivial" proper classes $\xi$ of triangles in the sense that $\tri_0 \subsetneqq \xi\subsetneqq \tri$ (see Example \ref{nontrivialxi}).
Morphism  $u$ is a $\xi$-proper monomorphism (respectively, a $\xi$-proper epimorphism)
if there exists a triangle $X \xra{u} Y \xra{} Z \xra{} X[1]$ (respectively, $Z[-1] \xra{} X \xra{u} Y \xra{} Z$) in $\xi$.

\vskip5pt

A new phenomenon of a proper class $\xi$ of triangles which differs from the class $\tri$ of all  triangles is that $\xi$ is not closed under rotations, in general,
i.e.,  if $X \xra{u} Y \xra{v} Z \xra{w} X[1]$ is a triangle in $\xi$, then
the triangles  $Y \xra{v} Z \xra{w} X[1] \xra{-u[1]} Y[1]$ and $Z[1] \xra{-w[1]} X \xra{u} Y\xra{v} Z$ are not necessarily in $\xi$.
For example,  $\tri_0$ is not closed under rotations. In fact, $\xi$ is closed under rotations if and only if $\xi = \tri$ (see Proposition \ref{rotation}).
In view of this fact one can say that triangles in a proper class $\xi$ are more similar to short exact sequences in abelian categories.
In some sense, this phenomenon makes the study of triangles in a proper class more challenging. The group $\Ext_\xi(X, Y)$ with respect to $\xi$ is
only a subgroup of $\Ext^1_\T(X, Y) = \Hom_\T(X, Y[1])$, and hence every result concerning triangles in $\xi$ needs to be carefully reinvestigated.

\vskip5pt

For these reasons,  the Beligiannis - Reiten correspondence should have its present version of triangulated categories with a proper class $\xi$ of triangles.
This is a main purpose of this paper.
If $\xi = \tri$, then  this is known in [HZZ]. However, in general, one should remark that
the corresponding result in an extriangulated category can not include the case for a triangulated category with a proper class $\xi$ of triangles.

\subsection{The homotopy category of a model structure} \ Quillen's {\it homotopy category} of a model structure $(\CoFib, \ \Fib,  \ \Weq)$ on category $\mathcal M$
is by definition the localization $\mathcal M[\Weq^{-1}]$ of $\mathcal M$ with respect to $\Weq$, which is denoted by ${\rm Ho}(\mathcal M)$.
This is an important object of study in algebra and topology. Although $\Ho(\mathcal M)$ is not a triangulated category in general,
the most known triangulated categories arise as the homotopy categories of model structures.

\vskip5pt

Let $\mathcal M$ be a model category. Then $\mathcal M$ has the initial object and the final object (and then there is the notion of  cofibrant objects and cofibrant objects),
and $\mathcal M$ has pull-backs and push-outs.
Let $\mathcal{M}_{cf}$ be the full subcategory of $\mathcal{M}$ of both cofibrant and fibrant objects.
It is proved that the left homotopy relation $\overset{l}\sim$ coincides with
the right homotopy relation $\overset{r}\sim$ \ on $\mathcal{M}_{cf}$, which is denoted by  $\sim$. Then $\sim$ is an equivalence relation ideal of $\mathcal{M}_{cf}$ ([Q1], Lemmas 4, 5, and their duals).
The corresponding quotient category is denoted by $\pi\mathcal{M}_{cf}$, and the composition of the embedding
 $\mathcal{M}_{cf}\hookrightarrow \mathcal{M}$ and the localization functor $\mathcal{M} \longrightarrow {\rm Ho}(\mathcal{M})$ induces an equivalence
\ $\pi\mathcal{M}_{cf} \cong {\rm Ho}(\mathcal{M})$.
See [Q1, Theorem 1'].

\vskip5pt

This important theorem has been named as the Fundamental Theorem of Model Categories. See M. Hovey [Hov1, Theorem 1.2.10]; also J. Gillespie [G1], and W. G. Dwyer - J. Spalinski [DS].

\vskip5pt

Some important categories can not be a model category. For example, the existence of pull-backs and push-outs
can not be guaranteed in triangulated categories. Thus, to apply the Fundamental Theorem of Model Categories in triangulated categories,
one has to weaken the condition that $\mathcal M$ is a model category.
In fact, when Gillespie treats exact model structures,
he has already noticed this problem and he gave a solution for the case ([G1, Section 4]).

\vskip5pt

Dropping the restriction of the uniqueness in the definition of pull-back and push-out, one has the following.

\vskip5pt

\begin{defn} \label{weakpullback} \ Let $\mathcal M$ be a category, and
 \begin{equation}\begin{gathered}\label{wp}
    \xymatrix@R=0.5cm{
      Y' \ar[r]^-{v'}\ar[d]_-{\alpha'} & Z'\ar[d]^-{\alpha} \\
      Y \ar[r]^-{v} & Z
    }
 \end{gathered}
  \end{equation}
a commutative square of morphisms in $\mathcal M$.

\vskip5pt

(1) \ The commutative square (\ref{wp}) is {\it a weak pull-back} of $v$ and $\alpha$, provided that
for any morphisms $f: W\longrightarrow Y$ and $g: W\longrightarrow Z'$ satisfying $v f = \alpha g$, there is a morphism
$h: W\longrightarrow Y'$ such that $\alpha' h = f$ and $v' h = g.$

\vskip5pt

(2) \ The commutative square (\ref{wp}) is {\it a weak push-out} of $v'$ and $\alpha'$, provided that
for any morphisms $f: Y\longrightarrow W$ and $g: Z'\longrightarrow W$ satisfying $f \alpha' = g v'$, there is a morphism
$h: Z\longrightarrow W$ such that $h v = f$ and $h \alpha = g.$
\end{defn}

Carefully checking Quillen's  proof for
[Q1, Theorem 1'], fortunately  one has

\begin{thm}\label{htpthm} \ {\rm (Fundamental Theorem of Model Categories, \ [Q1, Theorem 1'])}

\vskip5pt  Let $\mathcal M$ be a category satisfying the following conditions$:$

\vskip5pt

${\rm (i)}$ \ \  $\mathcal M$ has the initial object and the final object$;$

${\rm (ii)}$ \  $\mathcal M$ has finite coproducts and finite products$;$

${\rm (iii)}$ \ \ There is a model structure $(\CoFib, \ \Fib,  \ \Weq)$ on $\mathcal M;$

${\rm (iv)}$ \ \ For any trivial cofibration $i: A\longrightarrow B$ and any morphism $u: A\longrightarrow C$,  there exists a weak push-out square
$$\xymatrix@R=0.5cm{
A\ar[r]^-u\ar[d]_-{i} & C \ar[d]^-{i'}\\
B\ar[r] & D
}$$
such that $i'$ is also a trivial cofibration.

\vskip5pt

${\rm (v)}$ \ \ For any trivial fibration $p: C\longrightarrow D$ and any morphism $v: B\longrightarrow D$,
there exists a weak pull-back square $$\xymatrix@R=0.5cm{
A\ar[r]\ar[d]_-{p'} & C \ar[d]^-{p}\\
B\ar[r]^-v & D
}$$ such that $p'$ is also a trivial fibration.

\vskip5pt

Then the composition of the embedding
 $\mathcal{M}_{cf}\hookrightarrow \mathcal{M}$ and the localization functor $\mathcal{M} \longrightarrow {\rm Ho}(\mathcal{M})$ induces an equivalence
 \ $\pi\mathcal{M}_{cf} \cong {\rm Ho}(\mathcal{M})$  of categories. \end{thm}

If $\mathcal M$ is a model category, then all the five conditions in Theorem \ref{htpthm} are satisfied.
The proof of Theorem \ref{htpthm} is almost the same as Quillen's original proof.  For convenience, the modifications will be included as Appendix.

\vskip5pt

Since the homotopy cartesian squares induced by base changes and cobase changes in a triangulated category are weak pull-backs and weak push-outs (see Fact \ref{weakpp}),
the conditions (iv) and (v) in Theorem \ref{htpthm} are hopefully satisfied, for some model structures in triangulated categories.
This is indeed the case, for the triangulated $\xi$-model structures ([Y]),
i.e., the triangulated category version of exact model structures (see Theorem \ref{abelianhtpcat}),
and for the $(\xi, \omega)$-model structures (see Lemma \ref{homtopycat}, or Theorem \ref{mainthm}).

\vskip5pt

\subsection {Cotorsion pairs with respect to a proper class of triangles}

Let $(\T, [1], \tri)$ be a triangulated category, and $\xi$ a proper class of triangles. For objects $X$ and $Y$,
using the fact that $\xi$ is closed under base changes and cobase change,
one can define the Baer sum and then an abelian group $\Ext_{\xi}(X, Y)$ (see [B]). It is a subgroup of $\Ext^1_{\T}(X, Y) = \Hom_\T(X, Y[1])$.

\vskip5pt

\begin{thm} \label{intoextlifting} {\rm (The $\xi$-version of Extension-Lifting Lemma. See Theorem \ref{extlifting})} \ Let  $X$ and $Y$ be objects of $\T$.
Then $\Ext_{\xi} (X,Y) = 0$ if and only if for any commutative diagram
\[\xymatrix@R=0.5cm{
      & A \ar[r]^-{i} \ar[d]_-{\alpha} & B \ar[r]^-{d}\ar[d]_-{\beta} & X \ar[r]^-{u} & A[1] \\
    Y \ar[r]^-{c} & C \ar[r]^-{p} & D \ar[r]^-{v} & Y[1] &
    }\]
where the two rows are triangles in $\xi$, there is a lifting $\lambda: B \longrightarrow C$ such that $\alpha = \lambda i$ and $\beta = p \lambda$.
\end{thm}

As expected, one has the following exact sequences with $\Ext_{\xi}.$

\begin{thm}\label{introexactseq}  \ {\rm (See Theorem \ref{exactseq})} \   For any triangle $X_1 \xra{f} X \xra{g} X_2 \xra{h} X_1[1]$ in $\xi$ and object $Y \in \T$,
there are exact sequences of abelian groups $($where $(X_2, Y) = \Hom_\mathcal T(X_2, Y))$
$$(X_2, Y) \xra{(g, Y)} (X, Y) \xra{(f, Y)} (X_1, Y) \xra{\widetilde{h}} \Ext_{\xi}(X_2, Y) \xra{\widetilde{g}} \Ext_{\xi}(X, Y) \xra{\widetilde{f}} \Ext_{\xi}(X_1, Y)$$
and
$$(Y, X_1) \xra{(Y, f)} (Y, X) \xra{(Y, g)} (Y, X_2) \xra{\overline{{h}}} \Ext_{\xi}(Y, X_1) \xra{\overline{{f}}} \Ext_{\xi}(Y, X) \xra{\overline{{g}}} \Ext_{\xi}(Y, X_2).$$
\end{thm}

\vskip5pt

Note that if $\T$ has enough $\xi$-projective objects, Beligiannis [B, Corollary 4.12] gives a similar result:
to  claim the sameness of Theorem \ref{introexactseq} and [B, Corollary 4.12], as in an abelian category, one needs the Yoneda extensions and pull-backs and push-outs.
Since in a triangulated category the existence of pull-backs and push-outs can not be guaranteed, we need to use the present form of Theorem \ref{introexactseq}.

\vskip5pt

Then one can define a (hereditary, complete) cotorsion pair in $\T$ with respect to $\xi$ (see Definition \ref{cpt}).
For any proper class $\xi$ of triangles, $(\mathcal P(\xi), \T)$ is a cotorsion pair with respect to $\xi$, where $\mathcal P(\xi)$ is
the class of $\xi$-projective objects;
and it is complete if and only if $\mathcal T$ has enough $\xi$-projective objects, and in this case $(\mathcal P(\xi), \T)$ is hereditary. See Remark \ref{xiprojxictp}.

\vskip5pt

If $\xi = \tri$, then a hereditary complete cotorsion pair with respect to $\xi$ is precisely a co-$t$-structure,
which has been introduced by M. V. Bondarko [Bon] and D. Pauksztello [P];
but if $\xi\ne \tri$, then the two notion is completely different.
Also, If $\xi = \tri$, then a $t$-structure $(\X, \Y)$, which has been introduced by A. A. Beilinson, J. Bernstein and P. Deligne \cite{BBD},
gives rise to a complete cotorsion pair $(\X, \Y[-2])$ with respect to $\xi$ (but in general not hereditary).

\begin{thm} \label{introheredity} \ {\rm (See Theorem \ref{heredity})} \  A complete cotorsion pair $(\X, \Y)$ in $\T$ with respect to $\xi$
is hereditary if and only if $\mathcal X$ is closed under the hokernel of $\xi$-proper epimorphisms$;$ and if and only if
$\mathcal Y$ is closed under the hocokernel of $\xi$-proper monomorphisms.
\end{thm}

\subsection{The $(\xi, \omega)$-model structures on triangualted categories}

Let $(\T, [1], \tri)$ be a triangulated category,
and $\xi$ a proper class of triangles.  Suppose that $\mathcal{X}$ and $\mathcal{Y}$ are additive full subcategories of $\T$ which are closed under isomorphisms and direct summands, and $\omega:=\mathcal{X}\cap \mathcal{Y}$. Define three classes of morphisms in $\T$ as follows.

\begin{equation} \begin{aligned}\label{Introconstruction}\CoFib_\omega^\xi& = \{\ \xi\text{-proper monomorphism} \ f \ | \ \hoc f\in \mathcal X\}\\
\Fib_\omega^\xi & = \{\text{morphism} \ f \ | \ \Hom_{\mathcal T}(W, f) \ \mbox{is surjective}, \ \forall \ W\in \omega\}\\
\Weq_\omega^\xi & = \{\ A\stackrel f \longrightarrow B \ | \ \mbox{there is a} \ \xi\mbox{-proper epimorphism} \  A\oplus W\stackrel {(f, t)} \longrightarrow B
\\ & \ \ \ \ \ \ \ \ \ \ \ \ \ \ \ \ \ \ \ \ \ \mbox{such that} \ W\in \omega \ \ \mbox{and} \ \ \hok (f, t)\in \mathcal{Y}\}
\end{aligned}\end{equation}

\vskip5pt

\begin{thm}\label{mainthm} \ {\rm (See Theorems \ref{if} and \ref{onlyif})} \ Let $(\T, [1], \tri)$ be a triangulated category,
and $\xi$ a proper class of triangles.  Suppose that $\mathcal{X}$ and $\mathcal{Y}$ are additive full subcategories of $\T$ closed under isomorphisms and direct summands, and $\omega=\mathcal{X}\cap \mathcal{Y}$.
Then $({\rm CoFib}_{\omega}^\xi, \ {\rm Fib}_{\omega}^\xi, \ {\rm Weq}_{\omega}^\xi)$  as defined in {\rm(\ref{Introconstruction})} is a model structure on $\T$ if and only if
$(\mathcal{X}, \mathcal{Y})$ is a hereditary complete cotorsion pair in $(\T, [1], \tri)$ with respect to $\xi$,
and $\omega$ is contravariantly finite in $\T$.

\vskip5pt

If this is the case, then
\begin{align*}\TCoFib_\omega^\xi & = \{\mbox{splitting monomorphism with hocokernel in} \ \omega\} \\
\TFib_\omega^\xi & = \{\xi\mbox{-proper epimorphisms with hokernel in} \ \mathcal Y\};\end{align*}
the class of cofibrant objects $($respectively, fibrant objects, trivial objects$)$ is $\mathcal X$  $($respectively,  $\mathcal T,$ $\mathcal Y);$ and the homotopy category is equivalent to additive quotient $\mathcal X/\omega$.\end{thm}

Theorem \ref{mainthm} is discovered by A. Beligiannis and I. Reiten for abelian categories (see [BR, VIII, Theorem 4.2]).
It is generalized to weakly idempotent complete exact category in [CLZ], and to weakly idempotent complete extriangulated categories by J. S. Hu, D. D. Zhang and P. Y. Zhou [HZZ].
Note that extriangulated categories have been introduced by H. Nakaoka and Y. Palu [NP], as a common generalization of exact categories and triangulated categories.
Thus, if $\xi = \tri$, then  Theorem \ref{mainthm} is known in [HZZ].
However, since there indeed exists a ``non-trivial" proper class $\xi$ of triangles in the sense that $\tri_0 \subsetneqq \xi\subsetneqq \tri$ (see Example \ref{nontrivialxi})
and since a proper class of triangles has quite different properties from the class $\tri$ of all (distinguished) triangles (for example,
a proper class of triangles is not closed under rotations), Theorem \ref{mainthm} should have its present version of triangulated categories with a proper class of triangles.
The model structure given in Theorem \ref{mainthm} will be called the {\it $(\xi, \omega)$-model structure}.

\subsection{The Beligiannis - Reiten correspondence}
A model structure on triangulated category $(\T, [1], \tri)$ with a proper class $\xi$ of triangles is {\it weakly $\xi$-projective} if
cofibrations are precisely the $\xi$-proper monomorphisms with cofibrant hocokernel,
the trivial fibrations are precisely the $\xi$-proper monomorphisms with trivially fibrant hokernel, and each object is fibrant.
See Proposition \ref{weaklyproj} for equivalent characterizations of a weakly $\xi$-projective model structure
on a triangulated category with a proper class $\xi$ of triangles. As in weakly idempotent complete exact categories (\cite[VIII, 4.6]{BR}; [CLZ, Theorem 1.3]),
the $(\xi, \omega)$-model structures on triangulated category $\T$ with proper class $\xi$ are precisely the weakly $\xi$-projective model structures.

\begin{thm} \label{Introcorrespondence} \ {\rm (See Theorem \ref{brcorrespondence})} \ Let $(\T, [1], \tri)$ be a triangulated category, and $\xi$ a proper class of triangles.
Let $\Omega$ denote the class of hereditary complete cotorsion pairs $(\X, \Y)$ in $\T$ with respect to $\xi$, such that
$\omega=\X\cap \Y$ is contravariantly finite in $\T$, and $\Gamma$ the class of weakly $\xi$-projective model structures on $\T$.
Then the map $$\Phi: (\X, \Y)\mapsto (\CoFib_{\omega}^\xi, \ \Fib_{\omega}^\xi, \ \Weq_{\omega}^\xi)$$
as defined in ${\rm(\ref{Introconstruction})},$ gives a bijection between $\Omega$ and $\Gamma$, with inverse given by  $$\Psi: (\CoFib, \Fib, \Weq)\mapsto (\mathcal{C}, \rm T\mathcal{F})$$
where $\mathcal C$ and $\rm T\mathcal{F}$ are respectively the class of cofibrant objects and the class of trivially fibrant objects, of model structure $(\CoFib, \Fib, \Weq)$.
\end{thm}

\subsection{\bf When is the $(\xi, \omega)$-model structure $\xi$-triangulated?}

M. Hovey [Hov2] has introduced {\it abelian model structure} on abelian category, and this has been extended as {\it exact model structure} on exact category by J. Gillespie [G1].
For a triangulated category and a proper class $\xi$ of triangles,  X. Y. Yang [Y] has called the corresponding version as {\it a $\xi$-triangulated model structure}.
See Definition \ref {defnxitriangulatedmds}. The following Hovey correspondence is discovered in [Hov2] for abelian model structures, and
is generalized to exact model structures in [G1]; and further generalized to extriangulated categories
in [NP] (although  extriangulated category is a generalization of exact category and triangulated category,
it can not include the case of a proper class $\xi$ of triangles.)

\begin{thm} \label{introabelianhtpcat} {\rm ([Y, Theorem A]; see Theorem \ref{abelianhtpcat})} \ Let $(\mathcal T, [1], \tri)$ be a triangulated category, and $\xi$ a proper class of triangles.
Then there is a one-one correspondence between $\xi$-triangulated model structures and the {\rm Hovey} triples in $\mathcal T$, given by
$$({\rm CoFib}, \ {\rm Fib}, \ {\rm Weq})\mapsto (\mathcal{C}, \ \mathcal{F}, \ \mathcal W)$$
where \ $\mathcal C  = \{\mbox{cofibrant objects}\}, \ \
\mathcal F  = \{\mbox{fibrant objects}\}, \ \
\mathcal W  = \{\mbox {trivial objects}\}$, with the inverse \ \ $(\mathcal{C}, \ \mathcal{F}, \ \mathcal W) \mapsto ({\rm CoFib}, \ {\rm Fib}, \ {\rm Weq}),$ where
\begin{align*} &{\rm CoFib} = \{\mbox{$\xi$-proper monomorphism with hocokernel in} \ \mathcal{C}\}, \\ &
{\rm Fib}  = \{\mbox{$\xi$-proper epimorphism with hokernel in} \ \mathcal{F} \}, \\
& {\rm Weq}  = \{pi \ \mid \ i \ \mbox{is a $\xi$-proper monomorphism,} \ \hoc i\in \mathcal{C}\cap \mathcal W, \\ & \ \ \ \ \ \ \ \ \ \ \ \ \ \ \ \ \  p \ \mbox{is a $\xi$-proper epimorphism,} \ \hok p\in \mathcal{F}\cap \mathcal W\}.\end{align*}
If this is the case, then the homotopy category $\Ho(\mathcal{T})$ is equivalent to additive quotient \ $(\mathcal{C} \cap \mathcal{F})/(\mathcal{C} \cap \mathcal{F} \cap \mathcal W)$.
\end{thm}

Theorem \ref{abelianhtpcat} is  proved in [Y], except for the final  part for $\Ho(\mathcal{T})$; and this will be done in Subsection 8.1, as an application of Theorem \ref{htpthm}.

\vskip5pt

\begin{rem} \label{hotri} \ If $\xi = \tri$, then $\Ho(\mathcal{T})$ is in fact a triangulated category, by H. Nakaoka and Y. Palu [NP, Theorem 6.20].
We conjecture that $\Ho(\mathcal{T})$ is always a triangulated category, for an arbitrary proper class $\xi$ of triangles.
\end{rem}

\vskip5pt

It is natural to ask when the $(\xi, \omega)$-model structure $(\CoFib_{\omega}^\xi, \Fib_{\omega}^\xi, \Weq_{\omega}^\xi)$  is $\xi$-triangulated.
It turns out that, the $(\xi, \omega)$-model structure $(\CoFib_{\omega}^\xi, \ \Fib_{\omega}^\xi, \ \Weq_{\omega}^\xi)$
is $\xi$-triangulated if and only if $\mathcal T$ has enough $\xi$-projective objects and $\X\cap \Y = \mathcal P(\xi)$,
the class of $\xi$-projective objects. See Proposition \ref{thesame}, and hence
the intersection of the class of $\xi$-triangulated model structures and the class of weakly $\xi$-projective model structures
is the class of $\xi$-projective model structures.

\vskip5pt

\subsection{The organization} Section 2 recalls  necessary preliminaries on triangulated categories, properties of a proper class of triangles, and model structures.
In particular, an example of a ``non-trivial" proper class of triangles is given in Example \ref{nontrivialxi}.

\vskip5pt

In Section 3 group  $\Ext_{\xi}(X, Y)$ and hereditary complete cotorsion pairs in $\T$ with respect to $\xi$
are discussed, and  Theorems \ref {extlifting}, \ref{exactseq}, \ref{heredity} are proved.

\vskip5pt

Section 4 is devoted to the proof of the ``if" part of Theorem \ref{mainthm};
and Section 5 is to prove the ``only if" part of Theorem \ref{mainthm}, including the description of the homotopy category.
In Section 6, weakly $\xi$-projective model structures are characterized, and Theorem \ref{Introcorrespondence} is proved.
In Section 7, the dual version of Theorems \ref{mainthm} and \ref{Introcorrespondence} is stated.

\vskip5pt

In Section 8, the homotopy category of a $\xi$-triangulated model structure is clarified, as an application of Theorem \ref{htpthm};
and the model structures which are both $\xi$-triangulated and weakly $\xi$-projective are characterized.
The proof of Theorem \ref{htpthm} is almost the same as Quillen's original proof. For the convenience, the places of the modifications are
included as Appendix.

\section{\bf Preliminaries}

\subsection{Triangulated categories.} \ Let $\mathcal T$ be an additive category, and $[1]: \T \xra{} \T$ an autoequivalence of category.
A quasi-inverse of $[1]$ is denoted by $[-1]$. {\it A pre-triangulated category} is a triple $(\T, [1], \tri)$, where
$\tri$ is a class of  sequences $X \xra{} Y \xra{} Z \xra{} X[1]$ of morphisms in $\T$, satisfying axioms (TR1), (TR2), and (TR3).
An element $X \xra{} Y \xra{} Z \xra{} X[1]$ in $\tri$ is called {\it a triangle}. A pre-triangulated category $(\T, [1], \tri)$ is {\it triangulated}, if  the octahedral axiom (TR4) holds.
For details we refer to [V], [N2], or [Hap].

\vskip5pt

For a pre-triangulated category $(\T, [1], \tri)$, any morphism $u: X \xra{} Y$  in $\T$ can be embedded into a triangle $X \xra{u} Y \xra{v} Z \xra{w} X[1]$ in the unique way, up to an isomorphism of triangles.
The morphism $v: Y \longrightarrow Z$, or simply, the object $Z$, is called {\it the homotopy cokernel} of $u: X \xra{} Y$, and is denoted by $\hoc u$, i.e.,
$\hoc u = Z$.

\vskip5pt

Also, any morphism $u: X \xra{} Y$  in $\T$ can be embedded into a triangle $W \xra{f} X \xra{u} Y \xra{g} W[1]$ in the unique way, up to an isomorphism of triangles.
The morphism $f: W \longrightarrow X$, or simply, the object $W$, is called {\it the homotopy kernel} of $u: X \xra{} Y$, and is denoted by $\hok u$, i.e.,
$\hoc u = W$.

\vskip5pt

We need various equivalent forms of the octahedral axiom (TR4). The following well-known fact can be found in [B, Proposition 2.1] or [M].

\begin{lem}\label{bcandcbc} \ Let $(\T, [1], \tri)$ be a pre-triangulated category. Then the following are equivalent$:$

  \vskip5pt

  {\rm (1) \ {\bf (Tr4). \ (The octahedral axiom)}} \ \  Given triangles in the first and the second rows, and in second column, as in the diagram below,
  there exists the following commutative diagram such that the third column is also a triangle.
  \begin{equation}
  \begin{gathered}
    \label{tr4}
    \xymatrix@=0.6cm{
    X\ar[r]^-{u}\ar@{=}[d] & Y\ar[r]^-{u'}\ar[d]_-{v} & Z' \ar[r]^-{u''}\ar[d]^-{f} & X[1]\ar@{=}[d] \\
    X\ar[r]^-{vu} & Z\ar[r]^-{w}\ar[d]_-{v'} & Y'\ar[r]^-{w'}\ar[d]^-{g} & X[1]\ar[d]^-{u[1]} \\
      & X'\ar@{=}[r]\ar[d]_-{v''} & X'\ar[d]\ar[r]^-{v''} & Y[1] \\
      & Y[1]\ar[r]^-{u'[1]} & Z'[1] & \\
    }
   \end{gathered}
  \end{equation}

  \vskip5pt

  {\rm (2) \ ${\bf (Tr4').}$} \ \  For any morphisms $u: X \xra{} Y$ and $v: Y \xra{} Z$, there exists a commutative diagram {\rm (\ref{tr4})}
  such that the first and the second rows, and the second and the third columns are triangles.

  \vskip5pt

 {\rm (3) \ {\bf Base change.}} \ \ For any triangle $X \xra{u} Y \xra{v} Z \xra{w} X[1]$ and any morphism $\alpha: Z' \xra{} Z$, there exists a commutative diagram
  \begin{equation}
  \begin{gathered}
    \label{bc}
    \xymatrix@=0.6cm{
      & X' \ar@{=}[r]\ar[d]_-{\beta'} & X' \ar[d] & \\
    X \ar[r]^-{u'}\ar@{=}[d] & Y' \ar[d]_-{\alpha'}\ar[r]^-{v'} & Z' \ar[r]\ar[d]_-{\alpha} & X[1]\ar@{=}[d] \\
    X \ar[r]^-{u} & Y \ar[r]^-{v}\ar[d] & Z \ar[r]^-{w}\ar[d]_-{\gamma} & X[1]\ar[d]^-{u'[1]} \\
      & X'[1]\ar@{=}[r] & X'[1]\ar[r]^-{-\beta'[1]} &  Y'[1]\\
    }
  \end{gathered}
  \end{equation}
 such that the second row, and the second and the third columns are triangles.

  \vskip5pt

  {\rm (4) \ {\bf Cobase change.}} \ \ For any triangle $X \xra{u} Y \xra{v} Z \xra{w} X[1]$ and any morphism $\beta: X \xra{} X'$, there exists a commutative diagram
  \begin{equation}
    \begin{gathered}\label{cbc}
    \xymatrix@=0.6cm{
      & Z'[-1] \ar@{=}[r]\ar[d]^-{-\gamma[-1]} & Z'[-1] \ar[d] & \\
    Z[-1] \ar[r]^-{-w[-1]}\ar@{=}[d] & X \ar[d]^-{\beta}\ar[r]^-{u} & Y \ar[r]^-{v}\ar[d]^-{\beta'} & Z\ar@{=}[d] \\
    Z[-1] \ar[r] & X' \ar[r]^-{u'}\ar[d] & Y' \ar[r]^-{v'}\ar[d]^-{\alpha'} & Z\ar[d]^-{-w} \\
      & Z' \ar@{=}[r] &  Z'\ar[r]^-{\gamma} &  X[1]\\
    }
  \end{gathered}
  \end{equation}
  such that the third row, and the second and the third columns are triangles.
\end{lem}

The diagram (\ref{bc}) is said to be {\it the base change diagram} of triangle $X \xra{u} Y \xra{v} Z \xra{w} X[1]$ along morphism $\alpha : Z' \xra{} Z$.
Transposing the rows and the columns, it is also the base change diagram of triangle
$X' \xra{\beta} Z' \xra{\alpha} Z \xra{\gamma} X'[1]$ along morphism $v: Y \xra{} Z$.

\vskip5pt

Dually, the diagram (\ref{cbc}) is said to be {\it the cobase change diagram } of triangle $X \xra{u} Y \xra{v} Z \xra{w} X[1]$ along morphism $\beta : X \xra{} X'$.
Transposing the rows and the columns,  it is also the cobase change diagram of triangle $X \xra{\beta} X' \xra{\alpha} Z' \xra{\gamma} X[1]$ along morphism $u: X \xra{} Y$.

\vskip10pt

A commutative square in a pre-triangulated category $\T$
  \[
    \xymatrix@R=0.5cm{
      Y \ar[r]^-{v}\ar[d]_-{g} & Z\ar[d]^-{h} \\
      Y' \ar[r]^-{v'} & Z'
    }
  \]
is {\it homotopy cartesian} if there is a triangle of the form
  $$Y \xra{\left(\begin{smallmatrix} g\\ -v \end{smallmatrix}\right)} Y' \oplus Z \xra{(v', h)} Z' \xra{\delta} Y[1].$$

For the proof of the following result we refer to A. Neeman \cite{N1} and \cite{N2}, and A. Hubery \cite{Hub}.

\begin{lem}\label{addaxiom}
  \ Let $(\T, [1], \tri)$ be a pre-triangulated category. Then the following are equivalent$:$

  \vskip5pt

  $(1)$  \ \ The octahedral axiom.

  \vskip5pt

  $(2)$ \ \ Given a commutative diagram with two rows in $\tri$
  \[
    \xymatrix@=0.6cm{
    X \ar[r]^-{u}\ar[d]_-{f}& Y \ar[r]^-{v}\ar[d]_-{g}& Z\ar[r]^-{w} & X[1]\ar[d]^-{f[1]} \\
    X' \ar[r]^-{u'}& Y'\ar[r]^-{v'} & Z'\ar[r]^-{w'} & X'[1] \\
  }
  \]
  there exists a morphism $h: Z \xra{} Z'$ such that the mapping cone
  $$Y \oplus X' \xra{\left(\begin{smallmatrix} -v & 0 \\ g & u' \end{smallmatrix}\right)} Z \oplus Y' \xra{\left(\begin{smallmatrix} -w & 0 \\ h & v' \end{smallmatrix}\right)} X[1] \oplus Z' \xra{\left(\begin{smallmatrix} -u[1] & 0 \\ f[1] & w' \end{smallmatrix}\right)} Y[1] \oplus X'[1]$$
  is a triangle.

  \vskip5pt

  $(3)$ \ \ Given a commutative diagram with two rows in $\tri$
  \[
    \xymatrix@=0.6cm{
    X \ar[r]^-{u}\ar@{=}[d]& Y \ar[r]^-{v}\ar[d]_-{g}& Z\ar[r]^-{w} & X[1]\ar@{=}[d] \\
    X \ar[r]^-{u'}& Y'\ar[r]^-{v'} & Z'\ar[r]^-{w'} & X[1] \\
  }
  \]
  there exists a morphism $h: Z \xra{} Z'$ such that the middle square is homotopy cartesian, i.e.,
  $$Y \xra{\left(\begin{smallmatrix} g\\ -v \end{smallmatrix}\right)} Y' \oplus Z \xra{(v', h)} Z' \xra{u[1]w'} Y[1]$$
  is a triangle.

  \vskip5pt

  $(4)$ \ \ Given a triangle
  $$Y \xra{\left(\begin{smallmatrix} g\\ -v \end{smallmatrix}\right)} Y' \oplus Z \xra{(v', h)} Z' \xra{\delta} Y[1]$$
  in $\tri$, there exists the following commutative diagram
  \[
    \xymatrix@=0.6cm{
    X\ar[r]^-{u}\ar@{=}[d] & Y \ar[r]^-{v}\ar[d]^-{g} & Z \ar[r]^-{w}\ar[d]^-{h} & X[1]\ar@{=}[d] \\
    X\ar[r]^-{u'} & Y' \ar[r]^-{v'}\ar[d]^-{g'} & Z' \ar[r]^-{w'}\ar[d]^-{h'} & X[1] \ar[d]^-{u[1]} \\
      & W \ar@{=}[r]\ar[d]^-{f} & W \ar[d]^-{f'}\ar[r]^-{f} & Y[1] \\
      & Y[1]\ar[r]^-{v[1]} & Z[1] & \\
    }
  \]
such that the first and the second rows, and the second and the third columns are triangles, and that $u[1]w' = fh' = \delta$.
\vskip5pt

  $(5)$ \ \  For any triangle $X \xra{u} Y \xra{v} Z \xra{w} X[1]$ and any morphism $\alpha: Z' \xra{} Z$, there exists a commutative diagram
    $$\xymatrix@=0.6cm{
      & X' \ar@{=}[r]\ar[d]_-{\beta'} & X' \ar[d] & \\
    X \ar[r]^-{u'}\ar@{=}[d] & Y' \ar[d]_-{\alpha'}\ar[r]^-{v'} & Z' \ar[r]\ar[d]_-{\alpha} & X[1]\ar@{=}[d] \\
    X \ar[r]^-{u} & Y \ar[r]^-{v}\ar[d] & Z \ar[r]^-{w}\ar[d]_-{\gamma} & X[1]\ar[d]^-{u'[1]} \\
      & X'[1]\ar@{=}[r] & X'[1]\ar[r]^-{-\beta'[1]} &  Y'[1]\\
    }$$
  such that the second row, and the second and the third columns are triangles, and the middle square is homotopy cartesian, i.e.,
  $$Y' \xra{\left(\begin{smallmatrix} \alpha' \\ -v' \end{smallmatrix}\right)} Y \oplus Z' \xra{(v, \alpha)} Z \xra{u'[1]w} Y'[1]$$
  is a triangle.

  \vskip5pt

  $(6)$ \ \ For any triangle $X \xra{u} Y \xra{v} Z \xra{w} X[1]$ and any morphism $\beta: X \xra{} X'$, there exists a commutative diagram
    $$\xymatrix@=0.6cm{
      & Z'[-1] \ar@{=}[r]\ar[d]^-{-\gamma[-1]} & Z'[-1] \ar[d] & \\
    Z[-1] \ar[r]^-{-w[-1]}\ar@{=}[d] & X \ar[d]^-{\beta}\ar[r]^-{u} & Y \ar[r]^-{v}\ar[d]^-{\beta'} & Z\ar@{=}[d] \\
    Z[-1] \ar[r] & X' \ar[r]^-{u'}\ar[d] & Y' \ar[r]^-{v'}\ar[d]^-{\alpha'} & Z\ar[d]^-{-w} \\
      & Z' \ar@{=}[r] &  Z'\ar[r]^-{\gamma} &  X[1]\\
    }$$
  such that the third row, and the second and the third columns are triangles, and the middle square is homotopy cartesian, i.e.,
  $$X \xra{\left(\begin{smallmatrix} \beta \\ -u \end{smallmatrix}\right)} X' \oplus Y \xra{(u', \beta ')} Y' \xra{-wv'} X[1]$$
  is a triangle.\end{lem}

\begin{fact} \label{weakpp} \ Let $\mathcal T$ be a triangulated category. Then

 \vskip5pt

$(1)$ \ For any morphisms $v: Y \longrightarrow Z$ and $\alpha: Z'\longrightarrow Z$, there exists a weak pull-back of $v$ and $\alpha$.

\vskip5pt

$(2)$ \ For any morphisms $v': Y' \longrightarrow Z'$ and $\alpha': Y'\longrightarrow Y$, there exists a weak push-out of $v'$ and $\alpha'$.
\end{fact}
\begin{proof} We only justify (1). The assertion (2) can be similarly proved. The middle square in Lemma \ref{addaxiom}(5) is homotopy cartesian, and hence it is
a weak pull-back of $v$ and $\alpha$. In fact,
  $$Y' \xra{\left(\begin{smallmatrix} \alpha' \\ -v' \end{smallmatrix}\right)} Y \oplus Z' \xra{(v, \alpha)} Z \xra{u'[1]w} Y'[1]$$
  is a triangle.  For any morphisms $f: W\longrightarrow Y$ and $g: W\longrightarrow Z'$ satisfying $v f = \alpha g$, applying $\Hom_\mathcal T(W, -)$ to the triangle one has an exact sequence
  $$\xymatrix@C=0.7cm{\Hom_\mathcal T(W, Y') \ar[rr]^-{(W,  \binom{\alpha'}{-v'})} && \Hom_\mathcal T(W, Y \oplus Z') \ar[rr]^-{(W, (v, \alpha))} && \Hom_\mathcal T(W, Z)}.$$
Then the morphism $\binom{f}{-g}: W\longrightarrow Y \oplus Z'\in \Ker \Hom_\mathcal T(W, (v, \alpha))$, and hence one gets a morphism
$h: W\longrightarrow Y'$ such that $\alpha' h = f$ and $v' h = g.$
\end{proof}

\vskip5pt

Using cohomological functors $\Hom(Z, -)$ and $\Hom(-, X')$ one sees the following fact.

\begin{fact}\label{2factor} \ Given a morphism of triangles in a pre-triangulated category
    $$\xymatrix@R=0.6cm{
        X \ar[r]^-{x}\ar[d]_-{f} & Y \ar[r]^-{y}\ar[d]_-{g} & Z \ar[r]^-{z}\ar[d]_-{h} & X[1]\ar[d]^-{f[1]} \\
        X' \ar[r]^-{x'} & Y' \ar[r]^-{y'} & Z'\ar[r]^-{z'} & X'[1]
    }$$
then $f$ factors through $x$ if and only if $h$ factors through $y'$.
\end{fact}

\subsection{Splitting triangles.}

\begin{deffact} \label{split} \ Let $(\T, [1], \tri)$ be a pre-triangulated category. A triangle $X \stackrel{u}{\longrightarrow} Y
\stackrel{v}{\longrightarrow}  Z\stackrel{w}{\longrightarrow}
TX$ is called a splitting triangle, if the following equivalent conditions are satisfied.

\vskip5pt

${\rm (i)}$ \  \ $w=0.$

\vskip5pt

${\rm (ii)}$ \ \ $u$ is a splitting monomorphism.

\vskip5pt

${\rm (iii)}$ \  \ $v$ is a splitting epimorphism.

\vskip5pt

${\rm (iv)}$ \ \ There are morphisms $s: \ Y\longrightarrow X$ and $t: \ Z\longrightarrow Y$
such that
$$su = {\rm Id}_X, \ \ vt = {\rm Id}_Z,  \ \ us+tv = {\rm
Id}_Y.$$

${\rm (v)}$ \ \ There are morphisms $s: \ Y\longrightarrow X$ and $t: \ Z\longrightarrow Y$ such that
$$\xymatrix{X\ar[r]^u\ar@{=}[d] & Y \ar[r]^-v \ar@{-->}[d]_-{\binom{s}{v}}& Z \ar[r]^-w\ar@{=}[d] &  X[1]\ar@{=}[d]\\
X \ar[r]^-{\left( \begin{smallmatrix}1\\
0\end{smallmatrix}\right)}& X\oplus Z\ar[r]^-{(0, 1)} & Z\ar[r]^-0 &  X[1]
}$$
\vskip5pt
\noindent is a morphism of triangles, where the inverse of $\binom{s}{v}$ is $(u, t): X\oplus Z\longrightarrow Y$, i.e.,
$$su = {\rm Id}_X, \ \ vt = {\rm Id}_Z,  \ \ st = 0, \ \ us+tv = {\rm
Id}_Y.$$
Moreover, for any morphism $s': \ Y\longrightarrow X$ with  $s'u = {\rm Id}_X$, morphism $\binom{s'}{v}: Y\longrightarrow X\oplus Z$ is an isomorphism,
thus there is a unique morphism  $t': \ Z\longrightarrow Y$ such that $(u, t'): X\oplus Z\longrightarrow Y$ is the inverse of $\binom{s'}{v}$.

\vskip5pt

${\rm (vi)}$ \ \ $u$ is a splitting monomorphism, and $v$ is the cokernel of $u$.

\vskip5pt

${\rm (vii)}$ \  \ $v$ is a splitting epimorphism, and $u$ is the kernel of $v$.

\vskip5pt

If this is the case, then $\hoc u = \Coker u$ and $\hok v = \Ker v$.
\end{deffact}
\begin{proof}  \ This is well-known. For example, the equivalences of (i), (ii) and (iii) have been explicitly stated in D. Happel [Hap, p.7].
The implications (v) $\Longrightarrow$ (i), (v)$\Longrightarrow$ (iv) and (iv)$\Longrightarrow$ (i) are trivial.
For the implication (i) $\Longrightarrow$ (v), by the commutative square at the right hand side one gets morphism $\binom{s}{v}$ of triangles which must be an isomorphism.
Moreover, for any morphism $s': \ Y\longrightarrow X$ with  $s'u = {\rm Id}_X$, then by the commutative square at the left hand side one gets the desired conclusion.
The equivalences of (v), (vi) and (vii) are also clear.
\end{proof}

\subsection{Proper class of triangles.}

\begin{defn}\label{closed} {\rm ([B])} \ \ Let $(\T, [1], \tri)$ be a triangulated category, and $\xi$ a class of triangles, i.e., $\xi \subseteq\tri.$

\vskip5pt

$(1)$ \ \ The class $\xi$ is closed under base changes,
if for any triangle $X \xra{u} Y \xra{v} Z \xra{w} X[1]$  in $\xi$
and any morphism $\alpha: Z' \xra{} Z$,  as in the base change diagram (\ref{bc}),
the triangle $X \xra{u'} Y' \xra{v'} Z' \xra{w'} X[1]$ is in $\xi$.

\vskip5pt

(2) \ \ Dually, $\xi$  is closed under cobase changes,
if for any triangle $X \xra{u} Y \xra{v} Z \xra{w} X[1]$ in $\xi$
and any morphism $\beta : X \xra{} X'$,  as in the cobase change diagram (\ref{cbc}),
the triangle $X' \xra{u'} Y' \xra{v'} Z \xra{w'} X'[1]$ is in $\xi$.

\vskip5pt

(3) \ \ The class $\xi$  is closed under suspensions,
if for any triangle $X \xra{u} Y \xra{v} Z \xra{w} X[1]$ in $\xi$
and any $i \in \mathbb{Z}$,
the triangle $X[i] \xra{(-1)^i u[i]} Y[i] \xra{(-1)^i v[i]} Z[i] \xra{(-1)^i w[i]} X[i+1]$ is in $\xi$.

\vskip5pt

(4) \ \ The class $\xi$ is saturated, if in the situation of the base change diagram (\ref{bc}),
whenever the second horizontal triangle $X \xra{u'} Y' \xra{v'} Z' \xra{} X[1]$
and the third vertical triangle $X' \xra{} Z' \xra{\alpha} Z \xra{\gamma} X'[1]$
are in $\xi$,
then the third horizontal triangle $X \xra{u} Y \xra{v} Z \xra{w} X[1]$ is in $\xi$.
\end{defn}

\begin{defn}\label{properclass} {\rm ([B])} \ \ Let $(\T, [1], \tri)$ be a triangulated category, $\xi\subseteq \tri$, and $u: X \xra{} Y$ a morphism in $\T$.

\vskip5pt

(1) \ The class  $\xi$ is called a proper class of triangles, provided that the following conditions are satisfied$:$

\hskip18pt (i) \ \ \  \ $\xi \supseteq \tri_0,$ where $\tri_0$ is the class of splitting triangles;

\hskip18pt (ii) \ \ \ $\xi$ is closed under isomorphisms of triangles, suspensions and finite coproducts;

\hskip18pt (iii) \ \   $\xi$ is closed under the base changes and the cobase changes;

\hskip18pt (iv) \ \  $\xi$ is saturated.

\vskip5pt

(2) \ The morphism  $u$ is a $\xi$-proper monomorphism
if there exists a triangle $X \xra{u} Y \xra{} Z \xra{} X[1]$ in $\xi$.

\vskip5pt

(3) \ The morphism $u$ is a $\xi$-proper epimorphism
if there exists a triangle $Z[-1] \xra{} X \xra{u} Y \xra{} Z$ in $\xi$.
\end{defn}

\vskip5pt

Clearly, both $\tri$ and $\tri_0$ are proper classes of triangles.
To justify the existence of ``non-trivial" examples of  proper classes $\xi$ of triangles
(i.e., a proper class $\xi$ of triangles such that $\tri_0 \subsetneqq \xi\subsetneqq \tri$), we use the following observation by A. Beligiannis.

\begin{lem} {\rm([B, Example 2.3(3)])} \label{exmproperclass} \ Let $(\T, [1], \tri)$ be a triangulated category, $\mathcal{A}$ an abelian category and $F: \T \xra{} \mathcal{A}$ a cohomological functor. Then $\xi(F)$ is a proper class of triangles, where $\xi(F)$ is the class of all triangles $X \xra{} Y \xra{} Z \xra{} X[1]$ such that
the induced sequence $0 \xra{} F(X[i]) \xra{} F(Y[i]) \xra{}F(Z[i]) \xra{} 0$ is exact for any $i \in \mathbb{Z}$.
  \end{lem}

\begin{exm} \label{nontrivialxi} \ Let $\mathcal{A}$ be an abelian category, $K(\mathcal{A})$ its homotopy category, and $P$ a non-zero projective object in $\mathcal{A}$.
  Denote by $P^{\bullet}$ the stalk complex with $P^{\bullet 0} = P$ (thus $P^{\bullet i} = 0$ for all $i \neq 0$).
  Then the proper class of triangles $\xi = \xi (\Hom_{K(\mathcal{A})}(P^{\bullet}, -))$ is non-trivial, i.e.,
   $\tri_0 \subsetneqq \xi\subsetneqq \tri$.

    \vskip5pt

    In fact, if $X^{\bullet} \xra{f} Y^{\bullet} \xra{} {\rm Cone}(f) \xra{} X^{\bullet}[1]$ is a triangle with $X^{\bullet}$ and $Y^{\bullet}$ acyclic, then so is $Z^{\bullet} = {\rm Cone}(f)$.
Note that a chain map $\alpha: P^{\bullet} \xra{} X^{\bullet}$ is precisely a morphism $\alpha^0: P \xra{} X^0$ such that $d_X^0 \alpha^0 = 0$.
Since by assumption $X^{\bullet}$ is acyclic, $\alpha^0$ factors through $d_X^{-1}$, and thus any chain map $\alpha: P^{\bullet} \xra{} X^{\bullet}$ is a null homotopy.
Thus one has
$$\Hom_{K(\mathcal{A})}(P^{\bullet}, X^{\bullet}[i]) = \Hom_{K(\mathcal{A})}(P^{\bullet}, Y^{\bullet}[i]) = \Hom_{K(\mathcal{A})}(P^{\bullet}, Z^{\bullet}[i]) = 0.$$
Therefore $\xi$ contains all the triangles $X^{\bullet} \xra{f} Y^{\bullet} \xra{} {\rm Cone}(f) \xra{} X^{\bullet}[1]$ where $X^{\bullet}$ and $Y^{\bullet}$ are acyclic.
Such a triangle is not necessarily a splitting triangle in $K(\mathcal A)$. Thus $\tri_0 \subsetneqq \xi$.

    \vskip5pt

On the other hand, for example, consider the triangle
$$X^{\bullet} \xra{0} Y^{\bullet} \xra{\binom{0}{1}} {\rm Cone}(0) = X[1]\oplus Y \xra{(1, 0)} X^{\bullet}[1]$$ where $X^{\bullet}$ is the stalk complex with $X^{\bullet 0} = X^0$ (thus $X^{\bullet i} = 0$ for all $i \neq 0$) such that there is a non-zero morphism $\pi: P \xra{} X^0$, and $Y^{\bullet}$ is stalk complex with $Y^{\bullet 0} = Y^0$.
Then the map $\Hom_{K(\mathcal{A})}(P^{\bullet}, 0) = 0: 0\ne \Hom_{K(\mathcal{A})}(P^{\bullet}, X^{\bullet}) \xra{} \Hom_{K(\mathcal{A})}(P^{\bullet}, Y^{\bullet})$ is not injective.
Thus the triangle $X^{\bullet} \xra{0} Y^{\bullet} \xra{} {\rm Cone}(0) \xra{} X^{\bullet}[1]$ is not in $\xi$ and hence $\xi \subsetneqq \tri$.
  \end{exm}

\vskip5pt

An important phenomenon of a proper class $\xi$ of triangles is that $\xi$ is {\bf not} closed under rotations, in general. That is,
a triangle  $X \xra{u} Y \xra{v} Z \xra{w} X[1]$ in $\xi$
does not imply that the triangle $Y \xra{v} Z \xra{w} X[1] \xra{-u[1]} Y[1]$ is in $\xi$,
and does not imply that the triangle $Z[1] \xra{-w[1]} X \xra{u} Y\xra{v} Z$ is in $\xi$.
For example,  $\tri_0$ is not closed under rotations.
In this way one can say that triangles in a proper class $\xi$ are more similar to short exact sequences in abelian categories.

\vskip5pt

Indeed, one has the following fact.

\begin{prop}\label{rotation} \ Let $(\T, [1], \tri)$ be a triangulated category, and $\xi$ a proper class of triangles. Then $\xi$ is closed under rotations if and only if $\xi = \tri$.
\end{prop}

\begin{proof} \ We only need to prove the ``only if'' part. Since $\xi$ is closed under rotations and $\tri_0 \subseteq \xi$, the triangle $Y[-1] \xra{} 0 \xra{} Y \xra{=} Y$ is in $\xi$ for any object $Y \in \T$. Thus for any morphism $f: X \xra{} Y$, by taking the base change of $Y[-1] \xra{} 0 \xra{} Y \xra{=} Y$ along $f$ one gets the following commutative diagram
  $$\xymatrix@=0.6cm{
      & W \ar@{=}[r]\ar[d] & W \ar[d] & \\
    Y[-1] \ar[r]\ar@{=}[d] & W \ar[d]\ar[r] & X \ar[r]^-{f}\ar[d]^-{f} & Y\ar@{=}[d] \\
    Y[-1] \ar[r] & 0 \ar[r]\ar[d] & Y \ar[r]^-{=}\ar[d] & Y\ar[d] \\
      & W[1] \ar@{=}[r] &  W[1]\ar@{=}[r] &  W[1]\\
    }$$
    and hence the triangle $Y[-1] \xra{} W \xra{} X \xra{f} Y$ is in $\xi$. Since $\xi$ is closed under rotations, any morphism $f$ can be embedded into a triangle $X \xra{f} Y \xra{} W[1] \xra{} X[1]$ in $\xi$. So $\xi = \tri$.
\end{proof}

\vskip5pt

The proofs of the following two results can be found in \cite{Y}.

\begin{lem}\label{carsq} \ {\rm ([Y, Proposition 2.2])}
  \ Let $(\T, [1], \tri)$ be a triangulated category, and $\xi$ a proper class of triangles. Then

  \vskip5pt

  $(1)$ \ As in the base change diagram {\rm (\ref{bc})}, if the triangle $X \xra{u} Y \xra{v} Z \xra{w} X[1]$ is in $\xi$,
then  $Y' \xra{\left(\begin{smallmatrix} \alpha '\\ -v' \end{smallmatrix}\right)} Y \oplus Z' \xra{(v, \alpha)} Z \xra{u'[1]w} Y'[1]$ is a triangle in $\xi$.

\vskip5pt

In particular, if $v$ is a $\xi$-proper epimorphism,
then so is $(v, \alpha)$  for any morphism $\alpha: Z'\longrightarrow Z$.

  \vskip5pt

  $(2)$ \ The notion of saturability is self-dual. That is, $\xi$ is saturated if and only if in the situation of the cobase change diagram {\rm(\ref{cbc})},
the triangles $X' \xra{u'} Y' \xra{v'} Z \xra{w'} X'[1]$ and $X \xra{\beta} X' \xra{\alpha} Z' \xra{\gamma} X[1]$ are in $\xi$,
then the triangle $X \xra{u} Y \xra{v} Z \xra{w} X[1]$ is also in $\xi$.

  \vskip5pt

  $(3)$ \ As in the cobase change diagram {\rm (\ref{cbc})}, if the triangle $X \xra{u} Y \xra{v} Z \xra{w} X[1]$ is in $\xi$,
then $X \xra{\left(\begin{smallmatrix} \beta \\ u \end{smallmatrix}\right)} X' \oplus Y \xra{(u', -\beta ')} Y' \xra{wv'} X[1]$ is a triangle in $\xi$.

\vskip5pt

In particular, if $u$ is a $\xi$-proper monomorphism, then so is
$\left(\begin{smallmatrix}\beta \\ u \end{smallmatrix}\right)$ for any morphism $\beta: X\longrightarrow X'$.
\end{lem}

\begin{lem}\label{xiclose} \ {\rm ([Y, Proposition 2.3])}
  \ Let $(\T, [1], \tri)$ be a triangulated category, and $\xi$ a proper class of triangles. Given a commutative diagram
  \[
    \xymatrix@=0.6cm{
      & Z' \ar@{=}[r]\ar[d]^-{\alpha '} & Z' \ar[d]^-{\alpha} & \\
    Z[-1] \ar[r]^-{-w'[-1]}\ar@{=}[d] & X' \ar[d]^-{\beta '}\ar[r]^-{u'} & Y' \ar[r]^-{v'}\ar[d]^-{\beta} & Z\ar@{=}[d] \\
    Z[-1] \ar[r]^-{-w[-1]} & X \ar[r]^-{u}\ar[d]^-{\gamma'} & Y \ar[r]^-{v}\ar[d]^-{\gamma} & Z\ar[d]^-{-w'} \\
      & Z'[1] \ar@{=}[r] &  Z'[1]\ar[r]^-{-\alpha'[1]} &  X'[1]\\
    }
  \]
where the two rows and the two columns are triangles.

  \vskip5pt

$(1)$ \  If the triangle $X \xra{u} Y \xra{v} Z \xra{w} X[1]$ and  the triangle in the third column are in $\xi$,
then the triangle $X' \xra{u'} Y' \xra{v'} Z \xra{w'} X'[1]$ is also in $\xi$.

  \vskip5pt

$(2)$ \ If the triangles $Z' \xra{\alpha'} X' \xra{\beta'} X \xra{\gamma'} Z'[1]$ and $X' \xra{u'} Y' \xra{v'} Z \xra{w'} X[1]$ are in $\xi$,
then the triangle $Z' \xra{\alpha} Y' \xra{\beta} Y \xra{\gamma} Z'[1]$ is also in $\xi$.
\end{lem}

\vskip5pt

The following two results can be found in \cite{Y}, without details of proofs.

\begin{lem}\label{monicepimorphism} \ {\rm ([Y, Proposition 2.5])} \ Let $(\T, [1], \tri)$ be a triangulated category, $\xi$ a proper class of triangles, and $u: X \xra{} Y$ and $v: Y \xra{} Z$ morphisms in $\T$.

  \vskip5pt

  $(1)$ \  If $vu$ is a $\xi$-proper monomorphism, then so is $u$.

  \vskip5pt

  $(2)$ \  If $vu$ is a $\xi$-proper epimorphism, then so is $v$.
\end{lem}

\begin{proof}  \ We only prove $(1)$, and the assertion $(2)$ can be similarly proved.
Embedding  $u$ into a triangle $X \xra{u} Y \xra{} Z' \xra{} X[1]$ and taking the cobase change of the triangle
$Y \xra{} Z' \xra{} X[1] \xra{-u[1]} Y[1]$
along $v: Y \xra{} Z$ one gets the following commutative diagram
  \[
    \xymatrix@=0.6cm{
      & X'[-1] \ar@{=}[r]\ar[d] & X'[-1] \ar[d] & \\
    X \ar[r]^-{u}\ar@{=}[d] & Y \ar[d]^-{v}\ar[r] & Z' \ar[r]\ar[d]^-\alpha & X[1]\ar@{=}[d] \\
    X \ar[r]^-{vu} & Z \ar[r]\ar[d] & Y' \ar[r]\ar[d] & X[1]\ar[d]^-{u[1]} \\
      & X' \ar@{=}[r] &  X'\ar[r] & Y[1] &  \\
    }
  \]
Since  $vu$ is a $\xi$-proper monomorphism, the triangle in the third row is in $\xi$. Viewing this diagram as the base change diagram of
the triangle $X \xra{vu} Z \xra{} Y' \xra{} X[1]$ along $\alpha$ one sees that the triangle in the second row is in $\xi$, since $\xi$ is closed under the base changes.
Thus $u$ is a $\xi$-proper monomorphism.
\end{proof}

\begin{lem}\label{compclose} \ {\rm ([Y, Proposition 2.4])}
  \ Let $(\T, [1], \tri)$ be a triangulated category, and $\xi$ a proper class of triangles.

  \vskip5pt

  $(1)$ \ The class of $\xi$-proper monomorphisms is closed under compositions.

  \vskip5pt

  $(2)$ \ The class of $\xi$-proper epimorphisms is closed under compositions.
\end{lem}

\begin{proof}   \ We only prove $(1)$, and the assertion $(2)$ can be similarly proved. Let $u:X \xra{} Y$ and $v: Y \xra{} Z$ be $\xi$-proper monomorphisms.
Embed $u$ and $v$ into triangles $X \xra{u} Y \xra{w} W \xra{} X[1]$ and $Y \xra{v} Z \xra{} K \xra{} Y[1]$.
By definition both the triangles are in $\xi$. Taking the cobase change of the triangle $Y \xra{v} Z \xra{} K \xra{} Y[1]$ along $w: Y \xra{} W$ one gets the following commutative diagram
  \[
    \xymatrix@=0.6cm{
      & X \ar@{=}[r]\ar[d]^-{u} & X \ar[d]^-{vu} & \\
    K[-1] \ar[r]\ar@{=}[d] & Y \ar[d]^-{w}\ar[r]^-{v} & Z \ar[r]\ar[d] & K\ar@{=}[d] \\
    K[-1] \ar[r] & W \ar[r]\ar[d] & E \ar[r]\ar[d] & K\ar[d] \\
      & X[1] \ar@{=}[r] &  X[1]\ar[r]^-{-u[1]} & Y[1]}
  \]
where the two rows and the two columns are triangles.  By Lemma \ref{xiclose}(2) the triangle $X \xra{vu} Z \xra{} E \xra{} X[1]$ is in $\xi$, hence $vu$ is a $\xi$-proper monomorphism.
\end{proof}

\begin{lem}\label{4tri}
  \ Let $(\T, [1], \tri)$ be a triangulated category, $\xi$ a proper class of triangles, and $\alpha: A \xra{} B$ and $\beta: B \xra{} C$  morphisms in $\T$.

  \vskip5pt

  $(1)$ \ If $\alpha$ is a $\xi$-proper epimorphism, then there exists a triangle in $\xi:$
  $$\hok \alpha \xra{} \hok \beta \alpha \xra{} \hok \beta \xra{} \hok \alpha [1].$$

  \vskip5pt

  $(2)$ \ If $\beta$ is a $\xi$-proper monomorphism, then there exists a triangle in $\xi:$
  $$\hoc \alpha \xra{} \hoc \beta \alpha \xra{} \hoc \beta \xra{} \hoc \alpha [1].$$

  \vskip5pt

  $(3)$ \ If $\alpha$ is a $\xi$-proper monomorphism and $\beta \alpha$ is a $\xi$-proper epimorphism, then there exists a triangle in $\xi:$
  $$\hok \beta \alpha \xra{} \hok \beta \xra{} \hoc \alpha \xra{} \hok \beta \alpha [1].$$

  \vskip5pt

  $(4)$ \ If $\beta$ is a $\xi$-proper epimorphism and $\beta \alpha$ is a $\xi$-proper monomorphism, then there exists a triangle in $\xi:$
  $$\hok \beta \xra{} \hoc \alpha \xra{} \hoc \beta \alpha \xra{} \hok \beta [1].$$
\end{lem}

\begin{proof}
  \ We only prove $(1)$ and $(3)$. The assertions $(2)$ and $(4)$ can be proved similarly.

  \vskip5pt

  $(1)$ \ Embed $\alpha$ into a triangle $\hok \alpha \xra{} A \xra{\alpha} B \xra{} \hok \alpha [1]$.
  Since $\alpha$ is a $\xi$-proper epimorphism, this triangle is in $\xi$.
  Also, embed $\beta$ into a triangle $C[-1] \xra{} \hok \beta \xra{k} B \xra{\beta} C$.
  Taking the base change of triangle $\hok \alpha \xra{} A \xra{\alpha} B \xra{} \hok \alpha [1]$ along $k: \hok \beta \xra{} B$, one gets a commutative diagram
  \[
    \xymatrix@=0.6cm{
      & C[-1] \ar@{=}[r]\ar[d] & C[-1] \ar[d] & \\
    \hok \alpha \ar[r]\ar@{=}[d] & \hok \beta \alpha \ar[d]\ar[r] & \hok \beta \ar[r]\ar[d]^-{k} & \hok \alpha [1] \ar@{=}[d] \\
    \hok \alpha \ar[r] & A \ar[r]^-{\alpha}\ar[d]^-{\beta \alpha} & B \ar[r]\ar[d]^-{\beta} & \hok \alpha [1]\ar[d] \\
      & C \ar@{=}[r] &  C\ar[r] &  \hok \beta \alpha[1]
    }
  \]
 Since $\xi$ is closed under the base changes, it follows that the triangle
 $\hok \alpha \xra{} \hok \beta \alpha \xra{} \hok \beta \xra{} \hok \alpha [1]$ is in $\xi$ .

  \vskip5pt

  $(3)$ \ Embed $\alpha$ and $\beta \alpha$ into triangles $A \xra{\alpha} B \xra{} \hoc \alpha \xra{} A[1]$ and
  $\hok \beta \alpha \xra{\gamma} A \xra{\beta \alpha} C \xra{} \hok \beta \alpha [1]$. By the assumption the both triangles are  in $\xi$.
   Embed $\beta$ into a triangle $\hok \beta \xra{k} B \xra{\beta} C \xra{} \hok \beta [1]$.
Taking the base change of triangle $\hoc \alpha [-1] \xra{}  A \xra{\alpha} B \xra{} \hoc \alpha$
along $k: \hok \beta \xra{} B$ one gets a commutative diagram
  \[
    \xymatrix@=0.6cm{
      & C[-1] \ar@{=}[r]\ar[d] & C[-1] \ar[d] & \\
      \hoc \alpha [-1] \ar[r]\ar@{=}[d] & \hok \beta \alpha \ar[d]^-\gamma\ar[r] & \hok \beta \ar[r]\ar[d]^-{k} & \hoc \alpha \ar@{=}[d] \\
      \hoc \alpha [-1] \ar[r] & A \ar[r]^-{\alpha}\ar[d]^-{\beta \alpha} & B \ar[r]\ar[d]^-{\beta} & \hoc \alpha\ar[d] \\
      & C \ar@{=}[r] &  C\ar[r] &  \hok \beta \alpha [1]
    }
  \]
Viewing this diagram as the cobase change diagram of the triangle $\hok \beta \alpha \xra{} \hok \beta \xra{} \hoc \alpha \xra{} \hok \beta \alpha [1]$
along $\gamma$, by Lemma \ref{carsq}(2) one sees that the triangle $\hok \beta \alpha \xra{} \hok \beta \xra{} \hoc \alpha \xra{} \hok \beta \alpha [1]$ is in $\xi$.
\end{proof}

\subsection{\bf Model structures}
{\it A closed model structure} on a category $\mathcal M$ is a triple
\ $(\CoFib$, \ $\Fib$,  \ $\Weq)$ of classes of morphisms,
where the morphisms in the three classes are respectively called {\it cofibrations, fibrations, and weak equivalences},
satisfying {\bf Two out of three axiom}, {\bf Retract axiom},   {\bf Lifting axiom}, and {\bf Factorization axiom}. See D. Quillen [Q1] and M. Hovey [Hov1].
The morphisms in $\TCoFib: = \CoFib\cap \Weq$ and $\TFib: = \Fib\cap \Weq$
are called {\it trivial cofibrations} and  {\it trivial fibrations}, respectively.

\vskip5pt

Following [Hov1] and also [Hir], we will call a closed model structure just as {\it a model structure}. But then a model structure here
is different from  ``a model structure" in the sense of [Q1]:
it is ``a model structure" in [Q1], but the converse is not true. The following facts are in the axioms of ``a model structure" in [Q1], thus one has

\begin{fact} \label{elementpropmodel} \ Let  $(\CoFib$, \ $\Fib$,  \ $\Weq)$ be a
model structure on category $\mathcal M$. Then

$(1)$ \ Both  $\CoFib$ and $\Fib$ are closed under compositions.

$(2)$ \ Isomorphisms are fibrations, cofibrations, and weak equivalences.

$(3)$ \ Cofibrations are closed under push-outs, i.e., if
$$\xymatrix@R=0.4cm{\bullet\ar[r]^-i\ar[d] & \bullet \ar@{.>}[d] \\
\bullet \ar@{.>}[r]^-{i'} & \bullet}
$$
is a push-out square with $i\in \CoFib$, then $i'\in \CoFib$.

Also, trivial cofibrations are closed under push-outs.

$(4)$ \ Fibrations and trivial fibrations are closed under pull-backs.\end{fact}

\vskip5pt

A striking property of a model structure is that any two classes of $\CoFib$, \  $\Fib$, \  $\Weq$ uniquely determine the third.

\vskip5pt

\begin{lem}\label{quillenlifting} {\rm ([Q2, p.234])} \ Let $(\CoFib, \ \Fib,  \ \Weq)$ be a model structure on category $\mathcal M$. Then

\vskip5pt

$(1)$ \ Cofibrations are precisely those morphisms which have the left lifting property with respect to all the trivial fibrations. That is,
a morphism $i: A\longrightarrow B$ is in $\CoFib$ if and only if for any commutative square
$$
\xymatrix@R=0.5cm{
A\ar[r]^-{u}\ar[d]_-{i} & C \ar[d]^-{p}\\
B\ar[r]^-v & D
}$$
with $p\in \TFib$, there is a lifting $\lambda: B \longrightarrow C$ such that $u = \lambda \circ i$ and $v = p\circ \lambda$.

\vskip5pt

$(2)$ \ Trivial cofibrations are precisely those morphisms which have the left lifting property  with respect to all the fibrations.

\vskip5pt

$(3)$ \ Fibrations are precisely those morphisms which have the right lifting property with respect to all the trivial cofibrations.

\vskip5pt

$(4)$ \ Trivial fibrations are precisely those morphisms which have the right lifting property with respect to all the cofibrations.

\vskip5pt

$(5)$ \ $\Weq = \TFib\circ \TCoFib.$

\end{lem}

For a model structure $(\CoFib$, \ $\Fib$,  \ $\Weq)$ on pointed category $\mathcal M$ (i.e., a category with zero object),
an object $X$ is {\it trivial} if $0 \longrightarrow X $ is a weak equivalence, or, equivalently,
$X\longrightarrow 0$ is a weak equivalence.
It is {\it cofibrant} if $0\longrightarrow X$ is a cofibration, and it is {\it fibrant} if $X\longrightarrow 0$ is a fibration.
An object is {\it trivially cofibrant} (respectively, {\it trivially fibrant}) if it is both trivial and  cofibrant (respectively, fibrant).

\section{\bf Cotorsion pairs in a triangulated category with a proper class}

Throughout this section, let $(\T, [1], \tri)$ be a triangulated category, and $\xi$ a proper class of triangles.

\subsection{Abelian group $\Ext_{\xi}(X, Y)$} \ For objects $X$ and $Y$,
denote by $\xi^*(X, Y)$  the class of triangles $Y \xra{} E \xra{} X \xra{} Y[1]$ in $\xi$. Define a relation on $\xi^*(X, Y)$:
for elements $(T_i): Y \xra{} E_i \xra{} X \xra{} Y[1] \ (i = 1, 2)$  in $\xi^*(X, Y)$, $(T_1) \sim (T_2)$ if and only if there is a morphism of triangles:
\[\xymatrix@R=0.5cm{
    (T_1): & Y \ar[r]^-{}\ar@{=}[d] & E_1 \ar[d]^-{g}\ar[r]^-{} & X \ar[r]^-{}\ar@{=}[d] & Y[1]\ar@{=}[d] \\
    (T_2): & Y \ar[r]^-{} & E_2 \ar[r]^-{} & X \ar[r]^-{} & Y[1]. \\
    }
\]
This is an equivalence relation.  Denote by $\Ext_{\xi}(X, Y)$ the corresponding class $\xi^*(X, Y)/\sim$ of equivalence classes.
Using base and cobase change, one can define the Baer sum in $\Ext_{\xi}(X, Y)$
such that $\Ext_{\xi}(X, Y)$ becomes an abelian group and $\Ext_{\xi}(-,-) : \T^{op} \times \T \xra{} \mathcal{A}\mathrm{b}$ is an additive bifunctor. See [B]. In particular,
the zero element of $\Ext_{\xi}(X, Y)$ is the splitting triangle, and
$$\Ext_{\xi}(X_1 \oplus X_2, Y) \cong \Ext_{\xi}(X_1, Y) \oplus \Ext_{\xi}(X_2, Y)$$ and $$\Ext_{\xi}(X, Y_1 \oplus Y_2) \cong \Ext_{\xi}(X, Y_1) \oplus \Ext_{\xi}(X, Y_2).$$

\vskip5pt

The Extension-Lifting Lemma has been observed in  \cite[VIII, 3.1]{BR} for abelian categories.
It has independent interest.

\begin{thm} \label{extlifting} {\rm (Extension-Lifting Lemma for triangulated categories)} \
\ Let $(\T, [1], \tri)$ be a triangulated category,  $\xi$ a proper class of triangles, and $X$ and $Y$  objects of $\T$.
Then $\Ext_{\xi} (X,Y) = 0$ if and only if for any commutative diagram
    \begin{equation}\label{ext}
    \begin{gathered}
    \xymatrix@R=0.5cm{
      & A \ar[r]^-{i} \ar[d]_-{\alpha} & B \ar[r]^-{d}\ar[d]_-{\beta} & X \ar[r]^-{u} & A[1] \\
    Y \ar[r]^-{c} & C \ar[r]^-{p} & D \ar[r]^-{v} & Y[1] &
    }
    \end{gathered}
    \end{equation}
where the two rows are triangles in $\xi$, there is a lifting $\lambda: B \longrightarrow C$ such that $\alpha = \lambda i$ and $\beta = p \lambda$.
\end{thm}

\begin{proof} \ Assume that for any commutative diagram of triangles as in (\ref{ext}) admits a lifting.
For any triangle $Y \xra{c} C \xra{p} X \xra{} Y[1]$ in $\xi$, by the assumption the commutative diagram
    \[
      \xymatrix@R=0.5cm{
        & C \ar[r]^-{\left(\begin{smallmatrix} 1 \\ 0 \end{smallmatrix}\right)} \ar@{=}[d] & C \oplus X \ar[r]^-{(0,1)}\ar[d]_-{(p,1)} & X \ar[r]^-{0} & C[1] \\
      Y \ar[r]^-{c} & C \ar[r]^-{p} & X \ar[r]^-{v} & Y[1] &
      }
    \]
admits a lifting $\lambda=(a,b): C\oplus X \longrightarrow C$
such that $\Id_C = (a,b) \left(\begin{smallmatrix} 1 \\ 0 \end{smallmatrix}\right) = a$ and $(p,1) = p(1,b) = (p,pb)$.
Thus $pb = \Id_X$, which means the given triangle splits.

    \vskip5pt

Conversely, assume that $\Ext_{\xi} (X,Y)= 0$.  Given a commutative diagram \eqref{ext} of triangles in $\xi$,
then there is a morphism $\gamma: X \xra{} Y[1]$ such that the diagram
    \[
      \xymatrix{
      A \ar[r]^-{i} \ar[d]_-{\alpha} & B \ar[r]^-{d}\ar[d]_-{\beta} & X \ar[r]^-{u}\ar[d]^-{\gamma} & A[1]\ar[d]^-{\alpha [1]} \\
      C \ar[r]^-{p} & D \ar[r]^-{v} & Y[1]\ar[r]^-{-c[1]} & C[1]&
      }
    \]
commutes. Taking the base change of $Y \xra{c} C \xra{p} D \xra{v} Y[1]$ along $\beta :B \xra{} D$ one gets a commutative diagram
    $$\xymatrix@R=0.5cm{
       & L \ar@{=}[r]\ar[d] & L \ar[d] & \\
      Y\ar[r]^-{m}\ar@{=}[d] & K\ar[r]^-{k}\ar[d]_-{j} & B\ar[r]^-{n}\ar[d]^{\beta}& Y[1]\ar@{=}[d] \\
      Y\ar[r]^-{c} & C\ar[r]^-{p}\ar[d] & D\ar[r]^-{v}\ar[d] & Y[1]\ar[d]^-{m[1]} \\
      & L[1] \ar@{=}[r] & L[1]\ar[r]  &  K[1]
    }$$
    Since the triangle in the third row is in $\xi$,
the triangle $Y \xra{m} K \xra{k} B \xra{\eta} Y[1]$ in the second row is also in $\xi$.
Since $n = v\beta = \gamma d$, one has a commutative diagram of triangle
    \[
      \xymatrix@R=0.6cm{
      X[-1]\ar[d]_-{\gamma [-1]} \ar[r]^-{-u[-1]}  & A \ar[r]^-{i} & B \ar[r]^-{d}\ar@{=}[d] & X \ar[d]_-{\gamma} \\
      Y\ar[r]^-{m} & K\ar[r]^-{k} & B\ar[r]^-{n}& Y[1]
      }
    \]
and hence there is a morphism $\eta: A \xra{} K$ such that $i=k \eta$ and $- \eta (u[-1]) = m (\gamma [-1])$.
Since $i = k \eta$ is a $\xi$-proper monomorphism, it follows from Lemma \ref{monicepimorphism}(1) that $\eta$ is a $\xi$-proper monomorphism.
By the octahedral axiom one gets a commutative diagram
    \[
    \xymatrix@R=0.6cm{
    A \ar[r]^-{\eta}\ar@{=}[d] & K \ar[r]^-{\theta}\ar[d]_-{k} & E \ar[r]^-{\tau}\ar[d]^-{f} & A[1]\ar@{=}[d] \\
    A \ar[r]^-{i} & B \ar[r]^-{d}\ar[d]_-{n} & X \ar[r]^-{u}\ar[d]^-{g} & A[1]\ar[d]^-{\eta [1]}  \\
      & Y[1] \ar@{=}[r]\ar[d]_-{-m[1]} & Y[1] \ar[d]^-{-h[1]}\ar[r]^-{-m[1]} & K[1] \\
      & K[1]\ar[r]^-{\theta[1]} & E[1]
    }
    \]
    Then one gets the base change  diagram of $Y \xra{h} E \xra{f} X \xra{g} Y[1]$ along $d:B \xra{} X$:
    $$\xymatrix@R=0.6cm{
      & A \ar@{=}[r]\ar[d]^-{\eta} & A \ar[d]^-{i} & \\
     Y\ar[r]^-{m}\ar@{=}[d] & K\ar[r]^-{k}\ar[d]^{\theta} & B\ar[r]^-{n}\ar[d]^{d}& Y[1]\ar@{=}[d] \\
     Y\ar[r]^-{h} & E\ar[r]^-{f}\ar[d]^-{\tau} & X\ar[r]^-{g}\ar[d]^-{u} & Y[1]\ar[d]^-{m[1]} \\
     & A[1] \ar@{=}[r] & A[1]\ar[r]^-{- \eta [1]}  & K[1]
    }$$
Since the triangle $Y \xra{m} K \xra{k} B \xra{n} Y[1]$ in the second row
and the triangle $A \xra{i} B \xra{d} X \xra{u} A[1]$ in the third column
are in $\xi$, it follows from the saturation  of $\xi$ that the triangle  $Y \xra{h} E \xra{f} X \xra{g} Y[1]$ is in $\xi$.
Since $\Ext_{\xi} (X,Y) = 0$, one has $g = 0$.  Then $n = gd = 0$, which means $k$ is a splitting epimorphism,
i.e., there is $l: B \xra{} K$ such that $kl = \Id _B$. Thus $pjl = \beta kl = \beta$.

    \vskip5pt

Since $p(\alpha -jli) = p\alpha -pjli = p\alpha -\beta i =0$, by the commutative diagram
    \[
    \xymatrix@R=0.5cm{
      A \ar@{=}[r] & A \ar[r]\ar[d]_-{\alpha -jli} & 0 \ar[r]\ar[d] & A[1] \\
      Y \ar[r]^-{c} & C \ar[r]^-{p} & D \ar[r] & Y[1]
    }
    \]
    there exists a morphism $\delta : A \xra{} Y$ such that $c \delta = \alpha - jli$.
Taking the cobase change of $A \xra{i} B \xra{d} X \xra{u} A[1]$ along $\delta: A \xra{} Y$ one gets the following commutative diagram
    $$\xymatrix@R=0.6cm{
      & M \ar@{=}[r]\ar[d] & M \ar[d] & \\
     X[-1]\ar[r]^-{-u[-1]}\ar@{=}[d] & A \ar[r]^-{i}\ar[d]^{\delta} & B\ar[r]^-{d}\ar[d]^-{w}& X\ar@{=}[d] \\
     X[-1]\ar[r] & Y\ar[r]^-{\varepsilon}\ar[d] & G \ar[r]\ar[d] & X\ar[d]^-{-u} \\
     & M[1] \ar@{=}[r] & M[1]\ar[r]  & A[1]
    }$$
Since $\xi$ is closed under cobase changes, the triangle $Y \xra{\varepsilon} G \xra{} X \xra{} Y[1]$ is in $\xi$.
Since $\Ext_{\xi} (X,Y) = 0$, $\varepsilon$ is a splitting monomorphism, i.e., there is $\pi : G \xra{} Y$ such that $\pi \varepsilon = \Id_Y$.
Thus $\pi w i = \pi \varepsilon \delta = \delta$.

    \vskip5pt

Now, put $\lambda = jl + c \pi w: B\longrightarrow C$. Then one has $p\lambda = pjl+ pc\pi w = pjl = \beta$ and
$\lambda i = jli + c\pi wi = \alpha -c\delta + c \delta = \alpha$. Thus $\lambda$ is the desired lifting. This completes the proof. \end{proof}

\vskip5pt

For the following result, if $\T$ has enough $\xi$-projective objects, Beligiannis [B, Corollary 4.12] gives a similar result.
However, to  claim the sameness of Theorem \ref{exactseq} and [B, Corollary 4.12], as in an abelian category, one needs the Yoneda extensions and pull-backs and push-outs.
Since in a triangulated category the existence of pull-backs and push-outs can not be guaranteed, we need to use the present form of Theorem \ref{exactseq}, without using $\xi$-projective resolutions.

\begin{thm}\label{exactseq}  \   $(1)$ \  $\Ext_{\xi}(X, Y)$ is a subgroup of $\Hom_\T(X, Y[1])$, in particular, $\Ext_{\xi}(X, Y)$ is a set. In fact, the map
given by \ $(T)\mapsto w$  is an injective group homomorphism,
where $(T)$ is the equivalence class of triangle $Y \xra{} E \xra{} X \xra{w} Y[1]$ in $\xi$. Thus, if $\xi = \tri$ then $\Ext_{\xi}(X, Y)= \Hom_\T(X, Y[1])$.

  \vskip5pt

  $(2)$ \ For any triangle $X_1 \xra{f} X \xra{g} X_2 \xra{h} X_1[1]$ in $\xi$ and object $Y \in \T$, there is an exact sequence of abelian groups $($where $(X_2, Y) = \Hom_\mathcal T(X_2, Y))$
$$(X_2, Y) \xra{(g, Y)} (X, Y) \xra{(f, Y)} (X_1, Y) \xra{\widetilde{h}} \Ext_{\xi}(X_2, Y) \xra{\widetilde{g}} \Ext_{\xi}(X, Y) \xra{\widetilde{f}} \Ext_{\xi}(X_1, Y)$$
where $\widetilde{h}$ sends $\alpha\in\Hom_\mathcal T(X_1, Y)$ to the equivalence class of triangle \  $Y \xra{} E_2 \xra{} X_2 \xra{\alpha[1]\circ h} Y[1];$
and $\widetilde{g}$ is induced by the base change of \  $Y \xra{} E_2 \xra{} X_2 \xra{\eta} Y[1]$ along $g: X \xra{} X_2$
  \[
    \xymatrix@R=0.6cm{
      & K \ar@{=}[r]\ar[d] & K \ar[d] & \\
    Y \ar[r]\ar@{=}[d] & E \ar[d]\ar[r] & X \ar[r]^-{\eta\circ g}\ar[d]^-{g} & Y[1]\ar@{=}[d] \\
    Y \ar[r] & E_2 \ar[r]\ar[d] & X_2 \ar[r]^-\eta\ar[d] & Y[1]\ar[d] \\
      & K[1] \ar@{=}[r] &  K[1]\ar[r] &  E[1].
    }
  \]

\vskip5pt

$(3)$ \ For any triangle $Y_1 \xra{f} Y \xra{g} Y_2 \xra{h} Y_1[1]$ in $\xi$ and object $X \in \T$, there is an exact sequence of abelian groups
$$(X, Y_1) \xra{(X, f)} (X, Y) \xra{(X, g)} (X, Y_2) \xra{\overline{h}} \Ext_{\xi}(X, Y_1) \xra{\overline{f}} \Ext_{\xi}(X, Y) \xra{\overline{{g}}} \Ext_{\xi}(X, Y_2)$$
where $\overline{{h}}$ sends $\alpha\in\Hom_\mathcal T(X, Y_2)$ to the equivalence class of triangle \  $Y_1 \xra{} E_1 \xra{} X \xra{h\alpha} Y[1];$
and $\overline{{f}}$ is induced by the cobase change of \  $Y_1 \xra{} E_1 \xra{} X \xra{w} Y_1[1]$ along $f: Y_1 \xra{} Y$
  \[
    \xymatrix@R=0.6cm{
      & K \ar@{=}[r]\ar[d] & K \ar[d] & \\
    X[-1] \ar[r]^-{-w[-1]}\ar@{=}[d] & Y_1 \ar[d]^-{f}\ar[r] & E_1 \ar[r]\ar[d] & X\ar@{=}[d] \\
    X[-1] \ar[r]^-{-f\circ w[-1]} & Y \ar[r]\ar[d] & E \ar[r]\ar[d] & X\ar[d] \\
      & K[1] \ar@{=}[r] &  K[1]\ar[r] &  Y_1[1].
    }
  \]
\end{thm}

\begin{proof}   \ The assertion (1) directly follows from the definition of the Baer sum of $\Ext_{\xi}(X, Y)$.

\vskip5pt

We only prove $(2)$. The assertion $(3)$ can be similarly proved.

\vskip5pt

Since $\Hom_\mathcal T(-, Y)$ is a cohomological functor,
one has exact sequence $(X_2, Y) \xra{(g, Y)} (X, Y) \xra{(f, Y)} (X_1, Y)$.

\vskip5pt

{\bf Claim 1:} \ $\widetilde{h}: \Hom_\mathcal T(X_1, Y) \longrightarrow \Ext_{\xi}(X_2, Y)$  is well-defined, namely, the triangle
\ $Y \xra{} E_2 \xra{} X_2 \xra{\alpha[1]\circ h} Y[1]$ is in $\xi$.

\vskip5pt

In fact,  by the cobase change of \  $X_1 \xra{f} X \xra{g} X_2 \xra{h} X_1[1]$ along $\alpha: X_1 \xra{} Y$,
one gets a commutative diagram
\[
    \xymatrix@R=0.6cm{
      & Z[-1] \ar@{=}[r]\ar[d] & Z[-1] \ar[d] & \\
    X_2[-1] \ar[r]^-{-h[-1]}\ar@{=}[d] & X_1 \ar[d]^-{\alpha}\ar[r]^-{f} & X \ar[r]^-{g}\ar[d] & X_2\ar@{=}[d] \\
    X_2[-1] \ar[r]^-{-\alpha\circ h[-1]} & Y \ar[r]\ar[d] & E_2 \ar[r]\ar[d] & X_2\ar[d]^-{-h} \\
      & Z \ar@{=}[r] &  Z\ar[r] &  X_1[1]
    }
  \]
Since $\xi$ is closed under the cobase changes, it follows that the triangle
\ $Y \xra{} E_2 \xra{} X_2 \xra{\alpha[1]\circ h} Y[1]$ is in $\xi$.

\vskip5pt

{\bf Claim 2:} \ $\Im (f,Y) = \Ker \widetilde{h}.$  \

\vskip5pt

In fact,  it is clear that $\Im (f, Y) \subseteq \Ker \widetilde{h}$. For $\gamma \in \Hom_\mathcal T(X_1, Y)$ in $\Ker \widetilde{h}$, i.e., $\gamma [1] \circ h = 0: X_2 \longrightarrow Y[1]$.
Then one gets a commutative diagram
$$\xymatrix@R=0.6cm{
    X_1\ar[d]_-{\gamma}\ar[r]^-{f} & X\ar[r]^-{g}\ar@{-->}[d]_-\theta & X_2 \ar[r]^-{h}\ar[d] & X_1[1]\ar[d]^-{\gamma [1]} \\
    Y \ar@{=}[r] & Y \ar[r] & 0 \ar[r] & Y[1]
}$$
Thus $\gamma = \theta \circ f\in \Im (f,Y)$.

\vskip5pt

{\bf Claim 3:} \ $\Im \widetilde{h} = \Ker \widetilde{g}$.

\vskip5pt

In fact,  let $Y \xra{} E_2 \xra{} X_2 \xra{\alpha [1]\circ h} Y[1]$ be a triangle in $\Im \widetilde{h}$, where
$\alpha\in \Hom_\mathcal T(X_1, Y)$.  By the definition of $\widetilde{g}$ one has the following base change diagram
$$\xymatrix@R=0.6cm{
     & K\ar@{=}[r]\ar[d] & K \ar[d] & \\
    Y \ar[r]\ar@{=}[d] & E \ar[r]\ar[d] & X \ar[d]^-{g}\ar[r]^-{\alpha[1]\circ h\circ g} & Y[1]\ar@{=}[d] \\
    Y \ar[r]& E_2 \ar[r]\ar[d] & X_2 \ar[r]^-{\alpha[1]\circ h}\ar[d] & Y[1]\ar[d] \\
     & K[1] \ar@{=}[r] & K[1]\ar[r] & E[1]
}$$
Since $\alpha[1]\circ h\circ g = 0$ one sees that $\Im \widetilde{h} \subseteq \Ker \widetilde{g}$.

\vskip5pt

Conversely, let $Y \xra{} E_2 \xra{} X_2 \xra{\eta} Y[1]$ be a triangle in $\Ker \widetilde{g}$. Then $\eta \circ g = 0$. By the following commutative diagram
$$\xymatrix@R=0.6cm{
    X_1\ar@{-->}[d]_-{\delta} \ar[r]^-{f} & X \ar[d]\ar[r]^-{g} & X_2 \ar[r]^-{h}\ar[d]^-{\eta} & X_1[1]\ar@{-->}[d]^-{\delta[1]} \\
    Y\ar[r] & 0 \ar[r] & Y[1]\ar@{=}[r] & Y[1]
}$$
one sees that $\eta = \delta [1]\circ h \in \Im \widetilde{h}$. Thus $\Ker \widetilde{g} \subseteq \Im \widetilde{h}$.

\vskip5pt

{\bf Claim 4:} \ $\Im \widetilde{g} = \Ker \widetilde{f}$.

\vskip5pt

In fact,  for any element $Y \xra{} E_2 \xra{} X_2 \xra{\eta} Y[1]$ in $\Ext_{\xi}(X_2,Y)$,
taking the base change  along $g: X \xra{} X_2$ one gets a triangle $Y \xra{} E \xra{} X \xra{\eta\circ g} Y[1]$ in $\Ext_{\xi}(X, Y)$, and then
taking the base change  along $f: X_1 \xra{} X$ one gets a triangle $Y \xra{} E_1 \xra{} X_1 \xra{\eta\circ g\circ f} Y[1]$ in $\Ext_{\xi}(X_1, Y)$
with $\eta\circ g\circ f = 0$. This shows $\widetilde{f}\circ \widetilde{g} = 0$.

  \vskip5pt

Finally, we show $\Ker \widetilde{f} \subseteq \Im \widetilde{g}$.
Let $Y \xra{u} E \xra{v} X \xra{w} Y[1]$ be a triangle in $\Ker \widetilde{f}$. By taking the base change
one gets a splitting triangle $Y \xra{u_1} E_1 \xra{v_1} X_1 \xra{w\circ f = 0} Y[1]$. See commutative diagram
  \[
    \xymatrix@R=0.6cm{
      & X_2[-1] \ar@{=}[r]\ar[d] & X_2[-1] \ar[d]^-{-h[-1]} & \\
    Y \ar[r]^{u_1}\ar@{=}[d] & E_1 \ar[d]^-{f'}\ar[r]^-{v_1} & X_1 \ar[r]^-{0}\ar[d]^-{f} & Y[1]\ar@{=}[d] \\
    Y \ar[r]^-{u} & E \ar[r]^-{v}\ar[d] & X \ar[r]^-{w}\ar[d]^-{g} & Y[1]\ar[d]^-{u_1[1]} \\
      & X_2 \ar@{=}[r] &  X_2\ar[r] &  E_1 [1].
    }
  \]
Thus there is a morphism $v_1': X_1 \xra{} E_1$ such that $v_1\circ v_1 ' = \Id_{X_1}$.
Then $v f' v_1 ' = f v_1 v_1' = f$. By the octahedral axiom one has a commutative diagram
  \[
    \xymatrix@R=0.6cm{
    X_1 \ar[r]^-{f' v_1 '}\ar@{=}[d] & E \ar[r]^-{m}\ar[d]_-{v} & E_2 \ar[r]^-{n}\ar[d]^-{v_2} & X_1[1]\ar@{=}[d] \\
    X_1 \ar[r]^-{f} & X \ar[r]^-{g}\ar[d]_-{w} & X_2 \ar[r]^-{h}\ar[d]^-{w_2} & X_1[1]\ar[d]^-{(f'v_1')[1]} \\
      & Y[1] \ar@{=}[r]\ar[d]_-{-u[1]} & Y[1] \ar[d]^-{-u_2[1]}\ar[r]^-{-u[1]} & E[1]\\
      & E[1] \ar[r]^-{m[1]} & E_2[1] &
    }
  \]
  Then one gets the base change diagram
  \[
    \xymatrix@R=0.6cm{
      & X_1 \ar@{=}[r]\ar[d]_-{f' v_1 '} & X_1 \ar[d]_-{f} & \\
    Y \ar[r]^-{u}\ar@{=}[d] & E \ar[d]_-{m}\ar[r]^-{v} & X \ar[r]^-{w}\ar[d]_-{g} & Y[1]\ar@{=}[d] \\
    Y \ar[r]^-{u_2} & E_2 \ar[r]^-{v_2}\ar[d]_-{n} & X_2 \ar[r]^-{w_2}\ar[d]_-{h} & Y[1]\ar[d]^-{u[1]} \\
      & X_1[1] \ar@{=}[r] &  X_1[1]\ar[r]^-{-(f'v_1')[1]} &  E[1]
    }
  \]
The triangle $Y \xra{u} E \xra{v} X \xra{w} Y[1]$ in the second row is in $\Ker \widetilde{f}$, thus it is in $\xi$.
The triangle $X_1 \xra{f} X \xra{g} X_2 \xra{h} X_1[1]$ in the third column is in $\xi$, by the assumption.
Since $\xi$ is saturated, it follows that the triangle $Y \xra{u_2} E_2 \xra{v_2} X_2 \xra{w_2} Y[1]$ in
the third row is in $\xi$, and by definition its image under $\widetilde{g}$ is $Y \xra{u} E \xra{v} X \xra{w} Y[1]$.
This proves  $\Ker \widetilde{f} \subseteq \Im \widetilde{g}$.
\end{proof}

\subsection {Cotorsion pairs with respect to a proper class of triangles}
Let $(\T, [1], \tri)$ be a triangulated category, and $\xi$ a proper class of triangles. For classes $\X$ and $\Y$ of objects of $\T$,
we write $\Ext_{\xi} (\X, \Y)= 0$ if $\Ext_{\xi} (X,Y)= 0$ for all $X \in \X$ and $Y \in \Y$. Similarly,
for object $Z$ we use the expressions $\Ext_{\xi} (\X, Z)= 0$  and $\Ext_{\xi} (Z, \X)= 0$.

\vskip5pt

Denote by $\X^\perp$ the class of objects $Z$ satisfying $\Ext_{\xi} (\X, Z)= 0$; and by \ $^\perp\X$ the class of objects $Z$ satisfying $\Ext_{\xi} (Z, \X)= 0$.

\begin{defn}\label{cpt} \ Let $(\T, [1], \tri)$ be a triangulated category, and $\xi$ a proper class of triangles.

\vskip5pt

(1) \  A pair $(\X, \Y)$ of classes of objects of $\T$  is a cotorsion pair with respect to $\xi$,
provided that $\X^\perp = \Y$ and \ $^\perp\Y = \X$.
\vskip5pt

(2) \   A cotorsion pair $(\X, \Y)$ with respect to $\xi$ is complete,
provided that  for any $T \in \T$ the following conditions are satisfied:

\hskip18pt (i) \ \ there is a triangle \ $Y \xra{} X \xra{} T \xra{} Y[1]$ in $\xi$ with $X\in \X$ and $Y\in \Y;$

\hskip18pt (ii) \  there is a triangle \ $T \xra{} Y' \xra{} X' \xra{} T[1]$ in $\xi$ with $Y' \in \Y$ and $X' \in \X$.

\vskip5pt

(3) \ A cotorsion pair $(\X, \Y)$ with respect to $\xi$ is hereditary, provided that the following conditions are satisfied:

\hskip18pt (i) \ \ $\mathcal X$ is closed under the hokernel of $\xi$-proper epimorphisms,
i.e., if $K\xra{} X_1 \xra{} X_2 \xra{} K[1]$ is a triangle in $\xi$ with $X_1, X_2 \in \X$, then $K \in \X;$

\hskip18pt (ii) \ $\mathcal Y$ is closed under the hocokernel of $\xi$-proper monomorphisms, i.e., if $Y_1 \xra{} Y_2 \xra{} C\xra{} Y_1[1]$ is a triangle in $\xi$ with $Y_1, Y_2 \in \Y$, then $C\in \Y$.
\end{defn}

Sometimes it is convenient to use the following description of a complete cotorsion pair with respect to a proper class of triangles.

\begin{fact}\label{completectp} \ Let $(\T, [1], \tri)$ be a triangulated category, $\xi$ a proper class of triangles,
and $\mathcal{X}$ and $\mathcal{Y}$ classes of objects of $\T$ which are closed under isomorphisms and direct summands.
Then $(\X, \Y)$ is a complete cotorsion pair with respect to $\xi$ if and only if the following conditions are satisfied$:$

\hskip18pt {\rm (i)} \ \ \ $\Ext_\xi(\X, \Y) = 0;$

\hskip18pt {\rm (ii)} \ \ for any $T \in \T$, there is a triangle \ $Y \xra{} X \xra{} T \xra{} Y[1]$ in $\xi$ with $X\in \X$ and $Y\in \Y;$

\hskip18pt {\rm (iii)} \  for any $T \in \T$, there is a triangle \ $T \xra{} Y' \xra{} X' \xra{} T[1]$ in $\xi$ with $Y' \in \Y$ and $X' \in \X$.
\end{fact}

If $(\X, \Y)$ is a cotorsion pair with respect to $\xi$, then clearly $\X$ and $\Y$ are closed under isomorphisms, direct summands, and finite coproducts.
By Theorem \ref{exactseq} one has

\begin{cor}\label{extclose}
  \ Let $(\T, [1], \tri)$ be a triangulated category, $\xi$ a proper class of triangles, and $(\X, \Y)$ a cotorsion pair with respect to $\xi$.
  Then $\X$ and $\Y$ are closed under extensions, i.e., if $X_1 \xra{} X \xra{} X_2 \xra{} X_1[1]$ and $Y_1 \xra{} Y \xra{} Y_2 \xra{} Y_1[1]$ are triangles in $\xi$ with $X_1, X_2 \in \X$ and $Y_1, Y_2 \in \Y$, then $X \in \X$ and $Y \in \Y$.
\end{cor}

\begin{rem}\label{cotstructure} The notion of a co-$t$-structure on a triangulated category has been introduced independently by M. V. Bondarko [Bon] (under the name of weight structure) and D. Pauksztello [P].
Recall that a pair $(\X, \Y)$ of full subcategories of triangulated category $(\T, [1], \tri)$ is a co-$t$-structure on $\T$ if the following
conditions are satisfied$:$

\vskip5pt

$(0)$ \ $\X$ and $\Y$ are closed under isomorphisms and direct summands;

$(1)$ \ $\X[-1] \subseteq \X$ and $\Y[1] \subseteq \Y$;

$(2)$ \ $\Hom_{\T}(\X, \Y[1]) = 0$;

$(3)$ \ For any object $T \in \T$ there exists a triangle $X[-1] \xra{} T \xra{} Y \xra{} X$ with $X \in \X$ and $Y \in \Y$.

\vskip5pt

If $\xi = \tri$, then it is not hard to see that
a hereditary complete cotorsion pair with respect to $\xi$ is precisely a co-$t$-structure.
Thus, a classical example $(\X, \Y)$ of co-$t$-structure on homotopy category  $K(\A)$ of additive category $\A$,
is a hereditary complete cotorsion pair in $K(\A)$, where
\begin{align*}\X & = \{ \mbox{complexes in} \ K(\A) \  \mbox{isomorphic to complexes} \ C \ | \  C^i = 0 \ \forall \ i < 0 \}, \ \ \mbox{and} \\
\Y & = \{ \mbox{complexes in} \ K(\A) \  \mbox{isomorphic to complexes} \ C \ | \  C^i = 0 \ \forall \ i > 0 \}.\end{align*}

But we emphasize that, if $\xi\ne \tri$,
then a hereditary complete cotorsion pair with respect to $\xi$ is completely different from a co-$t$-structure.
\end{rem}

\begin{rem}\label{tstructure} The notion of a $t$-structure on a triangulated category has been introduced by A. A. Beilinson, J. Bernstein and P. Deligne \cite{BBD}.
Recall that a pair $(\X, \Y)$ of full subcategories of triangulated category $(\T, [1], \tri)$ is a $t$-structure on $\T$ if the following
  conditions are satisfied$:$

  \vskip5pt

  $(1)$ \ $\X[1] \subseteq \X$ and $\Y[-1] \subseteq \Y$;

  $(2)$ \ $\Hom_{\T}(\X, \Y[-1]) = 0$;

  $(3)$ \ For any object $T \in \T$ there exists a triangle $X \xra{} T \xra{} Y \xra{} X[1]$ with $X \in \X$ and $Y \in \Y[-1]$.

  \vskip5pt

  If $\xi = \tri$, then it is not hard to see that
  a $t$-structure $(\X, \Y)$ gives rise to a complete cotorsion pair $(\X, \Y[-2])$ with respect to $\xi$, but in general
  the complete cotorsion pair $(\X, \Y[-2])$ is not hereditary.

\vskip5pt

Thus, a classical example $(\X, \Y)$ of $t$-structure on bounded derived category $D^b(\A)$ of abelian category $\A$
  gives rise to a complete cotorsion pair $(\X, \Y[-2])$, where
  \begin{align*}\X & = \{\mbox{a bounded complexes} \ C \ \mbox{in} \ D^b(\A)  \ | \  {\rm H}^i(C) = 0,  \ \forall \ i > 0 \}, \ \ \mbox{and} \\
  \Y & = \{\mbox{a bounded complexes} \ C \ \mbox{in} \ D^b(\A)  \ | \  {\rm H}^i(C) = 0,  \ \forall \ i < 0 \}.\end{align*}
But this complete cotorsion pair $(\X, \Y[-2])$ is not hereditary.

  \end{rem}

\subsection{Description of heredity of a complete cotorsion pair.}

\begin{thm} \label{heredity} \  Let $(\T, [1], \tri)$ be a triangulated category,  $\xi$ a proper class of triangles,
and $(\X, \Y)$ be a complete cotorsion pair in $\T$ with respect to $\xi$.
Then the following are equivalent$:$

\vskip5pt

$(1)$ \ $(\X, \Y)$ is a hereditary cotorsion pair$;$

\vskip5pt

$(2)$ \ $\mathcal X$ is closed under the hokernel of $\xi$-proper epimorphisms$;$

\vskip5pt

$(3)$  \ $\mathcal Y$ is closed under the hocokernel of $\xi$-proper monomorphisms.
\end{thm}
\begin{proof} $(2) \Longrightarrow(3):$ \ Suppose that $Y_1 \xra{} Y_2 \xra{\sigma} C \xra{} Y_1[1]$ is a triangle in $\xi$ with $Y_1, Y_2 \in \Y$. We need to show $C \in \Y$.

    \vskip5pt

    {\rm \bf Claim:} \ Any morphism $f: X \xra{} C$ with $X \in \X$ factors through an object of $\omega$. In fact, by Theorem \ref{exactseq} one has an exact sequence
    $\Hom_\mathcal T(X, Y_2) \xra{(X,\sigma)} \Hom_\mathcal T(X, C) \xra{} \Ext_{\xi}(X, Y_1) = 0$, and thus there is a morphism $g\in \Hom_\mathcal T(X, Y_2)$ such that $f = \sigma g$.
    By the completeness of cotorsion pair $(\X, \Y)$, one can take a triangle $Y \xra{} X' \xra{} Y_2 \xra{} Y[1]$ in $\xi$ with $X' \in \X$ and $Y \in \Y$. Then $X' \in \X \cap \Y = \omega$.
From exact sequence $\Hom_\mathcal T(X, X') \xra{} \Hom_\mathcal T(X, Y_2) \xra{} \Ext_{\xi}(X,Y') = 0$ one sees that $g$ factors through $X'$, and hence $f$ factors through $X'\in \omega$.

    \vskip5pt

One can take a triangle $C \xra{p} Y_3 \xra{q} X_1 \xra{r} C[1]$ with $Y_3\in \Y$ and $X_1\in \X$,
and a triangle $Y_4 \xra{u} X_2 \xra{v} Y_3 \xra{w} Y_4[1]$ in $\xi$ with $X_2 \in \X$ and $Y_4 \in \Y$.
Then $qv$ is also a $\xi$-proper epimorphism. By the octahedral axiom one has a commutative diagram
    \[
      \xymatrix@R=0.6cm{
      X_2 \ar[r]^-{v}\ar@{=}[d] & Y_3 \ar[r]^-{w}\ar[d]_-{q} & Y_4[1] \ar[r]^-{-u[1]}\ar[d]_-{g[1]} & X_2[1]\ar@{=}[d] \\
      X_2 \ar[r]^-{qv} & X_1 \ar[r]^-{\tau}\ar[d]_-{r} & X_3[1] \ar[r]^-{-\gamma[1]}\ar[d]_-{f[1]} & X_2[1]\ar[d]^-{v[1]} \\
        & C[1] \ar@{=}[r]\ar[d]_-{-p[1]} & C[1] \ar[d]\ar[r]^-{-p[1]} & Y_3[1]\\
        & Y_3[1] \ar[r]^-{w[1]} & Y_4[2] &
      }
    \]
    Then one gets a triangle $X_3 \xra{\gamma} X_2 \xra{qv} X_1 \xra{\tau} X_3[1]$ in $\xi$ with $X_1, X_2\in \xi$.
Since by assumption $\mathcal X$ is closed under the hokernel of $\xi$-proper epimorphisms, one sees that $X_3 \in \X$.
By {\bf Claim} above, $f$ factors through $W \in \omega$, say $f = mn$ with $n\in \Hom_\mathcal T(X_3, W)$ and $m\in \Hom_\mathcal T(W, C)$.
From exact sequence $\Hom_\mathcal T(X_2, W) \xra{(\gamma, W)} \Hom_\mathcal T(X_3, W) \xra{} \Ext_{\xi}(X_1, W) = 0$,
one sees that  $n$ factors through $\gamma$,  and then $f$ factors through $\gamma$.
From the following commutative diagram and by Fact \ref{2factor} one sees that  $\Id_{X_1}$ factors through $q$.
    $$\xymatrix{
      X_3 \ar[r]^-{\gamma}\ar[d]_-{f} & X_2 \ar[r]^-{qv}\ar[d]^-{v} & X_1 \ar@{=}[d]\ar[r]^-{\tau} & X_3[1]\ar[d]^-{f[1]} \\
      C \ar[r]^-{p} & Y_3 \ar[r]^-{q} & X_1 \ar[r]^-{r} & C[1]
    }$$
That is $q$ is a splitting epimorphism, and $C$ is a direct summand of $Y_3$. Thus $C \in \Y$.

    \vskip5pt

    $(3) \Longrightarrow (2):$  Suppose that $K \xra{} X_1 \xra{} X_2 \xra{} K[1]$ is a triangle in $\xi$ with $X_1, X_2 \in \X$.
    One can show $K \in \X$, similarly as above. In fact, any morphism $f: K \xra{} Y$ with $Y \in \Y$ factors through an object in $\omega$.
Taking triangles $Y_1 \xra{t} X_3 \xra{} K \xra{} Y_1[1]$ and $X_3 \xra{s} Y_2 \xra{} X_4 \xra{} X_3[1]$ in $\xi$ with $X_3, X_4 \in \X$ and $Y_1, Y_2 \in \Y$,
one gets a  commutative diagram
    \[
      \xymatrix@R=0.6cm{
      Y_1 \ar[r]^-{t}\ar@{=}[d] & X_3 \ar[r]\ar[d]_-{s} & K \ar[r]\ar[d]^-{f} & Y_1[1]\ar@{=}[d] \\
      Y_1 \ar[r]^-{st} & Y_2 \ar[r]^-{\alpha}\ar[d] & Y_3 \ar[r]\ar[d] & Y_1[1]\ar[d] \\
        & X_4 \ar@{=}[r]\ar[d] & X_4 \ar[d]\ar[r] & X_3[1]\\
        & X_3[1] \ar[r] & K[1] &
      }
    \]
By (3), $Y_3 \in \Y$. Then $f$ factors through an object in $\omega$, and then $f$ factors through $\alpha$. By Fact \ref{2factor}, $\Id_{Y_1}$ factors through $t$, and hence
$K$ is a direct summand of $X_3$. This proves $K \in \X$. \end{proof}

\section{\bf Model structure induced by a hereditary complete cotorsion pair}
Let $(\T, [1], \tri)$ be a triangulated category, $\xi$ a proper class of triangles, and $\mathcal{X}$ and $\mathcal{Y}$ additive full subcategories of $\T$ which are closed under isomorphisms and direct summands. Put $\omega:=\mathcal{X}\cap \mathcal{Y}$.
Define five classes of morphisms in $\T$ as follows.

\begin{equation} \begin{aligned}\label{construction}\CoFib_\omega^\xi & = \{\ \xi\text{-proper monomorphism} \ f \ | \ \hoc f\in \mathcal X\}\\
\Fib_\omega^\xi & = \{\text{morphism} \ f \ | \ \Hom_{\mathcal T}(W, f) \ \mbox{is surjective}, \ \forall \ W\in \omega\}\\
\Weq_\omega^\xi & = \{\ A\stackrel f \longrightarrow B \ | \ \mbox{there is a} \ \xi\mbox{-proper epimorphism} \  A\oplus W\stackrel {(f, t)} \longrightarrow B
\\ & \ \ \ \ \ \ \ \ \ \ \ \ \ \ \ \ \ \ \ \ \ \mbox{such that} \ W\in \omega \ \ \mbox{and} \ \ \hok (f, t)\in \mathcal{Y}\} \\
\TCoFib_\omega^\xi & = \CoFib_\omega^\xi \cap \Weq_\omega^\xi\\
\TFib_\omega^\xi & = \Fib_\omega^\xi \cap \Weq_\omega^\xi \end{aligned}\end{equation}

\vskip5pt
\noindent Thus, a morphism $f: A\longrightarrow B$ is in $\Weq_{\omega}^\xi$ if and only if there is a commutative diagram
\[
\xymatrix@R=0.5cm {A\ar[dr]_-{\left(\begin{smallmatrix}1\\ 0  \end{smallmatrix}\right)}\ar[rr]^-f && B\\ & A\oplus W\ar[ur]_-{(f, t)}}\]
such that $W\in \omega$, \ $(f, t)$ is a $\xi$-proper epimorphism,   and $\hok (f, t)\in \mathcal{Y}$.

\begin{thm}\label{if} \ Let $(\T, [1], \tri)$ be a triangulated category and $\xi$ a proper class of triangles.
Suppose that  $(\mathcal{X}, \mathcal{Y})$ is a hereditary complete cotorsion pair in $(\T, [1], \tri)$ with respect to $\xi$,
such that $\omega = \X \cap \Y$ is contravariantly finite in $\T$. Then $(\CoFib_{\omega}^\xi, \ \Fib_{\omega}^\xi, \ \Weq_{\omega}^\xi)$,  as defined in {\rm (\ref{construction})}, is a model structure on $\T$.

\vskip5pt

Moreover, the class $\TCoFib_\omega^\xi$ of trivial cofibrations is precisely
the class of splitting monomorphisms with hocokernel in $\omega$, and the class $\TFib_\omega^\xi$
of trivial fibrations is precisely the class of $\xi$-proper epimorphisms with hokernel in $\mathcal Y$. \end{thm}

\subsection{Descriptions of $\TCoFib_\omega$ and $\TFib_\omega$}
\begin{lem} \label{tcofibandtfib}
    \ Let $(\T, [1], \tri)$ be a triangulated category and $\xi$ a proper class of triangles.
Let $\X$ and $\Y$ be full additive subcategories of $\T$ which are closed under isomorphisms and direct summands, and $\omega = \X \cap \Y$. If $\Ext_{\xi}(\X, \Y) = 0$, then
\begin{align*}\TCoFib_\omega^\xi &  = \{\text{splitting monomorphism} \ f \ | \ \hoc f\in \omega\}\\
\TFib_\omega^\xi & =\{\ \xi\text{-proper epimorphism} \ f \ | \ \hok f\in \mathcal Y\}\\
\Weq_\omega^\xi &  = \{\ A\stackrel f \longrightarrow B \ | \ \exists \ \ A\oplus W\stackrel{(f, t)} \longrightarrow B \ \text{in} \ \TFib_\omega^\xi \
\mbox{with} \ W\in \omega\} \\ & = \TFib_\omega^\xi \circ \TCoFib_\omega^\xi.\end{align*}

\end{lem}
\begin{proof}  \ First, we prove $\TCoFib_\omega^\xi = \{\text{splitting monomorphism} \ f \ | \ \hoc f\in \omega\}$.
Let $f: A \xra{} B$ be a splitting monomorphism with $\hoc f\in \omega$.
By  definition $f\in \CoFib_{\omega}^\xi$. Without loss of generality one may assume that $f$ is just $f=\left(\begin{smallmatrix} 1\\0\end{smallmatrix}\right): A\longrightarrow A\oplus W$ with $W\in \omega$.
By the definition one sees $f\in \Weq_{\omega}^\xi$, by taking $t=\left(\begin{smallmatrix} 0\\1\end{smallmatrix}\right): W\longrightarrow A\oplus W$.

    \vskip5pt

Conversely, let $f: A\longrightarrow B$ be a morphism in $\TCoFib_\omega^\xi = \CoFib_\omega^\xi \cap \Weq_\omega^\xi$.
Then $f\in \CoFib_\omega^\xi$, i.e., $f$ is a $\xi$-proper monomorphism with $X:= \hoc f \in \X$.
Also, $f\in \Weq_\omega^\xi$, by definition there is a triangle $Y \xra{} A \oplus W \xra{{(f,t)}} B\xra{} Y[1]$ in $\xi$ with $W\in \omega$ and $Y\in \mathcal{Y}$.
Since $\Ext_{\xi}(X, Y) = 0$, it follows from Extension-Lifting Lemma \ref{extlifting} that there is a lifting $\left(\begin{smallmatrix}u\\ v \end{smallmatrix}\right): B \longrightarrow A \oplus W$ such that $\left(\begin{smallmatrix}u\\ v \end{smallmatrix}\right)f = \left(\begin{smallmatrix}1\\ 0 \end{smallmatrix}\right)$ and $(f,t) \left(\begin{smallmatrix}u\\ v \end{smallmatrix}\right) = {\rm Id}_B$.
    \[
    \xymatrix@=0.8cm{
     & A\ar[r]^-{f}\ar[d]_{\left(\begin{smallmatrix}1\\ 0 \end{smallmatrix}\right)} & B \ar[r]\ar@{=}[d]\ar@{-->}[dl]_-{\left(\begin{smallmatrix}u\\ v \end{smallmatrix}\right)}& X \ar[r] & A[1]\\
    Y\ar[r]& A \oplus W \ar[r]^-{(f,t)}& B \ar[r] & Y[1] & \\
    }
    \]
    Thus $f$ is a splitting monomorphism. By the Octahedral Axiom one has a commutative diagram
    \[
    \xymatrix@=0.6cm{
    A\ar[r]^-{\left(\begin{smallmatrix}1\\ 0 \end{smallmatrix}\right)}\ar@{=}[d] & A \oplus W\ar[r]^-{(0,1)}\ar[d]^-{(f,t)} & W\ar[r]^-{0}\ar[d] & A[1]\ar@{=}[d] \\
    A\ar[r]^-{f} & B\ar[r]\ar[d] & X\ar[r]\ar[d] & A[1]\ar[d]^-{\left(\begin{smallmatrix}1\\ 0 \end{smallmatrix}\right)} \\
      & Y[1]\ar@{=}[r]\ar[d] & Y[1]\ar[d]\ar[r] & A[1] \oplus W[1] \\
      & A[1]\oplus W[1]\ar[r] & W[1] & \\
    }
    \]
    Thus one has the following commutative diagram
    \[
    \xymatrix@=0.6cm{
      & A \ar@{=}[r]\ar[d]_-{\left(\begin{smallmatrix}1\\ 0 \end{smallmatrix}\right)} & A \ar[d]^-{f} & \\
    Y\ar[r]\ar@{=}[d] & A \oplus W\ar[d]_-{(0,1)}\ar[r]^-{(f,t)} & B\ar[r]\ar[d] & Y[1]\ar@{=}[d] \\
    Y\ar[r] & W\ar[r]\ar[d]_-{0} & X\ar[r]\ar[d] & Y[1]\ar[d] \\
      & A[1]\ar@{=}[r] & A[1]\ar[r] &  A[1] \oplus W[1]
    }
    \]
    Since the second horizontal triangle and the third vertical triangle belong to $\xi$,
it follows from the saturation  of $\xi$ that the triangle $Y \xra{} W \xra{} X\xra{} Y[1]$ also belongs to $\xi$.
This triangle splits, since $\Ext_{\xi}(X, Y) = 0$. Thus $X$ is a direct summand of $W\in \omega$. Since by assumption $\omega$ is closed under direct summands,
$X \in \omega$. Thus
     $f\in \{\text{splitting monomorphism} \ f \ | \ \hoc f\in \omega\}$.
This proves the equality.
    \vskip5pt

    Next, we prove $\TFib_\omega^\xi =\{\ \xi\text{-proper epimorphism} \ f \ | \ \hok f\in \mathcal Y\}$.
    Let $f: A \xra{} B$ be a $\xi$-proper epimorphism with $Y: = \hok f\in \mathcal Y$.
    Thus there is a triangle $Y \xra{} A \xra{f} B \xra{} Y[1]$ in $\xi$ with $Y \in \mathcal{Y}$. Clearly $f\in \Weq_{\omega}^\xi$, by taking
    $t: 0\longrightarrow B$. To show $f\in \Fib_{\omega}^\xi,$ one needs to show that $\Hom_{\mathcal A}(W, f)$ is surjective for any object $W\in \omega$.
In fact, for any morphism $u: W \xra{} B$, since $\Ext_{\xi}(W,Y) = 0$,
it follows from Extension-Lifting Lemma \ref{extlifting}  that there is a lifting $s: W \xra{} A$ such that $u = fs$.
    \[
    \xymatrix@=0.6cm{
     & 0\ar[r]\ar[d] & W \ar[d]^-{u}\ar@{=}[r]\ar@{-->}[dl]_-{s}& W \ar[r] & 0\\
    Y\ar[r]& A \ar[r]^-{f}& B \ar[r] & Y[1] & \\
    }
    \]
    This shows that $\Hom_{\mathcal A}(W, f)$ is surjective, and hence $f\in \Fib_\omega^\xi \cap \Weq_\omega^\xi = {\rm TFib}_{\omega}^\xi$.

    \vskip5pt

    Conversely, let $f: A \xra{} B$ be a morphism in ${\rm TFib}_{\omega}^\xi = \Fib_{\omega}^\xi\cap \Weq_{\omega}^\xi$.
   Since $f\in \Weq_{\omega}^\xi$, by definition there is a triangle
   $Y \xra{} A\oplus W \xra{(f,t)} B \xra{} Y[1]$ in $\xi$ with $W\in \omega$ and $Y\in \mathcal{Y}$.
   Since $f\in \Fib_{\omega}^\xi$, by definition $\Hom_{\mathcal A}(W, f)$ is surjective for any object $W\in \omega$.
   Thus there is a morphism $s: W\longrightarrow A$ such that $t=fs$.
   Then $(f,t)=f\circ (1,s)$. Since $(f, t)$ is a $\xi$-proper epimorphism, it follows from Lemma \ref{monicepimorphism}(2) that
   $f$ is a $\xi$-proper epimorphism. So one has a triangle $K \xra{k} A \xra{f} B \xra{} K[1]$ in $\xi$, and a commutative diagram
    \[
    \xymatrix@=0.6cm{
    K \ar[r]^-{k}\ar@{-->}[d]_-{g}& A\ar[r]^-{f}\ar[d]_-{\left(\begin{smallmatrix}1\\0 \end{smallmatrix}\right)}& B\ar@{=}[d]\ar[r] & K[1]\ar@{-->}[d]^-{g[1]}\\
    Y \ar[r]\ar@{-->}[d]_-{h}& A\oplus W\ar[r]^-{(f,t)}\ar[d]_{(1,s)}& B\ar@{=}[d] \ar[r] & Y[1]\ar@{-->}[d]^-{h[1]}\\
    K \ar[r]^-{k}& A\ar[r]^-{f}&B\ar[r] & K[1]\\
    }
    \]
    Then one gets a morphism of triangles
    \[
    \xymatrix@=0.6cm{
    K \ar[r]^-{k}\ar[d]_-{hg}& A\ar[r]^-{f}\ar@{=}[d]& B\ar@{=}[d]\ar[r] & K[1]\ar[d]^-{(hg)[1]}\\
    K \ar[r]^-{k}& A\ar[r]^-{f}&B\ar[r] & K[1]\\
    }
    \]
    which shows that $hg$ is an isomorphism. Thus $K$ is a direct summand of $Y \in \Y$, and  $K \in \Y$.
   This proves $f\in \{\ \xi\text{-proper epimorphism} \ f \ | \ \hok f\in \mathcal Y\}$, and hence the equality.

\vskip5pt

Finally, using $\TFib_\omega^\xi  =\{\ \xi\text{-proper epimorphism} \ f \ | \ \hok f\in \mathcal Y\}$ one can rewrite
$$\Weq_\omega^\xi  = \{\ A\stackrel f \longrightarrow B \ | \ \exists \ \mbox{a} \ \xi\mbox{-proper epimorphism} \  A\oplus W\stackrel {(f, t)} \longrightarrow B
 \ \mbox{with} \ W\in \omega,  \ \hok (f, t)\in \mathcal{Y}\}$$
as
$$\Weq_\omega^\xi  = \{\ A\stackrel f \longrightarrow B \ | \ \exists \ \mbox{a morphism} \ A\oplus W\stackrel{(f, t)} \longrightarrow B \ \text{in} \ \TFib_\omega^\xi \
\mbox{with} \ W\in \omega\}.$$
It is clear that $\Weq_\omega^\xi  = \TFib_\omega^\xi \circ \TCoFib_\omega^\xi$.  In fact, if $f: A\longrightarrow B$ is in $\Weq_{\omega}^\xi$,
then $f$ is the composition of
$A \xra{\left(\begin{smallmatrix}1\\0 \end{smallmatrix}\right)} A\oplus W \xra{(f,t)} B$,
where $\left(\begin{smallmatrix}1\\0 \end{smallmatrix}\right)$ is in $\TCoFib_\omega^\xi$
and $(f,t)$ is in $\TFib_\omega^\xi$. Thus $\Weq_\omega^\xi  \subseteq\TFib_\omega^\xi \circ \TCoFib_\omega^\xi$. Conversely, if $f = p\circ i: A\longrightarrow B$ with $i\in \TCoFib_\omega^\xi$ and $p\in \TFib_\omega^\xi$,
without loss of generality, one may assume that $i = \left(\begin{smallmatrix}1\\0 \end{smallmatrix}\right): A\longrightarrow A\oplus W$ with $W\in \omega$. Then $p = (f, t): A\oplus W \longrightarrow B$
is in  $\TFib_\omega^\xi$. By definition $f\in \Weq_\omega^\xi$.
\end{proof}

\subsection{Factorization Axiom}
\begin{lem}\label{cofibandtfibcomp} \ Keep the assumptions in {\rm Theorem \ref{if}}. Then the classes  $\CoFib_{\omega}^\xi$ and $\TFib_{\omega}^\xi$ are closed under compositions.
\end{lem}
\begin{proof} Let $\alpha: A \xra{} B$ and $\beta: B \xra{} C$ be morphisms in $\CoFib_{\omega}^\xi$.
Then one has two triangles $A \xra{\alpha} B \xra{} \hoc \alpha \xra{} A[1]$ and $B \xra{\beta} C \xra{} \hoc \beta \xra{} B[1]$ in $\xi$,
where $\hoc \alpha\in \X$ and $\hoc \beta \in \X$. By Lemma \ref{compclose}(1),
$\beta \alpha$ is a $\xi$-proper monomorphism. It remains to show that $\hoc \beta \alpha \in \X$.
By Lemma \ref{4tri}(2),  $\hoc \alpha \xra{} \hoc \beta \alpha \xra{} \hoc \beta \xra{} \hoc \alpha [1]$ is a triangle in $\xi$,
and then  by Corollary \ref{extclose} one knows that $\hoc \beta \alpha \in \X$. By definition $\beta\alpha\in \CoFib_{\omega}^\xi.$

    \vskip5pt
Similarly one can show that $\TFib_{\omega}^\xi$ is closed under compositions, by using Lemma \ref{4tri}(1).
\end{proof}

\begin{lem}\label{factor} \ \ Keep the assumptions in {\rm Theorem \ref{if}}.  Then each morphism $f: A \xra{} B$ in $\T$ admits factorizations $f = pi = qj$,
 where $i\in \TCoFib_{\omega}^\xi$ and $p \in \Fib_{\omega}^\xi$, and $j\in \CoFib_{\omega}^\xi$ and $q \in \TFib_{\omega}^\xi$.
\end{lem}

\begin{proof} \  Since $\omega$ is contravariantly finite in $\T$, there is a right $\omega$-approximation $\tau_B: W_B \longrightarrow B$ of $B$.
Then $(f, \tau_B): A\oplus W_B \xra{} B$ is in $\Fib_{\omega}^\xi$:  in fact, for each morphism $g: W' \xra{} B$ with $W' \in \omega$,
there is a morphism $h: W' \xra{} W_B$ such that $g = \tau_B h$; and then $g = (f,\tau_B) \left(\begin{smallmatrix} 0 \\ h \end{smallmatrix}\right)$
with $\left(\begin{smallmatrix} 0 \\ h \end{smallmatrix}\right): W' \xra{} A \oplus W_B$.
Thus one has the factorization $f = (f,\tau_B) \left(\begin{smallmatrix} 1 \\ 0 \end{smallmatrix}\right)$, where $\left(\begin{smallmatrix} 1 \\ 0 \end{smallmatrix}\right): A \xra{} A \oplus W_B$ is in $\TCoFib_{\omega}^\xi$ and $(f,\tau_B):A\oplus T_B \xra{} B$ is in $\Fib_{\omega}^\xi$. This gives the first factorization.

\vskip5pt

By the assumption that $(\mathcal{X}, \mathcal{Y})$ is a complete cotorsion pair in $\T$ with respect to $\xi$, one has
a triangle $Y_B \longrightarrow X_B \xra{t_B} B \longrightarrow Y_B[1]$ in $\xi$ with $X_B \in \X$ and $Y_B \in \Y$.
By Lemma \ref{carsq}(1), $(f,t_B): A \oplus X_B \xra{} B$ is also a $\xi$-proper epimorphism,
so one gets a triangle $K \xra{} A \oplus X_B \xra{(f,t_B)} B \xra{} K[1]$ in $\xi$.
Similarly,  there is a triangle $K \xra{\sigma} Y \xra{} X \xra{} K[1]$ in $\xi$ with $X \in \X$ and $Y \in \Y$. By the cobase change one has the following commutative diagram
\[
  \xymatrix@R=0.6cm{
    & X[-1] \ar@{=}[r]\ar[d] & X[-1] \ar[d] & \\
  B[-1]\ar[r]\ar@{=}[d] & K \ar[d]^-{\sigma}\ar[r] & A \oplus X_B \ar[r]^-{(f,t_B)}\ar[d]^-{i} & B\ar@{=}[d] \\
  B[-1]\ar[r] & Y \ar[r]\ar[d] & E\ar[r]^-{p}\ar[d] & B\ar[d] \\
    & X\ar@{=}[r] & X\ar[r] &  K[1]
  }
\]
Since $\xi$ is closed under cobase changes,
the triangle $Y \xra{} E \xra{p} B \xra{} Y[1]$ is in $\xi$. By Lemma \ref{tcofibandtfib}, \ $p \in \TFib_{\omega}^\xi$.
Transposing the rows and columns of the diagram, for the same reason one sees that
the triangle $A\oplus X_B \xra{i} E \xra{} X \xra{} A[1]\oplus X_B [1]$ is in $\xi$. By definition $i \in \CoFib_{\omega}^\xi$.
Thus
$$f = (f, t_B) \left(\begin{smallmatrix} 1 \\ 0 \end{smallmatrix}\right) = p \circ (i \circ \left(\begin{smallmatrix} 1 \\ 0 \end{smallmatrix}\right)),$$
where the composition $i \circ \left(\begin{smallmatrix} 1 \\ 0 \end{smallmatrix}\right): A \xra{\left(\begin{smallmatrix} 1 \\ 0 \end{smallmatrix}\right)} A \oplus X_B \xra{i} E$
is  in $\CoFib_{\omega}^\xi$,  by Lemma \ref{cofibandtfibcomp}(1). This gives the second factorization.
\end{proof}

\subsection{Two out of three axiom.}

The aim of this subsection is to prove that ${\rm Weq}_\omega^\xi$ satisfies Two out of three axiom. For this we need some preparations.

\vskip5pt

\begin{lem} \label{weqcomp} \ Keep the assumptions in {\rm Theorem \ref{if}}. Then ${\rm Weq}_\omega^\xi$ is closed under compositions.
\end{lem}

\begin{proof}
  \ Let $\alpha: A \longrightarrow B$ and $\beta: B\longrightarrow C$ be in ${\rm Weq}_\omega^\xi$. Then there is a morphism $(\alpha, t_1): A\oplus W_1 \longrightarrow B$ in
  ${\rm TFib}_\omega^\xi$ with $W_1\in \omega$, and there is a morphism $(\beta, t_2):  B\oplus W_2 \longrightarrow C$ in ${\rm TFib}_\omega^\xi$ with $W_2\in \omega$. Then $\beta\alpha$ has a factorization
  $$\beta\alpha =(\beta\alpha, \beta t_1 , t_2)\left(\begin{smallmatrix} 1\\0\\0\end{smallmatrix}\right)$$
where $\left(\begin{smallmatrix} 1\\0\\0\end{smallmatrix}\right): A \longrightarrow A\oplus W_1\oplus W_2$
is in ${\rm TCoFib}_\omega^\xi$,
and $(\beta\alpha, \beta t_1 , t_2): A\oplus W_1\oplus W_2\longrightarrow C$. See the following diagram.
  \[\xymatrix@R=0.6cm{
    A\ar[dr]_-{\left(\begin{smallmatrix} 1\\ 0 \end{smallmatrix}\right)}\ar[rr]^-{\alpha}& & B\ar[dr]_-{\left(\begin{smallmatrix}1\\ 0   \end{smallmatrix}\right)}\ar[rr]^-{\beta}&&C\\
    &A\oplus W_1\ar[dr]_(.4){\left(\begin{smallmatrix}1&0\\0& 1\\0&0   \end{smallmatrix}\right)}\ar[ur]_-{\left(\begin{smallmatrix} \alpha, t_1 \end{smallmatrix}\right)}&&B\oplus W_2\ar[ur]_-{\left(\begin{smallmatrix} \beta, t_2 \end{smallmatrix}\right)}&\\
    &&A\oplus W_1\oplus W_2\ar[ur]_(.65){\left(\begin{smallmatrix}\alpha& t_1 & 0\\ 0 &0 &1   \end{smallmatrix}\right)} &&
  }
  \]
By definition  $\hok (\alpha, t_1) \xra{} A \oplus W_1 \xra{(\alpha, t_1)} B \xra{} \hok (\alpha, t_1) [1]$ is a triangle in $\xi$.
Since $\hok (\alpha, t_1) \xra{} A \oplus W_1 \oplus W_2 \xra{\left(\begin{smallmatrix}\alpha& t_1 & 0\\ 0 &0 &1   \end{smallmatrix}\right)} B \oplus W_2 \xra{} \hok (\alpha, t_1) [1]$  is the coproduct of the triangle $\hok (\alpha, t_1) \xra{} A \oplus W_1 \xra{(\alpha, t_1)} B \xra{} \hok (\alpha, t_1) [1]$
and $0 \xra{} W_2 \xra{1} W_2 \xra{} 0$, it is in $\xi$. So $\left(\begin{smallmatrix}\alpha& t_1 & 0\\ 0 &0 &1   \end{smallmatrix}\right)\in \TFib_{\omega}^\xi$, and then $(\beta\alpha, \beta t_1 , t_2) = (\beta, t_2)\left(\begin{smallmatrix}\alpha& t_1 & 0\\ 0 &0 &1   \end{smallmatrix}\right)\in \TFib_\omega^\xi$ by Lemma \ref{cofibandtfibcomp}. Hence $\beta\alpha=(\beta\alpha, \beta t_2 , t_1)\left(\begin{smallmatrix} 1\\0\\0\end{smallmatrix}\right)\in \TFib_\omega^\xi\circ {\rm TCoFib}_\omega^\xi = {\rm Weq}_\omega^\xi$.
\end{proof}

\begin{lem}\label{necessityofWeq} \ Keep the assumptions in {\rm Theorem \ref{if}}. Let $\alpha: A\longrightarrow B$ be a morphism in ${\rm Weq}_\omega^\xi.$  Then for an arbitrary right $\omega$-approximation $t: W\longrightarrow B$ of $B$, the morphism \ $(\alpha, t): A\oplus W\longrightarrow B$ is in ${\rm TFib}_\omega^\xi$.
\end{lem}
\begin{proof} \ Since $\alpha$ is in ${\rm Weq}_\omega^\xi$, there is a triangle $Y \xra{} A \oplus W' \xra{(\alpha, t')} B \xra{} Y[1]$ in $\xi$ with $W'\in \omega$ and $Y \in \Y$.
Since $t: W\longrightarrow B$ is a right $\omega$-approximation,  there is a morphism $s: W'\longrightarrow W$ such that $t'=ts$.
Since  $(\alpha, t')$ is a $\xi$-proper epimorphism and $(\alpha, t')=(\alpha, t)\left(\begin{smallmatrix} 1 & 0\\0 & s\end{smallmatrix}\right)$, it follows from Lemma \ref{monicepimorphism}(2) that $(\alpha, t): A\oplus W \rightarrow B$ is also a $\xi$-proper epimorphism. It remains to prove that $\hok (\alpha, t)\in\Y$.

  \vskip5pt

  Since $\left(\begin{smallmatrix} 1 & 0\\0 & 1\\ 0 & 0 \end{smallmatrix}\right): A\oplus W'\longrightarrow A\oplus W'\oplus W$ is a splitting monomorphism with $\hoc \left(\begin{smallmatrix} 1 & 0\\0 & 1\\ 0 & 0 \end{smallmatrix}\right) = W \in \omega$, it is in $\TCoFib_{\omega}^\xi$ and is a $\xi$-proper monomorphism. Now that $(\alpha, t')=(\alpha, t', t)\left(\begin{smallmatrix} 1 & 0\\0 & 1\\ 0 & 0 \end{smallmatrix}\right)$ is a $\xi$-proper epimorphism, it follows from Lemma \ref{4tri}(3) that there is a triangle in $\xi$
  $$Y \longrightarrow \hok (\alpha, t', t) \longrightarrow W \longrightarrow Y[1].$$
  Since $Y \in \Y$ and $W \in \Y$, $\hok (\alpha, t', t) \in \Y$, by Corollary \ref{extclose}. Note that
  $$W' \xra{\left(\begin{smallmatrix}  0\\ 1\\ -s \end{smallmatrix}\right)} A \oplus W' \oplus W \xra{\left(\begin{smallmatrix} 1 & 0 & 0\\0 & s & 1 \end{smallmatrix}\right)} A \oplus W \longrightarrow W'[1]$$
  is a splitting triangle,  thus it is in $\xi$, so $\left(\begin{smallmatrix} 1 & 0 & 0\\0 & s & 1 \end{smallmatrix}\right)$ is a $\xi$-proper epimorphism with $\hok \left(\begin{smallmatrix} 1 & 0 & 0\\0 & s & 1 \end{smallmatrix}\right) = W' \in \omega$. Hence by $(\alpha, t', t) = (\alpha, t)\left(\begin{smallmatrix}1&0&0\\0&s&1 \end{smallmatrix}\right)$
and Lemma \ref{4tri}(1) there is a triangle in $\xi$
  $$W' \xra{} \hok (\alpha, t', t) \xra{} \hok (\alpha, t) \xra{} W'[1].$$
  By the assumption that $(\X, \Y)$ is a hereditary cotorsion pair, one gets $\hok (\alpha, t) \in \Y$. This completes the proof.
\end{proof}

\begin{lem}\label{TFibWeq} \ Keep the assumptions in {\rm Theorem \ref{if}}. Let $\alpha: A\longrightarrow B$ and $\beta: B\longrightarrow C$ be morphisms
in $\T$ with $\alpha\in {\rm TFib}_\omega^\xi$ and $\beta\alpha\in {\rm Weq}_\omega^\xi$. Then $\beta \in {\rm Weq}_\omega^\xi$.
\end{lem}
\begin{proof}
  \ Take a right $\omega$-approximation $t: W\longrightarrow C$ of $C$. Since $\beta\alpha\in {\rm Weq}_\omega^\xi$, it follows from Lemma \ref{necessityofWeq} that $(\beta\alpha, t): A\oplus W \longrightarrow C$ is in ${\rm TFib}_\omega^\xi$, i.e., $(\beta\alpha, t)$ is a $\xi$-proper epimorphism with $\hok (\beta\alpha, t)\in \Y$. Since $(\beta\alpha, t)=(\beta, t)\left(\begin{smallmatrix} \alpha & 0\\ 0 & 1\end{smallmatrix}\right)$, it follows from Lemma \ref{monicepimorphism}(2) that $(\beta, t): B\oplus W \longrightarrow C$ is a $\xi$-proper epimorphism.

  \vskip5pt

Since the triangle $\hok \alpha \xra{} A \oplus W \xra{\left(\begin{smallmatrix} \alpha & 0\\ 0 & 1\end{smallmatrix}\right)} B\oplus W \xra{} \hok \alpha [1]$ is
the coproduct of triangles $\hok \alpha \xra{} A  \xra{\alpha} B \xra{} \hok \alpha [1]$ and $0 \xra{} W \xra{1} W \xra{} 0$ in $\xi$, it is in $\xi$,
i.e., $\left(\begin{smallmatrix} \alpha & 0\\ 0 & 1\end{smallmatrix}\right)$ is a $\xi$-proper monomorphism.
By $(\beta\alpha, t)=(\beta, t)\left(\begin{smallmatrix} \alpha & 0\\ 0 & 1\end{smallmatrix}\right)$
and Lemma \ref{4tri}(1) one gets a triangle in $\xi$:
  $$\hok \alpha \xra{} \hok (\beta\alpha, t) \xra{} \hok (\beta, t) \xra{} \hok \alpha [1].$$
Note that $\hok \alpha \in \Y$ and $\hok (\beta\alpha, t) \in \Y$. Since $(\X, \Y)$ is a hereditary cotorsion pair, $\hok (\beta, t) \in \Y$.  By definition $\beta\in {\rm Weq}_\omega^\xi$.
\end{proof}

\begin{lem} \label{231}  \ Keep the assumptions in {\rm Theorem \ref{if}}. Let $\alpha: A\longrightarrow B$ and $\beta: B\longrightarrow C$ be morphisms in $\T$ such that $\alpha$ and $\beta\alpha$ are in ${\rm Weq}_\omega^\xi$. Then $\beta \in {\rm Weq}_\omega^\xi$.
\end{lem}

\begin{proof}
  \ Since $\alpha\in {\rm Weq}_\omega^\xi$, there is a morphism $(\alpha, t): A\oplus W\longrightarrow B$
  with $W\in \omega$ such that $(\alpha, t)\in {\rm TFib}_\omega^\xi$.
  To prove $\beta \in {\rm Weq}_\omega^\xi$, by the factorization
  \[
  \xymatrix@R=0.5cm {A\oplus W\ar[dr]_-{(\alpha, t)}\ar[rr]^-{(\beta\alpha, \beta t)} & & C\\
   & B\ar[ur]_-{\beta} & }
  \]
  and Lemma \ref{TFibWeq}, it suffices to prove $(\beta\alpha, \beta t)\in {\rm Weq}_\omega^\xi$.

  \vskip5pt

  Take a right $\omega$-approximation $t': W'\longrightarrow C$ of $C$. Since $\beta\alpha\in {\rm Weq}_\omega^\xi$,
  it follows from Lemma \ref{necessityofWeq} that $(\beta\alpha, t'): A\oplus W'\longrightarrow C$ is in ${\rm TFib}_\omega^\xi$, i.e.,
  $(\beta\alpha, t')$ is a $\xi$-proper epimorphism and $\hok (\beta\alpha, t')\in \Y$.

 \vskip5pt
 Consider the factorization \[
  \xymatrix@R=0.5cm {A\oplus W'\ar[dr]_-{\left(\begin{smallmatrix}1&0\\0&1\\0&0 \end{smallmatrix}\right)}\ar[rr]^-{(\beta\alpha, t')} && C\\ & A\oplus W\oplus W'\ar[ur]_-{(\beta\alpha, \beta t, t')}}\]
 Since $(\beta\alpha, t') = (\beta\alpha, \beta t, t')\left(\begin{smallmatrix}1&0\\0&0\\0&1 \end{smallmatrix}\right)$ is a $\xi$-proper epimorphism, it follows that
 $(\beta\alpha, \beta t, t'): A\oplus W\oplus W'\longrightarrow C'$ is a $\xi$-proper epimorphism (cf. Lemma \ref{monicepimorphism}).
 Since the splitting monomorphism $\left(\begin{smallmatrix}1&0\\0&0\\0&1 \end{smallmatrix}\right): A\oplus W'\longrightarrow A\oplus W\oplus W'$
  is a $\xi$-proper monomorphism with $\hoc \left(\begin{smallmatrix}1&0\\0&0\\0&1 \end{smallmatrix}\right) = W \in \omega$ and since $(\beta\alpha, t') = (\beta\alpha, \beta t, t')\left(\begin{smallmatrix}1&0\\0&0\\0&1 \end{smallmatrix}\right)$ is a $\xi$-proper epimorphism, it follows from Lemma \ref{4tri}(3) that there is a triangle in $\xi$
  $$\hok (\beta\alpha, t') \longrightarrow \hok (\beta\alpha, \beta t, t') \longrightarrow
  W \longrightarrow \hok (\beta\alpha, t') [1].$$
  Since $\hok (\beta\alpha, t')\in\Y$ and $W\in \mathcal{Y}$, it follows from Corollary \ref{extclose} that $\hok (\beta\alpha, \beta t, t')\in Y$,
  and hence $(\beta\alpha, \beta t, t')\in {\rm TFib}_\omega^\xi.$
  By the factorization
  \[
  \xymatrix@R=0.5cm {A\oplus W\ar[dr]_-{\left(\begin{smallmatrix}1&0\\0&1\\0&0 \end{smallmatrix}\right)}\ar[rr]^-{(\beta\alpha, \beta t)} && C\\ & A\oplus W\oplus W'\ar[ur]_-{(\beta\alpha, \beta t, t')}}\]
 and $\left(\begin{smallmatrix}1&0\\0&0\\0&1 \end{smallmatrix}\right)\in\TCoFib_\omega^\xi$,
 we see that $(\beta\alpha, \beta t)\in \TFib_\omega^\xi\circ \TCoFib_\omega^\xi  = {\rm Weq}_\omega^\xi$. This completes the proof.
  \end{proof}

  \begin{lem}\label{CofTFib} \ Keep the assumptions in {\rm Theorem \ref{if}}. Let $\alpha=pi: A \xra{} B$ be a morphism in $\Weq_{\omega}^\xi$, where $i: A \xra{} C$ is in $\CoFib_{\omega}^\xi$ and $p: C \xra{} B$ is in ${\rm TFib}_{\omega}^\xi$. Then $i\in {\rm TCoFib}_{\omega}^\xi$.
  \end{lem}
  \begin{proof}
  \ We first show that $i$ is a splitting monomorphism. Since $i : A \xra{} C \in \CoFib_{\omega}^\xi$, $i$ is a $\xi$-proper monomorphism with $\hoc i\in \mathcal X$. Since $\alpha\in \Weq_{\omega}^\xi$, by definition there is a $\xi$-proper epimorphism $(\alpha, t): A\oplus W\longrightarrow B$ with $W \in \omega$ and $\hok (\alpha, t)\in \Y$.
  Consider the following commutative diagram
  \[
  \xymatrix{
    & A \ar[r]^-{i}\ar[d]_-{\left( \begin{smallmatrix} 1\\ 0  \end{smallmatrix} \right)} & C \ar[d]^-{p}\ar[r]\ar@{-->}[dl]_-{\left( \begin{smallmatrix} \sigma_1 \\ \sigma_2 \end{smallmatrix} \right)}& \hoc i \ar[r] & A[1]\\
    \hok (\alpha, t) \ar[r]& A \oplus W \ar[r]^-{(\alpha, t)}& B \ar[r] & \hok (\alpha, t)[1] & \\
  }
  \]
  By Extension-Lifting Lemma \ref{extlifting}, there is a lifting $\left( \begin{smallmatrix} \sigma_1 \\ \sigma_2 \end{smallmatrix} \right): C\rightarrow A\oplus W$ such that $\sigma_1 i={\rm Id}_A$. Thus $i$ is a splitting monomorphism.

  \vskip5pt

  Thus, without loss of generality, one may write $\alpha=pi$ as
  \[
  \xymatrix@R=0.5cm {A\ar[dr]_-{i=\left(\begin{smallmatrix}1\\ 0  \end{smallmatrix}\right)}\ar[rr]^-\alpha && B\\ & A\oplus X\ar[ur]_-{p=(\alpha, \alpha')}}
  \]
  with $X = \hoc i\in \mathcal X$ and $p=(\alpha, \alpha')\in {\rm TFib}_\omega^\xi$.
  It remains to prove that $X \in \mathcal{Y}$.

  \vskip5pt

  Since $\left(\begin{smallmatrix}1&0\\0&0\\0&1 \end{smallmatrix}\right) : A \oplus X \xra{} A \oplus W \oplus X$ is a splitting monomorphism with $\hoc \left(\begin{smallmatrix}1&0\\0&0\\0&1 \end{smallmatrix}\right) = W \in \omega$ and $p= (\alpha, \alpha') = (\alpha, t, \alpha')\left(\begin{smallmatrix}1&0\\0&0\\0&1 \end{smallmatrix}\right)$ is a $\xi$-proper epimorphism, where $(\alpha, t, \alpha'): A \oplus W \oplus X \xra{} B$,
  it follows that $(\alpha, t, \alpha')$ is a $\xi$-proper epimorphism and  that  there is a triangle in $\xi$ (cf. Lemma \ref{4tri}(3))
  $$\hok (\alpha, \alpha') \longrightarrow \hok (\alpha, \alpha', t) \longrightarrow
  W \longrightarrow \hok (\alpha, \alpha')[1].$$
  Since $\hok (\alpha, \alpha') \in \Y$ and $W\in \mathcal{Y}$, it follows from Lemma \ref{extclose} that $\hok (\alpha, \alpha', t)\in \Y$.

  \vskip5pt

  Since $(\alpha, t, \alpha')$ is a $\xi$-proper epimorphism, there is a triangle in $\xi$
  $$\hok (\alpha, t, \alpha') \xra{\left(\begin{smallmatrix} k_1 \\ k_2 \\ k_3 \end{smallmatrix}\right)} A \oplus W \oplus X \xra{(\alpha, t, \alpha')} B \xra{} \hok (\alpha, t, \alpha')[1].$$
  By Lemma \ref{addaxiom}(4) one has a commutative diagram
  \[
    \xymatrix@R=0.6cm{
    \hok (\alpha, t) \ar[r]\ar@{=}[d] & \hok (\alpha, t, \alpha') \ar[r]^-{-k_3}\ar[d]_-{\left(\begin{smallmatrix} k_1 \\ k_2 \end{smallmatrix}\right)} & X \ar[r]\ar[d]^-{\alpha '} & \hok (\alpha, t)[1]\ar@{=}[d] \\
    \hok (\alpha, t) \ar[r] & A \oplus W \ar[r]^-{(\alpha, t)}\ar[d] & B \ar[r]\ar[d] & \hok (\alpha, t)[1]\ar[d] \\
      & U \ar@{=}[r]\ar[d] & U \ar[d]\ar[r] &  \hok (\alpha, t, \alpha') [1]\\
      & \hok (\alpha, t, \alpha')[1]\ar[r] & X[1] &
    }
  \]
  By base change one gets a commutative diagram
  \[
  \xymatrix@R=0.6cm{
    & U[-1] \ar@{=}[r]\ar[d] & U[-1] \ar[d] & \\
    \hok (\alpha, t)\ar[r]\ar@{=}[d] & \hok (\alpha, t, \alpha') \ar[d]_-{\left(\begin{smallmatrix} k_1 \\ k_2 \end{smallmatrix}\right)}\ar[r]^-{-k_3} & X \ar[r]\ar[d]^-{\alpha '} & \hok (\alpha, t)[1]\ar@{=}[d] \\
    \hok (\alpha, t)\ar[r] & A\oplus W \ar[r]^-{(\alpha, t)}\ar[d] & B \ar[r]\ar[d] & \hok (\alpha, t)[1]\ar[d] \\
    & U\ar@{=}[r] & U\ar[r] &  \hok (\alpha, t, \alpha')[1]
      }
  \]
  where  the third horizontal triangle is in $\xi$. Since $\xi$ is closed under the base changes,
  the second horizontal triangle is also in $\xi$. Since $(\mathcal{X}, \Y)$ is a hereditary cotorsion pair with $\hok (\alpha, t) \in \Y$ and
$\hok (\alpha, t, \alpha')\in \mathcal{Y}$, it follows that  $X\in\Y$. This completes the proof.
  \end{proof}

\begin{lem} \label{232} \ Keep the assumptions in {\rm Theorem \ref{if}}. Let $\alpha: A\longrightarrow B$ and $\beta: B\longrightarrow C$ be morphisms in $\T$ with  $\beta\in {\rm Weq}_\omega^\xi$ and $\beta\alpha\in {\rm Weq}_\omega^\xi$. Then $\alpha \in {\rm Weq}_\omega^\xi$.
\end{lem}
\begin{proof}
  \ Since $\beta\in {\rm Weq}_\omega^\xi$, there is a morphism $(\beta, t): B \oplus W \xra{} C$ in $\TFib_{\omega}^\xi$ with $W\in \omega$. By Factorization axiom, which has been already proved, one can factorize  $\left(\begin{smallmatrix} \alpha \\ 0\end{smallmatrix}\right): A\longrightarrow B\oplus W$ as $\left(\begin{smallmatrix} \alpha \\ 0\end{smallmatrix}\right)=\left(\begin{smallmatrix} p_1 \\ p_2\end{smallmatrix}\right) i$ with  $i\in \CoFib_{\omega}^\xi$ and $\left(\begin{smallmatrix} p_1 \\ p_2 \end{smallmatrix}\right)\in {\rm TFib}_{\omega}^\xi$.
  By Lemma \ref{cofibandtfibcomp}, $(\beta, t)\left(\begin{smallmatrix} p_1 \\ p_2 \end{smallmatrix}\right)\in {\rm TFib}_{\omega}^\xi.$ Write
  $$\beta\alpha = (\beta, t)\left(\begin{smallmatrix} \alpha\\ 0 \end{smallmatrix}\right) = (\beta, t)\left(\begin{smallmatrix} p_1 \\ p_2 \end{smallmatrix}\right) i$$ where $\beta\alpha\in \Weq_{\omega}^\xi,$ $i\in \CoFib_{\omega}^\xi$ and $(\beta, t)\left(\begin{smallmatrix} p_1 \\ p_2 \end{smallmatrix}\right)\in {\rm TFib}_{\omega}^\xi.$
  Then $i\in {\rm TCofib}_{\omega}^\xi$, by Lemma \ref{CofTFib}. It follows that $\left(\begin{smallmatrix} \alpha \\ 0\end{smallmatrix}\right)\in {\rm TFib}_{\omega}^\xi\circ {\rm TCofib}_{\omega}^\xi=\Weq_{\omega}^\xi$.
  Since  $(1,0): B\oplus W \longrightarrow B$ is in ${\rm TFib}_{\omega}^\xi \subseteq \Weq_{\omega}^\xi$, it follows from Lemma \ref{weqcomp}
  that $\alpha=(1,0)\left(\begin{smallmatrix} \alpha \\ 0\end{smallmatrix}\right)\in \Weq_{\omega}^\xi$.
\end{proof}

  \vskip5pt

By Lemma \ref{weqcomp}, Lemma \ref{231} and Lemma \ref{232} one gets

\begin{lem}\label{2outof3} {\bf (Two out of three axiom)} \ Keep the assumptions in {\rm Theorem \ref{if}}. Let $\alpha: A \longrightarrow B$ and $\beta: B\longrightarrow C$ be morphisms in $\mathcal A$.
  If two of the three morphisms $\alpha, \ \beta, \ \beta\alpha$ are in ${\rm Weq}_\omega^\xi$, then so is the third.
\end{lem}

\subsection{Retract axiom.} Suppose that $g: A' \longrightarrow B'$ is a retract of $f: A \longrightarrow B$, i.e., there is a commutative diagram of morphisms
$$\xymatrix@R=0.5cm{
    A'\ar[r]^{\varphi_1}\ar[d]_g & A\ar[r]^{\psi_1}\ar[d]_{f} & A'\ar[d]_{g} \\
    B'\ar[r]^{\varphi_2} & B\ar[r]^{\psi_2} & B'
}$$
with  $\psi_1 \varphi_1 =\Id_{A'}$ and $\psi_2 \varphi_2 = \Id_{B'}$.
Then $\varphi_1$ and $\varphi_2$ are splitting monomorphisms which are $\xi$-proper monomorphisms,
$\psi_1$ and $\psi_2$ are splitting epimorphisms which are $\xi$-proper epimorphisms, and $A'$ and $B'$ are respectively direct summands of $A$ and $B$.

\begin{lem}\label{retract}
    The classes $\CoFib_{\omega}^\xi, \ \Fib_{\omega}^\xi$ and $\Weq_{\omega}^\xi$ are closed under retracts.
\end{lem}

\begin{proof} First,  $\CoFib_{\omega}^\xi$ is closed under retracts. In fact, let $f \in \CoFib_{\omega}^\xi$,
i.e., $f$ is a $\xi$-proper monomorphism with $\hoc f \in \mathcal{X}$.
Since $\varphi_2 g = f \varphi_1$ is a $\xi$-proper monomorphism (cf. Lemma \ref{compclose}), it follows from Lemma \ref{monicepimorphism} that $g$ is a $\xi$-proper monomorphism.
Then one has a commutative diagram of triangles in $\xi$
\[
\xymatrix{
A'\ar[r]^-{g}\ar@<-.5ex>[d]_-{\varphi_1} & B'\ar[r]^-{c_g}\ar@<-.5ex>[d]_-{\varphi_2} & \hoc g\ar[r]\ar@<-.5ex>@{-->}[d]_-{\widetilde{\varphi}} & A'[1]\ar@<-.5ex>[d]_-{\varphi_1 [1]} \\
A\ar[r]^-{f}\ar@<-.5ex>[u]_-{\psi_1} & B\ar[r]^-{c_f}\ar@<-.5ex>[u]_-{\psi_2} & \hoc f\ar[r]\ar@<-.5ex>@{-->}[u]_-{\widetilde{\psi}} & A[1] \ar@<-.5ex>[u]_-{\psi_1 [1]}
}
\]
Thus $\widetilde{\psi}\widetilde{\varphi} c_g = \widetilde{\psi} c_f \varphi_2 = c_g \psi_2 \varphi_2 = c_g$, so one has a morphism of triangles
\[
\xymatrix{
A'\ar[r]^-{g}\ar@{=}[d] & B'\ar[r]^-{c_g}\ar@{=}[d] & \hoc g\ar[r]\ar[d]_-{\widetilde{\psi} \widetilde{\varphi}} & A'[1]\ar@{=}[d] \\
A'\ar[r]^-{g} & B'\ar[r]^-{c_g} & \hoc g\ar[r] & A'[1]
}
\]
So $\widetilde{\psi} \widetilde{\varphi}$ is an isomorphism, and hence $\hoc g$ is a direct summand of $\hoc f \in \X$. Thus $\hoc g \in \X$ and $g \in \CoFib_{\omega}^\xi$.

\vskip5pt

Next, $\Fib_{\omega}^\xi$ is closed under retracts. In fact, let $f \in \Fib_{\omega}^\xi$, i.e., $\Hom_\mathcal T(W, f)$ is surjective for any $W\in \omega$.
It is clear that $\Hom_\mathcal T(W, g)$ is surjective for any $W\in \omega$, i.e., $g\in \Fib_{\omega}^\xi$.

\vskip5pt

Finally, we prove that $\Weq_{\omega}^\xi$ is closed under retracts. Suppose that $f \in \Weq_{\omega}^\xi$,
i.e., there is a triangle $Y \xra{} A \oplus W \xra{(f,\alpha)} B \xra{} Y[1]$ in $\xi$ with $Y \in \Y$ and $W \in \omega$.
Since $(g, \psi_2 \alpha) \left(\begin{smallmatrix} \psi_1 & 0 \\ 0 & 1 \end{smallmatrix}\right) = \psi_2 (f,\alpha)$ is a $\xi$-proper epimorphism (cf. Lemma \ref{compclose}),
$(g, \psi_2 \alpha): A'\oplus W\longrightarrow B'$ is a $\xi$-proper epimorphism. So there is a triangle $K \xra{\left(\begin{smallmatrix} k_1 \\ k_2 \end{smallmatrix}\right)} A' \oplus W \xra{(g, \psi_2 \alpha)} B' \xra{} K[1]$ in $\xi$. To show that $g \in \Weq_{\omega}^\xi$, it suffices to show that $K \in \Y$.

\vskip5pt

Consider a splitting triangle $B' \stackrel {\varphi_2}\longrightarrow B \stackrel {\partial_2}\longrightarrow B'' \stackrel {0}\longrightarrow B'[1]$. By Definition - Fact \ref{split}(v)
there is a morphism $\delta_2: B''\longrightarrow B$ such that  $\varphi_2 \psi_2 + \delta_2 \partial_2 = \Id_B$ and $\psi_2\delta_2 = 0$.

\vskip5pt

Since $(\X, \Y)$ is a complete cotorsion pair, there is a triangle $K \xra{i} Y_K \xra{d} X_K \xra{} K[1]$ in $\xi$ with $X_K \in \X$ and $Y_K \in \Y$.
Since $W \in \omega \subseteq \mathcal{Y}$, by Extension-Lifting Lemma \ref{extlifting} there exists a morphism $s: Y' \longrightarrow W$ such that $k_2 = si$. See the diagram
\[
\xymatrix@R=0.6cm{
  & K\ar[r]^-{i}\ar[d]_-{k_2} & Y_K \ar[d]\ar[r]^-{d}\ar@{-->}[dl]_-{s}& X_K \ar[r] & K[1]\\
  W\ar@{=}[r]& W \ar[r]& 0 \ar[r] & W[1]& \\
}
\]
Consider the following diagram
\[
\xymatrix{
  & K\ar[r]^-{i}\ar[d]_-{\left(\begin{smallmatrix} \varphi_1 k_1 \\ k_2 \end{smallmatrix}\right)} & Y_K \ar[d]^-{\delta_2 \partial_2 \alpha s}\ar[r]^-{d}\ar@{-->}[dl]_-{\left(\begin{smallmatrix} m \\ n \end{smallmatrix}\right)}& X_K \ar[r] & K[1]\\
  Y\ar[r]& A \oplus W \ar[r]^-{(f,\alpha)}& B \ar[r] & Y[1] & \\
}
\]
Since
\begin{equation*}
  \begin{aligned}
      (f,\alpha) \left(\begin{smallmatrix} \varphi_1 k_1 \\ k_2 \end{smallmatrix}\right) - \delta_2\partial_2\alpha si & = f \varphi_1 k_1 + \alpha k_2 - \delta_2\partial_2\alpha k_2 \\
       & = \varphi_2 g k_1 + \alpha k_2 - (\Id_B-\varphi_2 \psi_2)\alpha k_2 \\ & = \varphi_2 g k_1 + \varphi_2 \psi_2 \alpha k_2 = \varphi_2 (g, \psi_2 \alpha) \left(\begin{smallmatrix} k_1 \\ k_2 \end{smallmatrix}\right) = 0
  \end{aligned}
\end{equation*}
it follows from Extension-Lifting Lemma \ref{extlifting} that there exists a lifting $\left(\begin{smallmatrix} m  \\ n \end{smallmatrix}\right): Y_K \longrightarrow A \oplus W$ such that $fm + \alpha n  = \delta_2\partial_2\alpha s, \  mi = \varphi_1 k_1, \ ni = k_2$. Consider the following diagram
\[
\xymatrix{
Y_K \ar@{=}[r] & Y_K \ar[r]\ar[d]^-{\left(\begin{smallmatrix} \psi_1 m  \\ n \end{smallmatrix}\right)} & 0 \ar[r]\ar[d] & Y_K[1] \\
K \ar[r]^-{\left(\begin{smallmatrix} k_1 \\ k_2 \end{smallmatrix}\right)} & A' \oplus W \ar[r]^-{(g,\psi_2 \alpha)} & B' \ar[r] & K[1] \\
}
\]
Since
$$(g,\psi_2\alpha)\left(\begin{smallmatrix} \psi_1 m  \\ n \end{smallmatrix}\right)  = g\psi_1 m +\psi_2\alpha n = \psi_2 fm +\psi_2 \alpha n = \psi_2 \delta_2 \partial_2 \alpha s = 0$$
there exists a morphism $t: Y_K \longrightarrow K$ such that $\left(\begin{smallmatrix} \psi_1 m  \\ n \end{smallmatrix}\right) = \left(\begin{smallmatrix} k_1 \\ k_2 \end{smallmatrix}\right) t$. Then
$$\left(\begin{smallmatrix} k_1 \\ k_2 \end{smallmatrix}\right) ti = \left(\begin{smallmatrix} \psi_1 mi \\ ni \end{smallmatrix}\right) = \left(\begin{smallmatrix} \psi_1 \varphi_1 k_1 \\ ni \end{smallmatrix}\right) = \left(\begin{smallmatrix} k_1 \\ k_2 \end{smallmatrix}\right) .$$
So one has a morphism of triangles
\[
\xymatrix{
K \ar[r]^-{\left(\begin{smallmatrix} k_1 \\ k_2 \end{smallmatrix}\right)}\ar[d]_-{ti} & A' \oplus W \ar[r]^-{(g,\psi_2 \alpha)}\ar@{=}[d] & B' \ar[r]\ar@{=}[d] & K[1]\ar[d]^-{(ti)[1]} \\
K \ar[r]^-{\left(\begin{smallmatrix} k_1 \\ k_2 \end{smallmatrix}\right)} & A' \oplus W \ar[r]^-{(g,\psi_2 \alpha)} & B' \ar[r] & K[1] \\
}
\]
Thus $ti$ is an isomorphism and then $i: K\longrightarrow Y_K$ is a splitting monomorphism. Therefore $K$ is a direct summand of $Y_K \in \Y$,  and hence $K \in \Y$. This completes the proof.
\end{proof}

\subsection{Proof of Theorem \ref{if}} \ To prove that $({\rm CoFib}_{\omega}^\xi, \ {\rm Fib}_{\omega}^\xi, \ {\rm Weq}_{\omega}^\xi)$
is a model structure on $\T$, by Lemma \ref{factor}, \ref{2outof3} and \ref{retract} it remains to prove Lifting axiom.
Let
\[
\xymatrix@R=0.5cm{
A\ar[r]^-{f}\ar[d]_-{i} & C\ar[d]^-{p}\\
B\ar[r]^-{g} & D
}
\]
be a commutative square with $i\in \CoFib_\omega^\xi$ and $p\in \Fib_\omega^\xi$.

\vskip5pt

{\bf Case 1. \ Suppose that $p \in \TFib_{\omega}^\xi$}. By Lemma \ref{tcofibandtfib}, $p$ is a $\xi$-proper epimorphism with $\hok p \in \Y$.
Then the lifting $\lambda: B\longrightarrow C$ indeed exists, directly by Extension-Lifting Lemma \ref{extlifting}.

\vskip5pt

{\bf Case 2. \ Suppose that $i \in \TCoFib_{\omega}^\xi$}.
By Lemma \ref{tcofibandtfib}, $i$ is a splitting monomorphism with $\hoc i \in \omega$. Without loss of generality, one can rewrite the commutative square as
\[
\xymatrix{
A\ar[r]^-{f}\ar[d]_-{\left(\begin{smallmatrix}1\\0\end{smallmatrix}\right)} & C\ar[d]^-{p}\\
A\oplus W\ar[r]^-{(pf,g')} & D
}
\]
where $W=\hoc i\in \omega$. By definition $\Hom_\mathcal T(W, p): \Hom_\mathcal T(W, C)\longrightarrow \Hom_\mathcal T(W, D)$ is surjective, therefore there is a morphism $s:W\longrightarrow C$ such that $g'=ps$. Then  $(f, s): A\oplus W\longrightarrow C$ is the desired lifting.  This prove that $({\rm CoFib}_{\omega}^\xi, \ {\rm Fib}_{\omega}^\xi, \ {\rm Weq}_{\omega}^\xi)$
is a model structure on $\T$.

\vskip5pt

By Lemma \ref{tcofibandtfib}, the class $\TCoFib_\omega^\xi$ of trivial cofibrations is precisely
the class of splitting monomorphisms with hocokernel in $\omega$, and the class $\TFib_\omega^\xi$
of trivial fibrations is precisely the class of $\xi$-proper epimorphisms with hokernel in $\mathcal Y$. This completes the proof.

\section{\bf Hereditary complete cotorsion pair arising from a model structure}

The aim of this section is to prove the following result, which in particular gives the inverse of Theorem \ref{if}.

\begin{thm}\label{onlyif} \
Let $(\T, [1], \tri)$ be a triangulated category, $\xi$ a proper class of triangles,
$\mathcal{X}$ and $\mathcal{Y}$ additive full subcategories of $\T$ which are closed under isomorphisms and direct summands, and $\omega=\mathcal{X}\cap \mathcal{Y}$.
If  \ $(\CoFib_{\omega}^\xi, \ \Fib_{\omega}^\xi, \ \Weq_{\omega}^\xi)$ as defined in {\rm(\ref{construction})} is a  model structure on $\T$,
then  $(\mathcal{X},\mathcal{Y})$ is a hereditary complete cotorsion pair with respect to $\xi$, and $\omega$ is contravariantly finite in $\mathcal T$.
Moreover, the class of cofibrant objects is $\mathcal X,$ the class of fibrant objects is $\mathcal T,$
the class of trivial objects is $\mathcal Y;$ and the homotopy category is equivalent to additive quotient $\mathcal X/\omega$.
\end{thm}

\subsection{The homotopy category} \

\begin{lem} \label{homtopycat}   \ Let $(\T, [1], \tri)$ be a triangulated category, $\xi$ a proper class of triangles,
 $(\mathcal{X},\mathcal{Y})$ a complete cotorsion pair with respect to $\xi$, and $\omega = \X\cap \Y$.
    If $(\CoFib_{\omega}^\xi, \ \Fib_{\omega}^\xi, \ \Weq_{\omega}^\xi)$ is a model structure,
as given in {\rm(\ref{construction})}, then the homotopy category $\Ho(\mathcal{T})$ is equivalent to additive quotient $\mathcal X/\omega$.
    \end{lem}

\begin{proof} \ We need to apply Theorem \ref{htpthm}. For this we need to verify the conditions in Theorem \ref{htpthm}. Since a triangulated category is additive,
the conditions (i) and  (ii) automatically hold.  It remains to verify the conditions (iv) and (v).

\vskip5pt

We first verify the condition (iv) in Theorem \ref{htpthm}.
For this, let $u: X \xra{} Y$ be an arbitrary trivial cofibration and $\beta: X \xra{} X'$ an arbitrary morphism.
By definition $u$ is a splitting monomorphism with $Z = \hoc u \in \omega$. By Lemma \ref{addaxiom}(6) there exists a commutative diagram
$$\xymatrix@=0.6cm{
  & Z'[-1] \ar@{=}[r]\ar[d]^-{-\gamma[-1]} & Z'[-1] \ar[d] & \\
Z[-1] \ar[r]^-{-w[-1]}\ar@{=}[d] & X \ar[d]^-{\beta}\ar[r]^-{u} & Y \ar[r]^-{v}\ar[d]^-{\beta'} & Z\ar@{=}[d] \\
Z[-1] \ar[r]^-{w'} & X' \ar[r]^-{u'}\ar[d] & Y' \ar[r]^-{v'}\ar[d]^-{\alpha'} & Z\ar[d]^-{-w} \\
  & Z' \ar@{=}[r] &  Z'\ar[r]^-{\gamma} &  X[1]\\
}$$
such that the third row, and the second and the third columns are triangles, and the middle square is homotopy cartesian, i.e.,
$$X \xra{\left(\begin{smallmatrix} \beta \\ -u \end{smallmatrix}\right)} X' \oplus Y \xra{(u', \beta ')} Y' \xra{-wv'} X[1]$$
is a triangle. Since the middle square is homotopy cartesian, it follows that this square is a weak push-out of $u$ and $\beta$ (see Fact \ref{weakpp}).
Since $u$ is a splitting monomorphism, it follows that $w = 0$, and hence $w' = 0$. Thus $u'$ is a splitting monomorphism, and hence $u'$ is a trivial cofibration.
This shows that the condition (iv) in Theorem \ref{htpthm} is satisfied.

\vskip5pt

Now, we verify the condition (v) in Theorem \ref{htpthm}.
For this, let $v: Y \xra{} Z$ be an arbitrary trivial fibration and $\alpha: Z' \xra{} Z$ an arbitrary morphism.
By definition $v$ is a $\xi$-proper epimorphism with $X = \hok v\in\mathcal Y$. By Lemma \ref{addaxiom}(5) there exists a commutative diagram
$$\xymatrix@=0.6cm{
      & X' \ar@{=}[r]\ar[d]_-{\beta'} & X' \ar[d] & \\
    X \ar[r]^-{u'}\ar@{=}[d] & Y' \ar[d]_-{\alpha'}\ar[r]^-{v'} & Z' \ar[r]\ar[d]_-{\alpha} & X[1]\ar@{=}[d] \\
    X \ar[r]^-{u} & Y \ar[r]^-{v}\ar[d] & Z \ar[r]^-{w}\ar[d]_-{\gamma} & X[1]\ar[d]^-{u'[1]} \\
      & X'[1]\ar@{=}[r] & X'[1]\ar[r]^-{-\beta'[1]} &  Y'[1]\\
    }$$
  such that the second row, and the second and the third columns are triangles, and the middle square is homotopy cartesian, i.e.,
  $$Y' \xra{\left(\begin{smallmatrix} \alpha \\ -v' \end{smallmatrix}\right)} Y \oplus Z' \xra{(v, \alpha)} Z \xra{u'[1]w} Y'[1]$$
  is a triangle. Since the middle square is homotopy cartesian, it follows that this square is a weak pull-back of $v$ and $\alpha$ (see Fact \ref{weakpp}). By definition the triangle in the third row is in $\xi$. Since $\xi$ is closed under the base changes, it follows that
the triangle in the second row is also in $\xi$, thus $v'$ is a $\xi$-proper epimorphism with $X = \hok v\in\mathcal Y$,  and hence
$v'$ is a trivial fibration.  This shows that the condition (v) in Theorem \ref{htpthm} is satisfied.

\vskip5pt

Thus, one can apply Theorem \ref{htpthm}, and it suffices to prove $\pi \T_{cf} = \mathcal{X}/ \omega$. For the model structure $(\CoFib_{\omega}^\xi$, \ $\Fib_{\omega}^\xi,$ \ $\Weq_{\omega}^\xi)$,
it is clear that $\mathcal{T}_{cf} = \mathcal{X}$. Let $f, g: A \longrightarrow B$ be morphisms with $A, B \in \mathcal{X}$. It suffices to prove the claim:
$f \overset{l}{\sim} g \iff f-g$ factors through an object in $\omega$.

\vskip5pt

If $f \overset{l}{\sim} g$, i.e., one has a commutative diagram
$$\xymatrix@R=0.4cm{
    A \oplus A\ar[rr]^-{(f,g)}\ar[dd]_{(1,1)}\ar[rrdd]^-{(\partial_1,\partial_2)} && B \\ \\
    A && \widetilde{A}\ar[ll]_-{\sigma}\ar[uu]_h
}$$
with $\sigma \in \Weq_{\omega}^\xi$. We claim that $\sigma$ can be chosen in $\TFib_\omega^\xi$. By definition, there is a $\xi$-proper epimorphism $(\sigma, t): \widetilde{A}\oplus W \longrightarrow A$
with $W\in \omega$ and $\hok (\sigma, t)\in \Y$. Then there is a commutative diagram
 $$\xymatrix@R=0.4cm{
    A \oplus A\ar[rr]^-{(f,g)}\ar[dd]_{(1,1)}\ar[rrdd]^-{\left(\begin{smallmatrix}\partial_1 &\partial_2\\ 0 & 0\end{smallmatrix}\right)} && B \\ \\
    A && \widetilde{A}\oplus W\ar[ll]_-{(\sigma,t)}\ar[uu]_{(h,0)}
}$$
Since $\Ext_\xi(\X, \Y) = 0$, applying Lemma \ref{tcofibandtfib} one sees that $(\sigma, t)\in \TFib_{\omega}^\xi$.
Thus, without loss of generality,  we may assume that $\sigma\in \TFib_{\omega}^\xi$. Note that $f-g = h(\partial_1 - \partial_2)$ and $\sigma (\partial_1 -\partial_2) = 0$. It suffices to show that $\partial_1 -\partial_2$ factors through $\omega$. Since $(\mathcal{X}, \mathcal{Y})$ is a complete cotorsion pair with respect to $\xi$, one can take a triangle $A \overset{i}{\longrightarrow} I \longrightarrow X \longrightarrow A[1]$ in $\xi$ with $I\in\mathcal Y$ and $X\in\mathcal X$.
Then $i \in \CoFib_{\omega}^\xi$.  Since $A\in \mathcal{X}$,  $I \in \mathcal{X} \cap \mathcal{Y} = \omega$. Consider the commutative diagram
$$\xymatrix@R=0.5cm{A\ar[r]^{\partial_1 -\partial_2}\ar[d]_i & \widetilde{A}\ar[d]^{\sigma} \\
    I\ar[r]_0\ar@{-->}[ru] & A
}$$
By Lifting axiom one sees that $\partial_1 -\partial_2$ factors through  object $I$ in $\omega$.

\vskip5pt

Conversely,  if $f-g$ factors through $W \in \omega$, say $f-g = v\circ u$ with $A \overset{u}{\longrightarrow} W \overset{v}{\longrightarrow} B$,
then one has a  commutative diagram
$$\xymatrix@R=0.4cm{
    A \oplus A\ar[rr]^-{(f,g)}\ar[dd]_{(1,1)}\ar[rrdd]^{\left(\begin{smallmatrix} 1 & 1 \\ u & 0 \end{smallmatrix}\right)} && B \\ \\
    A && A \oplus W\ar[ll]^{\sigma = (1,0)}\ar[uu]_{(g,v)}
}$$
where $\sigma = (1, 0) \in \TFib_{\omega}^\xi$ by Lemma \ref{tcofibandtfib}. Thus $\sigma \in \Weq_{\omega}^\xi$. By definition $f \overset{l}{\sim} g$. This proves the claim.
\end{proof}

\subsection{Heredity}
\begin{lem} \label{resolvingcoresolving}   \ Let $(\T, [1], \tri)$ be a triangulated category, $\xi$ a proper class of triangles,
    $(\mathcal{X},\mathcal{Y})$ a complete cotorsion pair with respect to $\xi$, and $\omega = \X\cap \Y$.
    If $(\CoFib_{\omega}^\xi, \ \Fib_{\omega}^\xi, \ \Weq_{\omega}^\xi)$ is a model structure,
as given in {\rm(\ref{construction})}, then the cotorsion pair $(\mathcal{X}, \mathcal{Y})$ is hereditary.
    \end{lem}
    \begin{proof} \ By Theorem \ref{heredity}, it suffices to prove that
$\Y$ is closed under the hocokernels of $\xi$-proper monomorphisms. Suppose that there is a triangle
        \[
        \xymatrix{
        Y_1\ar[r]& Y_2\ar[r]^-{d} & C\ar[r]& Y_1[1]
        }
        \]
    with  $Y_i\in \mathcal{Y}$ for $i = 1, 2.$  By the definition given in {\rm(\ref{construction})}, the morphism $0: Y_2\longrightarrow 0$ is in $\Weq_{\omega}^\xi$, since
    $(0, 0): Y_2\oplus 0\longrightarrow 0$ is a $\xi$-proper epimorphism with $0\in \omega$ and $\hok (0, 0) = Y_2\in\Y$.
    In the similar way, $d: Y_2\longrightarrow C$ is in $\Weq_{\omega}^\xi$, since
    $(d, 0): Y_2\oplus 0\longrightarrow C$ is a $\xi$-proper epimorphism with $\hok (d, 0) = Y_1\in\Y$.

    \vskip5pt

    Since $(Y_2\longrightarrow 0) = (C\longrightarrow 0) \circ d$,  by Two out of three axiom
     the morphism $0: C\longrightarrow 0$ is in $\Weq_{\omega}^\xi$. By definition there is a $\xi$-proper epimorphism $0:  C\oplus W\longrightarrow 0$
    with  $W\in \omega$ and $C\oplus W\in \mathcal{Y}$. Thus $C\in \mathcal  Y$.
    \end{proof}

\subsection{Proof Theorem \ref{onlyif}.} \ By the definition in (\ref{construction})
one easily sees that the class of cofibrant objects and the class of fibrant objects of model structure \ $(\CoFib_{\omega}^\xi, \ \Fib_{\omega}^\xi, \ \Weq_{\omega}^\xi)$ on $\T$ are respectively $\mathcal{X}$ and $\mathcal{T}$.
Also,  the class of trivial objects is $\mathcal{Y}$.
In fact, for any $Y \in \mathcal{Y}$, since $(0, 0): Y\oplus 0\longrightarrow 0$ is a $\xi$-proper epimorphism with $0\in \omega$ and $\hok (0, 0) = Y\in \Y$,
by definition $0: Y\longrightarrow 0$ is a weak equivalence, i.e., $Y$ is a trivial object.
Conversely, if $0: W\longrightarrow 0$ is a weak equivalence,
then there is a $\xi$-proper epimorphism $(0, 0): W \oplus W' \longrightarrow 0$
with $W' \in \omega$ and $\hok(0, 0) = W \oplus W' \in \mathcal{Y}$, thus $W \in \mathcal{Y}$.

\vskip5pt

For any triangle $Y \overset{i}{\longrightarrow} L \overset{d}{\longrightarrow} X\longrightarrow Y[1]$ in $\xi$
with $Y \in \mathcal{Y}$ and $X \in \mathcal{X}$, by definition $i\in\CoFib_\omega^\xi$. Since $Y \in \mathcal{Y}$ is a trivially fibrant object,
$Y \longrightarrow 0$ is a trivial fibration. Thus commutative diagram
$$\xymatrix@R=0.5cm{Y\ar[r]^{{\rm Id}_Y}\ar[d]_i & Y \ar[d]\\
L\ar[r]& 0}$$
admits a lifting, which means that $i$ is a splitting monomorphism. This proves
$\Ext_{\xi}(\X, \Y) = 0$.

\vskip5pt

For any object $T\in\mathcal T$, by Factorization axiom there is a commutative diagram
$$\xymatrix@R=0.5cm{0\ar[rr]\ar[rd] & & T \\ & X\ar[ru]_-f}$$
such that $0\longrightarrow X$ is in $\CoFib_\omega^\xi$ and $f\in \TFib_\omega^\xi$.  By definition $X\in\mathcal X$ is a cofibrant object;
and by Lemma \ref{tcofibandtfib}, $f$ is a $\xi$-proper epimorphism with $Y = \hok f\in \mathcal Y$.
Thus there is a triangle $Y \xra{} X \xra{f} T \xra{} Y[1]$ in $\xi$ with $X\in \X$ and $Y\in \Y$.
Similarly, there is a triangle $T \xra{} Y' \xra{} X' \xra{} T[1]$ in $\xi$ with $X' \in \X$ and $Y' \in \Y$.

\vskip5pt

Now,  $(\X, \Y)$ is a complete cotorsion pair with respect to $\xi$, by Fact \ref{completectp}; and it is also a hereditary cotorsion pair, by Lemma \ref{resolvingcoresolving}.
Also, by Lemma \ref{homtopycat}, the homotopy category is equivalent to additive quotient $\mathcal X/\omega$.

\vskip5pt

\vskip5pt

It remains to see that $\omega$ is contravariantly finite in $\mathcal{T}$.  For any object $T\in\mathcal T$, by Factorization axiom there is a commutative diagram
$$\xymatrix@R=0.5cm{0\ar[rr]\ar[rd] & & T \\ & W\ar[ru]_-g}$$
such that $0\longrightarrow W$ is in $\TCoFib_\omega^\xi$ and $g\in \Fib_\omega^\xi$.
Thus $W$ is a trivially cofibrant object, and hence $W\in\omega$.
For any morphism $h: W'\longrightarrow T$ with $W'\in \omega$, then $0\longrightarrow W'$ is in $\TCoFib_\omega^\xi$, thus commutative diagram
$$\xymatrix@R=0.5cm{0\ar[r]\ar[d] & W \ar[d]^-g\\
W'\ar[r]^-h& T}$$
admits a lifting, i.e., $h$ factors through $g$. Thus $g: W\longrightarrow T$ is a right $\omega$-approximation of $T$.
This proves that $\omega$ is contravariantly finite in $\mathcal{T}$, and hence completes the proof. \hfill $\square$

\section{\bf The Beligiannis - Reiten correspondence}

This section aims to give the Beligiannis - Reiten correspondence, based on Theorems \ref{if} and \ref{onlyif},
a one-one correspondence between
hereditary complete cotorsion pairs  in triangulated category $\T$ with respect to proper class such that
the heart is contravariantly finite, and weakly projective model structures on $\T$.

\subsection{Some necessary facts} Let $(\CoFib ,\  \Fib ,\  \Weq )$ be a model structure on pointed category $\mathcal{A}$ (i.e., a category with zero object).  Put
\begin{align*} & \mathcal C: = \{\mbox{cofibrant objects}\}, \ \ \ \mathcal F: = \{\mbox{fibrant objects}\},
\ \ \ \mathcal W: = \{\mbox{trivial objects}\} \\
& {\rm T}\mathcal C: = \{\mbox{trivially cofibrant objects}\} , \ \ \ {\rm T}\mathcal F: = \{\mbox{trivially fibrant objects}\}.\end{align*}

Using Lifting axiom one easily see the following fact ([BR, VIII, 1.1]).

\begin{fact} \label{factorthrough} \ Let $(\CoFib ,\  \Fib ,\  \Weq )$ be a model structure on pointed category $\mathcal{A}$. Then

\vskip5pt

$(1)$ \ If $p: B \longrightarrow C$ is a trivial fibration $($respectively, a fibration$)$,
then any morphism $\gamma: X \longrightarrow C$ factors through $p$, where $X \in \mathcal{C}$ $($respectively, $X\in {\rm T}\mathcal{C})$.

\vskip5pt

$(2)$ \ If $i: A \longrightarrow B$ is a trivial cofibration $($respectively, a cofibration$)$, then any morphism $\alpha: A \longrightarrow Y$ factors through $i$,
where $Y \in \mathcal{F}$ $($respectively, $Y\in {\rm T}\mathcal{F})$.
\end{fact}

By Fact \ref{factorthrough} one can see the following lemma.

\begin{lem}\label{contravariantlyfinite} \ {\rm ([BR, VIII, 2.1])} \ Let $(\CoFib ,\  \Fib ,\  \Weq )$ be a model structure on pointed category $\mathcal{A}$. Then

\vskip5pt

$(1)$  \ The full subcategory $\mathcal{C}$ $($respectively, ${\rm T}\mathcal{C})$ is contravariantly finite in $\mathcal{A}$. Explicitly,
for any object    $A$ of $\mathcal{A}$, there is a right $\mathcal{C}$-approximation
$($respectively, a right ${\rm T}\mathcal{C}$-approximation$)$ $f: C \longrightarrow A$ with $f \in \TFib$ $($respectively, $f\in \Fib)$.

\vskip5pt

$(2)$  \ The full subcategory $\mathcal{F}$ $($respectively, ${\rm T}\mathcal{F})$ is covariantly finite in $\mathcal{A}$. Explicitly,
for any object    $A$ of $\mathcal{A}$, there is
a left $\mathcal{F}$-approximation $($respectively, a left ${\rm T}\mathcal{F}$-approximation$)$ $g: A \longrightarrow F$ with $g \in \TCoFib$
$($respectively, $g\in \CoFib)$.
\end{lem}

For abelian categories, the following result is in [BR, VIII, Lemma 3.2], with a slight difference.

\begin{lem}\label{extinclusion} \ {\rm ([BR, VIII, 3.2])} \ Let  $(\CoFib ,\  \Fib ,\  \Weq )$ be a model structure on a triangulated category $(\T, [1], \tri)$,
and $\xi$ a proper class of triangles.

\vskip5pt

$(1)$ \ \ If \ $\{\xi\mbox{-proper monomorphism} \ f \ | \ \hoc f \in \mathcal C\} \subseteq \CoFib$, then $\Ext_{\xi} (\mathcal{C}, {\rm T}\mathcal{F}) = 0.$

\vskip5pt

$(2)$  \ \ If \ $\TFib\subseteq \{\xi\mbox{-proper epimorphism} \ f \ | \ \hok f \in {\rm T}\mathcal F\}$,
then \ $^{\perp}{\rm T}\mathcal{F} \subseteq \mathcal{C}$.

\vskip5pt

$(3)$ \ \ If \ $\CoFib \subseteq \{\xi\mbox{-proper monomorphism} \ f \ | \ \hoc f \in \mathcal C\},$
then $\mathcal{C}^{\perp} \subseteq {\rm T}\mathcal{F}.$
\end{lem}

\begin{proof} \ We only justify (1) and (2).

\vskip5pt

(1) \ For any triangle $Y \overset{i}{\longrightarrow} L \overset{d}{\longrightarrow} C\longrightarrow Y[1]$ in $\xi$
with $Y \in {\rm T}\mathcal{F}$ and $C \in \mathcal{C}$, by the assumption $i$ is a cofibration.
Thus by Fact \ref{factorthrough}(2), $\Id_Y : Y \longrightarrow Y$ factors through $i$, i.e., $i$ is a splitting monomorphism.

\vskip5pt

(2) \ For any object $A \in \ ^{\perp}{\rm T}\mathcal{F}$, by Lemma \ref{contravariantlyfinite}(1)
there is a right $\mathcal{C}$-approximation  $f: C \longrightarrow A$ with $f \in \TFib$.
By the assumption one gets a triangle $Y \longrightarrow C \overset{f}{\longrightarrow} A\longrightarrow Y[1]$ in $\xi$
with $Y \in {\rm T}\mathcal{F}$. Since $A \in \ ^{\perp}{\rm T}\mathcal{F}$ and $Y \in {\rm T}\mathcal{F}$, this triangle splits.
Hence $A$ is a direct summand of $C$, and thus $A \in \mathcal{C}$.
\end{proof}

\subsection{\bf Weakly $\xi$-projective model structures on triangulated category with proper class $\xi$ of triangles}
We keep the notations in Subsection 6.1. So $\mathcal C$ (respectively, $\mathcal F$, ${\rm T}\mathcal C$, and ${\rm T}\mathcal F$)
is the class of cofibrant objects (respectively, fibrant objects, trivially cofibrant objects, and trivially fibrant objects).

\begin{prop}\label{weaklyproj} \ Let $(\T, [1], \tri)$ be a triangulated category, $\xi$ a proper class of triangles,
and  $(\CoFib ,\  \Fib ,\  \Weq )$  a model structure on  $\T$. Then the following are equivalent.

\vskip5pt

$(1)$ \ $\CoFib = \{\mbox{$\xi$-proper monomorphism} \ f \ | \  \hoc f\in \mathcal C\}$,  and \
$${\rm T}\Fib \subseteq \{\mbox{$\xi$-proper epimorphism} \ f \ | \ \hok f\in {\rm T}\mathcal F\}.$$

\vskip5pt

$(2)$   \ $\Ext_{\xi} (\mathcal{C}, {\rm T}\mathcal{F}) = 0,$ \ $\CoFib \subseteq \{\mbox{$\xi$-proper monomorphism} \ f \ | \  \hoc f\in \mathcal C\}$,
and \
$${\rm T}\Fib \subseteq \{\mbox{$\xi$-proper epimorphism} \ f \ | \ \hok f\in {\rm T}\mathcal F\}.$$

\vskip5pt

$(3)$ \ $\CoFib \subseteq \{\mbox{$\xi$-proper monomorphism} \ f \ | \  \hoc f\in \mathcal C\}$, and  $${\rm T}\Fib = \{\mbox{$\xi$-proper epimorphism} \ f \ |  \  \hok f\in {\rm T}\mathcal F\}.$$

\vskip5pt

$(4)$ \ $\CoFib = \{\mbox{$\xi$-proper monomorphism} \ f \ | \  \hoc f\in \mathcal C\}$, and
$${\rm T}\Fib = \{\mbox{$\xi$-proper epimorphism} \ f \ | \  \hok f\in {\rm T}\mathcal F\}.$$

\vskip5pt

$(5)$ \ $(\mathcal{C}, {\rm T}\mathcal{F})$ is a complete cotorsion pair with respect to $\xi$, \ $\CoFib \subseteq \{\mbox{morphism} \ f \ | \  \hoc f\in \mathcal C\}$,
and \
${\rm T}\Fib \subseteq \{\mbox{morphism} \ f \ | \ \hok f\in {\rm T}\mathcal F\}.$

\end{prop}
\begin{proof} \ $(1) \Longrightarrow (2)$: This  follows from Lemma
\ref{extinclusion}(1).

\vskip5pt

$(2) \Longrightarrow (1)$:   Let $i: A\longrightarrow B$ be a $\xi$-proper monomorphism with $\hoc f\in \mathcal C$. For any $p\in \TFib$,
by assumption $p$ is a $\xi$-proper epimorphism with $\hok p\in \rm T\mathcal F$.
Since $\Ext_\xi(\mathcal C, \rm T\mathcal F) = 0$, one can apply
Extension-Lifting Lemma \ref{extlifting} to see that $i$ has
the left lifting property with respect to $p$. Thus $i$ is a cofibration, by Lemma \ref{quillenlifting}.
Thus
$$\{\xi\mbox{-proper monomorphism} \ f \ | \ \hoc f \in \mathcal C\} \subseteq \CoFib.$$
Then by assumption one has equality \ $\CoFib = \{\xi\mbox{-proper monomorphism} \ f \ | \ \hoc f \in \mathcal C\}.$

\vskip5pt

The equivalence $(2)\Longleftrightarrow (3)$ can be proved similarly.

\vskip5pt

$(4) \Longrightarrow (1)$ is clear; and $(1) \Longrightarrow (4)$ is also clear, since $(1)$ and $(3)$ imply $(4)$.

\vskip5pt

$(1) \Longrightarrow (5):$  By the assumptions and Lemma \ref{extinclusion}(1), (2) and (3),
$(\mathcal{C}, {\rm T}\mathcal{F})$ is a cotorsion pair.

\vskip5pt

By Lemma \ref{contravariantlyfinite}(1), for any object $A \in \T$, there exists a right $\mathcal{C}$-approximation  $f: C \longrightarrow A$ such that $f \in \TFib$. Then by assumption there is a triangle $Y \longrightarrow C \stackrel f\longrightarrow A \longrightarrow Y[1]$ in $\xi$ with $Y \in {\rm T}\mathcal{F}$.

\vskip5pt

Similarly, by Lemma \ref{contravariantlyfinite}(2) one has a triangle $A \longrightarrow Y' \longrightarrow C' \longrightarrow A[1]$ in $\xi$ with $Y' \in {\rm T}\mathcal{F}$ and $C' \in \mathcal{C}$. Thus, the cotorsion pair $(\mathcal{C}, {\rm T}\mathcal{F})$ is complete.

\vskip5pt

$(5) \Longrightarrow (2):$  Let $f: A\longrightarrow B$ be a cofibration. By assumption $\hoc f\in \mathcal C$. By the completeness of the cotorsion pair $(\mathcal{C}, {\rm T}\mathcal{F})$,
there is a triangle $A\stackrel i \longrightarrow Y \longrightarrow X \longrightarrow A[1]$ in $\xi$ such that $Y\in {\rm T}\mathcal{F}$ and $X\in\mathcal C$.
By Fact \ref{factorthrough}(2), $i$ factors through $f$,
say $i = gf$ for some $g: B\longrightarrow Y$. By definition $i$ is a $\xi$-proper monomorphism. It follows from Lemma \ref{monicepimorphism}
that $f$ is a $\xi$-proper monomorphism. This proves $\CoFib \subseteq \{\mbox{$\xi$-proper monomorphism} \ f \ | \  \hoc f\in \mathcal C\}$.

\vskip5pt

Similarly, one has \ ${\rm T}\Fib \subseteq \{\mbox{$\xi$-proper epimorphism} \ f \ | \ \hok f\in {\rm T}\mathcal F\}.$ This completes the proof.
\end{proof}

\begin{defn}\label{defwp} \ Let $(\T, [1], \tri)$ be a triangulated category, and $\xi$ a proper class of triangles.
A model structure on $\T$  is  {\it weakly $\xi$-projective,} provided
that each object is fibrant and it satisfies the equivalent conditions in {\rm Proposition \ref{weaklyproj}.}
\end{defn}

\subsection{Beligiannis - Reiten correspondence}

By the following result,
the $(\xi, \omega)$-model structures on triangulated category $\T$ with proper class $\xi$
are precisely the weakly $\xi$-projective model structures.

\begin{thm} \label{brcorrespondence} \ Let $(\T, [1], \tri)$ be a triangulated category, and $\xi$ a proper class of triangles.
Let $\Omega$ denote the class of hereditary complete cotorsion pairs $(\X, \Y)$ in $\T$ with respect to $\xi$, such that
$\omega=\X\cap \Y$ is contravariantly finite in $\T$, and $\Gamma$ the class of weakly $\xi$-projective model structures on $\T$.
Then the map $$\Phi: (\X, \Y)\mapsto (\CoFib_{\omega}^\xi, \ \Fib_{\omega}^\xi, \ \Weq_{\omega}^\xi)$$
as given in ${\rm(\ref{construction})},$ is a bijection between $\Omega$ and $\Gamma$, with the inverse given by
$\Psi: (\CoFib, \Fib, \Weq)\mapsto (\mathcal{C}, \rm T\mathcal{F})$.
\end{thm}

\begin{proof}  \ By Theorem \ref{if}, $\rm{Im}\Phi\in \Gamma$ and $\Psi\Phi={\rm Id}_\Omega$.
It remains to prove $\rm{Im}\Psi\in \Omega$ and $\Phi \Psi ={\rm Id}_\Gamma$.

\vskip5pt

For this purpose, let $(\CoFib, \Fib, \Weq)$ be a  weakly $\xi$-projective model structure on $\T$.
By Proposition \ref{weaklyproj}(5), $(\mathcal{C}, {\rm T}\mathcal{F})$ is a complete cotorsion pair.
By definition $\mathcal{F} = \mathcal{T}$. Thus
$$\mathcal{C}\cap \rm T\mathcal{F} = \mathcal{C}\cap \mathcal{F} \cap \mathcal{W} = \mathcal{C}\cap \mathcal{W} = {\rm T}\mathcal{C}.$$
Then,  by Lemma \ref{contravariantlyfinite}(1),  $\mathcal{C}\cap \rm T\mathcal{F} = {\rm T}\mathcal{C}$ is contravariantly finite in $\mathcal{T}$.

\vskip5pt

It remains to prove that cotorsion pair $(\mathcal{C}, {\rm T}\mathcal{F})$ is hereditary (and hence $\rm{Im}\Psi\in \Omega$), and that
$(\CoFib, \Fib, \Weq) = (\CoFib_{\omega}^\xi, \ \Fib_{\omega}^\xi, \ \Weq_{\omega}^\xi)$, where $\omega = \mathcal{C}\cap \rm T\mathcal{F}= {\rm T}\mathcal{C}$.
This will be done in several steps.

\vskip5pt

Since $(\CoFib, \Fib, \Weq)$ is a weakly $\xi$-projective model structure, by Proposition \ref{weaklyproj}(4) one has already
$$\CoFib = \{\mbox{$\xi$-proper monomorphism} \ f \ | \  \hoc f\in \mathcal C\} = \CoFib_{\omega}^\xi$$ and
$${\rm T}\Fib = \{\mbox{$\xi$-proper epimorphism} \ f \ | \  \hok f\in {\rm T}\mathcal F\} = \TFib_{\omega}^\xi.$$

\vskip5pt

{\bf Claim 1:} \ $\TCoFib= \{\mbox{splitting monomorphism} \ f \ | \ \hoc f\in {\rm T}\mathcal C\} = \TCoFib_\omega^\xi$.

\vskip5pt

In fact, let $f: A\longrightarrow B$ be a splitting monomorphism with $\hoc f\in \rm T\mathcal{C}$. Then
there are morphisms $i: B\longrightarrow A$ and $p: \hoc f \longrightarrow B$ such that
$i\circ f= {\rm Id}_A, \ \ \pi \circ p = {\rm Id}_{\hoc f}, \ \ i \circ p = 0$, where $\pi: B \longrightarrow \hoc f.$ Then it is clear that the square
\[
\xymatrix@R=0.8cm{
0\ar[d]\ar[r] &A\ar[d]^{f}\\
\hoc f\ar[r]^-p &B = A\oplus \hoc f
}
\]
is a push-out square. Since $\hoc f$ is a trivially cofibrant object, $0\longrightarrow \hoc f$ is in $\TCoFib$. It follows from Fact \ref{elementpropmodel}(3) that $f\in \TCoFib$.

\vskip5pt

Conversely, let $f: A\longrightarrow B$ be a trivial cofibration. Then  $f\in \CoFib$, and hence $f$ is a $\xi$-proper monomorphism and $\hoc f\in \mathcal{C}$. Since $A\in \mathcal{T}=\mathcal F$, it follows from  Fact \ref{factorthrough}(2) that ${\rm Id}_A: A\longrightarrow A$ facts through $f$, i.e., $f$ is a splitting monomorphism. Hence $\hoc f = \Coker f$ (cf. Definition - Fact \ref{split}). Then one has the push-out square
\[
\xymatrix@R=0.8cm{
A\ar[d]_-{f}\ar[r] &0\ar[d] \\
B\ar[r]^-{\pi} & \hoc f = \Coker f
}
\]
Since $f$ is a trivial cofibration, it follows from Fact \ref{elementpropmodel}(3) that $0\longrightarrow \hoc f$ is also a trivial cofibration, i.e.,  $\hoc f \in {\rm T}\mathcal C$.
This completes the proof of {\bf Claim 1}.

\vskip5pt

{\bf Claim 2:} \ $\Weq= \Weq_\omega^\xi$. This follow from
\ $\Weq=\TFib \circ \TCoFib
= \TFib^\xi_\omega\circ \TCoFib^\xi_\omega= \Weq_\omega^\xi$.

\vskip5pt

{\bf Claim 3:} \ $\Fib = \{\text{morphism} \ f \ | \ \Hom_{\mathcal T}(W, f) \ \mbox{is surjective}, \ \forall \ W\in \omega\} = \Fib_\omega^\xi$.

\vskip5pt

In fact, by {\bf Claim 1} and using the fact that $\Fib$ is precisely the class of morphisms which have the right lifting property with respect to all the trivial cofibrations
(cf. Lemma \ref{quillenlifting}(3)) one can easily see this: because that trivial cofibrations are splitting monomorphisms with hocokernel in $\rm T\mathcal C$,
and that a morphism $p$ has
the right lifting property with respect to trivial cofibrations is amount to say that $\Hom_{\mathcal T}(W, f)$ is surjective for all $W\in \omega$.

\vskip5pt

Noe we have proved $\CoFib = \CoFib_\omega^\xi, \ \ \Fib = \Fib_\omega^\xi, \ \ \Weq = \Weq_\omega^\xi.$ Thus
$(\CoFib_{\omega}^\xi, \Fib_{\omega}^\xi, \Weq_{\omega}^\xi)$ is also a model structure.
It follows from Lemma \ref{resolvingcoresolving} that cotorsion pair $(\mathcal{C}, {\rm T}\mathcal{F})$ is hereditary.
Thus $\rm{Im}\Psi\in \Omega$ and $\Phi \Psi ={\rm Id}_\Gamma$.  This completes the proof.\end{proof}

\section{\bf The dual version} \ We state the dual version of the main results without proofs.
Let $(\T, [1], \tri)$ be a triangulated category,
and $\xi$ a proper class of triangles.  Suppose that $\mathcal{X}$ and $\mathcal{Y}$ are additive full subcategories of $\T$ which are closed under isomorphisms and direct summands, and $\omega:=\mathcal{X}\cap \mathcal{Y}$. Define three classes of morphisms in $\T$ as follows.

\begin{equation}\begin{aligned}\label{coconstruction}\CoFib^\omega_\xi & = \{\text{morphism} \ f \ | \ \Hom_{\mathcal T}(f, W) \ \mbox{is surjective}, \ \forall \ W\in \omega\}\\
\Fib^\omega_\xi& = \{\ \xi\text{-proper epimorphism} \ f \ | \ \hok f\in \mathcal Y\}\\
\Weq^\omega_\xi & = \{\ A\stackrel f \longrightarrow B \ | \ \mbox{there is a} \ \xi\mbox{-proper monomorphism} \ \left(\begin{smallmatrix} f\\t \end{smallmatrix}\right): A \longrightarrow B\oplus W
\\ & \ \ \ \ \ \ \ \ \ \ \ \ \ \ \ \ \ \ \ \ \ \mbox{such that} \ W\in \omega \ \ \mbox{and} \ \ \hoc \left(\begin{smallmatrix} f\\t \end{smallmatrix}\right)\in \mathcal{X}\}
\end{aligned}\end{equation}

\vskip5pt

\begin{thm}\label{dualmainthm} \ Let $(\T, [1], \tri)$ be a triangulated category,
and $\xi$ a proper class of triangles.  Suppose that $\mathcal{X}$ and $\mathcal{Y}$ are additive full subcategories of $\T$ which are closed under isomorphisms and direct summands, and $\omega=\mathcal{X}\cap \mathcal{Y}$.
Then $({\rm CoFib}^{\omega}_\xi, \ {\rm Fib}^{\omega}_\xi, \ {\rm Weq}^{\omega}_\xi)$  as defined in {\rm(\ref{coconstruction})} is a model structure on $\T$ if and only if
$(\mathcal{X}, \mathcal{Y})$ is a hereditary complete cotorsion pair in $(\T, [1], \tri)$ with respect to $\xi$,
and $\omega$ is covariantly finite in $\T$.
\vskip5pt

If this is the case, then the class $\TCoFib^\omega_\xi$ of trivial cofibrations is precisely
the class of $\xi$-proper monomorphisms with hocokernel in $\mathcal X$, and the class $\TFib^\omega_\xi$
of trivial fibrations is precisely the class of splitting epimorphisms with hokernel in $\omega$;
the class of cofibrant objects is $\mathcal T,$ the class of fibrant objects is $\mathcal Y,$
the class of trivial objects is $\mathcal X;$ and the homotopy category is equivalent to additive quotient $\mathcal Y/\omega$.\end{thm}

\vskip5pt

A model structure $(\CoFib, \Fib, \Weq)$ on $\mathcal A$ is {\it weakly $\xi$-injective} if
$\TCoFib$ is precisely the class of $\xi$-proper monomorphisms with trivially cofibrant hockernel,
$\Fib$ is precisely the class of $\xi$-proper epimorphisms with fibrant hocokernel, and each object is cofibrant.

\vskip5pt

\begin{thm} \label{dualbrcorrespondence} \ Let $(\T, [1], \tri)$ be a triangulated category, and $\xi$ a proper class of triangles.
Let $\Omega$ denote the class of hereditary complete cotorsion pairs $(\X, \Y)$ in $\T$ with respect to $\xi$, such that
$\omega=\X\cap \Y$ is covariantly finite in $\T$, and $\Gamma$ the class of weakly $\xi$-injective model structures on $\T$.
Then the map $$\Phi: (\X, \Y)\mapsto (\CoFib^{\omega}_\xi, \ \Fib^{\omega}_\xi, \ \Weq^{\omega}_\xi)$$
as defined in ${\rm(\ref{coconstruction})},$ gives a bijection between $\Omega$ and $\Gamma$, with  inverse given by  \
$\Psi: (\CoFib, \Fib, \Weq)$ $\mapsto (\rm T\mathcal{C}, \mathcal{F})$,
where $\rm T\mathcal C$ and $\mathcal{F}$ are respectively the class of trivially cofibrant objects and the class of  fibrant objects.
\end{thm}

\section {\bf A $\xi$-triangulated model structure}

\subsection{A $\xi$-triangulated model structure}  M. Hovey [Hov2] has introduced {\it abelian model structure} on abelian category.
This has been extended as {\it exact model structure} on exact category by J. Gillespie [G1].
For a triangulated category and a proper class $\xi$ of triangles,  X. Y. Yang [Y] called the corresponding version as a $\xi$-triangulated model structure.

\begin{defn} \ \label {defnxitriangulatedmds}{\rm ([Y])} \
    \ Let $\mathcal T$ be a triangulated category, and $\xi$ a proper class of triangles.

\vskip5pt

(1) \ A model structure on $\T$ is {\it a $\xi$-triangulated model structure}, provided that the following conditions are satisfied:

    \vskip5pt

\hskip20pt     (i) \ A morphism is a  cofibration (respectively, trivial cofibration) if and only if it is a $\xi$-proper monomorphism with cofibrant (respectively, trivially cofibrant)  hocokernel;

\hskip20pt    (ii) \ A morphism is a fibration (respectively, trivial fibration) if and only if it is a $\xi$-proper epimorphism with  fibrant (respectively, trivially fibrant) hokernel.

\vskip5pt

(2) \ A triple $(\mathcal{C}, \mathcal{F}, \mathcal{W})$ of classes of objects of $\T$ is {\it a {\rm Hovey} triple in $\T$ with respect to $\xi$}, provided that the following conditions are satisfied:

\hskip20pt     (i) \   $\mathcal{W}$ is a $\xi$-thick subcategory of $\T$, i.e., $\mathcal{W}$ is closed under isomorphisms and direct summands,
and if $X \longrightarrow Y\longrightarrow Z\longrightarrow X[1]$ is a triangle in $\xi$ and if two of the three terms $X$, $Y$ and $Z$ are in $\mathcal{W}$, then so is the third;

\hskip20pt    (ii) \ $(\mathcal{C}, \mathcal{F} \cap \mathcal W)$ and $(\mathcal{C}\cap \mathcal W, \mathcal{F})$ are complete cotorsion pairs with respect to $\xi$.

\vskip5pt

(3) \ A {\rm Hovey} triple $(\mathcal{C}, \mathcal{F}, \mathcal{W})$ in $\T$ with respect to $\xi$ is {\it hereditary}, provided that the complete cotorsion pairs $(\mathcal{C}, \mathcal{F} \cap \mathcal W)$ and $(\mathcal{C}\cap \mathcal W, \mathcal{F})$ with respect to $\xi$ are hereditary.
\end{defn}

\vskip5pt

The following Hovey correspondence is discovered by M. Hovey [Hov1] for abelian model structures.
It is generalized to exact model structures by J. Gillespie [G1]; and further generalized to extriangulated categories
by H. Nakaoka and Y. Palu [NP, Section 5]; and moreover, quite surprising, the corresponding homotopy category is a triangulated category,
even if the two complete cotorsion pairs induced by the Hovey triple are not necessarily hereditary,
see [NP, Theorem 6.20] and [G2, Theorem 6.34].
Extriangulated category is a common generalization of exact category and triangulated category,
however, it can not include the case for a triangulated category with a proper class $\xi$ of triangles.

\begin{thm} \label{abelianhtpcat} {\rm ([Y, Theorem A])} \ Let $\mathcal T$ be a triangulated category, and $\xi$ a proper class of triangles.
Then there is a one-one correspondence between $\xi$-triangulated model structures and the {\rm Hovey} triples in $\mathcal T$ with respect to $\xi$, given by
$$({\rm CoFib}, \ {\rm Fib}, \ {\rm Weq})\mapsto (\mathcal{C}, \ \mathcal{F}, \ \mathcal W)$$
where \ $\mathcal C  = \{\mbox{cofibrant objects}\}, \ \
\mathcal F  = \{\mbox{fibrant objects}\}, \ \
\mathcal W  = \{\mbox {trivial objects}\}$, with the inverse \ \ $(\mathcal{C}, \ \mathcal{F}, \ \mathcal W) \mapsto ({\rm CoFib}, \ {\rm Fib}, \ {\rm Weq}),$ where
\begin{align*} &{\rm CoFib} = \{\mbox{$\xi$-proper monomorphism with hocokernel in} \ \mathcal{C}\}, \\ &
{\rm Fib}  = \{\mbox{$\xi$-proper epimorphism with hokernel in} \ \mathcal{F} \}, \\
& {\rm Weq}  = \{pi \ \mid \ i \ \mbox{is a $\xi$-proper monomorphism,} \ \hoc i\in \mathcal{C}\cap \mathcal W, \\ & \ \ \ \ \ \ \ \ \ \ \ \ \ \ \ \ \  p \ \mbox{is a $\xi$-proper epimorphism,} \ \hok p\in \mathcal{F}\cap \mathcal W\}.\end{align*}
If this is the case, then
\begin{align*} &{\rm TCoFib} = \{\mbox{$\xi$-proper monomorphism with hocokernel in} \ \mathcal C\cap \mathcal W\}
\\ & {\rm TFib} = \{\mbox{$\xi$-proper epimorphism with hokernel in}  \ \mathcal{F}\cap \mathcal W \}\end{align*}
and the homotopy category $\Ho(\mathcal{T})$ is equivalent to additive quotient \ $(\mathcal{C} \cap \mathcal{F})/(\mathcal{C} \cap \mathcal{F} \cap \mathcal W)$.
\end{thm}

\begin{proof} \ Everything has been proved in [Y], except for the part for $\Ho(\mathcal{T})$.

\vskip5pt

Similarly to the proof of Lemma \ref{homtopycat},
one can see that the conditions (iv) and (v) in Theorem \ref{htpthm} also hold for triangulated $\xi$-model structures,
and hence one can apply Theorem \ref{htpthm}. Thus, it remains to prove $\pi \T_{cf} = (\mathcal{C} \cap \mathcal{F})/(\mathcal{C} \cap \mathcal{F}  \cap \mathcal W)$.
Let $f, g: A \longrightarrow B$ be morphisms with $A, B \in \T_{cf}= \mathcal{C} \cap \mathcal{F}$. It suffices to prove the claim:
$f \overset{l}{\sim} g \iff f-g$ factors through an object in $\mathcal{C} \cap \mathcal{F}  \cap \mathcal W$.

\vskip5pt

Assume that $f \overset{l}{\sim} g$. By definition there is a commutative diagram
$$\xymatrix@R=0.4cm{
    A \oplus A\ar[rr]^-{(f,g)}\ar[dd]_{(1,1)}\ar[rrdd]^-{(\partial_1,\partial_2)} && B \\ \\
    A && \widetilde{A}\ar[ll]_-{\sigma}\ar[uu]_h}$$ with $\sigma \in \Weq$.
By Factorization axiom there is a commutative diagram
$$\xymatrix@R=0.4cm{
    \widetilde{A} \ar[rr]^-{\sigma}\ar[rd]_{i} && A \\ & D\ar[ru]_p}$$
i.e., $\sigma = p \circ i$,  where $i\in \CoFib$  and $p\in {\rm T}\Fib\subseteq \Weq$. Since $\sigma\in \Weq$ and $p\in \Weq$, it follows that $i\in\Weq$ and hence $i\in {\rm T}\CoFib$.
Since $B\in \mathcal F$, $B\longrightarrow 0$ is in $\Fib$. By  Lifting axiom one gets a morphism $h': D\longrightarrow B$ such that $h = h'\circ i.$
$$\xymatrix@R=0.5cm{\widetilde{A}\ar[r]^{h}\ar[d]_i & B\ar[d] \\
    D\ar[r]\ar@{-->}[ru]^-{h'} &  0}$$

\vskip5pt

Since $(\mathcal{C}, \mathcal{F}\cap \mathcal W)$ a complete cotorsion pair with respect to $\xi$,
one can take a triangle $A \overset{j}{\longrightarrow} W \longrightarrow C \longrightarrow A[1]$ in $\xi$ with $W\in\mathcal{F}\cap \mathcal W$ and $C\in\mathcal C$.
By definition $j \in \CoFib$.
Since $A\in \mathcal C$ and $C\in\mathcal C$, it follows that $W\in\mathcal{C}$ and hence $W\in\mathcal{C}\cap \mathcal{F}\cap\mathcal W$.
Since $p\circ i \circ (\partial_1 - \partial_2) = \sigma \circ (\partial_1 - \partial_2) = 0$, one has commutative diagram
$$\xymatrix@R=0.5cm{A\ar[rr]^{i\circ (\partial_1 -\partial_2)}\ar[d]_j && D\ar[d]^{p} \\
    W\ar[rr]_0\ar@{-->}[rru]^-t & & A
}$$
Since $j \in \CoFib$ and $p\in {\rm T}\Fib$, it follows from Lifting axiom that $i\circ (\partial_1 -\partial_2)$ factors through  object $W\in \mathcal{C}\cap \mathcal{F}\cap\mathcal W$.

\vskip5pt

Now, one has
$$f-g = h(\partial_1 - \partial_2) = h'\circ i \circ (\partial_1 - \partial_2) = h'\circ t \circ j.$$
This proves that $f - g$ factors through object $W\in \mathcal{C} \cap \mathcal{F}  \cap \mathcal W$.

\vskip5pt

Conversely,  if $f-g$ factors through object $W \in \mathcal{C} \cap \mathcal{F}  \cap \mathcal W$, say $f-g = v\circ u$ with $A \overset{u}{\longrightarrow} W \overset{v}{\longrightarrow} B$,
then one has a  commutative diagram
$$\xymatrix@R=0.4cm{
    A \oplus A\ar[rr]^-{(f,g)}\ar[dd]_{(1,1)}\ar[rrdd]^{\left(\begin{smallmatrix} 1 & 1 \\ u & 0 \end{smallmatrix}\right)} && B \\ \\
    A && A \oplus W\ar[ll]^{\sigma = (1,0)}\ar[uu]_{(g,v)}
}$$
where $\sigma = (1, 0) \in \TFib$. Thus $\sigma \in \Weq$. By definition $f \overset{l}{\sim} g$. This proves the claim.
\end{proof}

\subsection{\bf When is the $(\xi, \omega)$-model structure $\xi$-triangulated?}

It is natural to ask when the $(\xi, \omega)$-model structure $(\CoFib_{\omega}^\xi, \Fib_{\omega}^\xi, \Weq_{\omega}^\xi)$  is $\xi$-triangulated.
To answer this question we need some terminologies.

\begin{defn} \ Let $(\T, [1], \tri)$ be a triangulated category, and $\xi$ a proper class of triangles.

\vskip5pt

(1) \  {\rm ([B, 4.1];  also [AS])} \ An object $P\in \mathcal T$ is {$\xi$-projective},
provided that for any triangle $X \longrightarrow Y\longrightarrow Z \longrightarrow X[1]$ in $\xi$, the induced sequence
$$0 \longrightarrow \Hom_\mathcal T(P, X)\longrightarrow \Hom_\mathcal T(P, Y)\longrightarrow \Hom_\mathcal T(P, Z)\longrightarrow 0 $$
is an exact sequence of abelian groups.

\vskip5pt

Let $\mathcal P(\xi)$ denote the class of $\xi$-projective objects.

\vskip5pt

(2) \  $\mathcal T$ has enough $\xi$-projective objects, if for any object $X\in \mathcal T$ there is
a triangle $K \longrightarrow P\longrightarrow X \longrightarrow K[1]$ in $\xi$ such that $P$ is $\xi$-projective.

\vskip5pt

(3) \ {\rm ([G, 4.5])} \ A model structure $(\CoFib, \Fib, \Weq)$ on $\mathcal T$ is {\it $\xi$-projective} if it is $\xi$-triangulated, and each object is fibrant.

\vskip5pt

(4) \ {\rm (\cite[1.6, 7.11]{CRZ}, \cite[1.1.9]{Bec})} \ A complete cotorsion pair $(\X,\Y)$ with respect to $\xi$
is {\it generalized $\xi$-projective} if $\Y$ is $\xi$-thick and $\X\cap \Y=\mathcal P(\xi)$.
\end{defn}

\vskip5pt

\begin{rem} \label{xiprojxictp} \ Let $(\T, [1], \tri)$ be a triangulated category, and $\xi$ a proper class of triangles.

\vskip5pt

{\rm (1)}  \ For the proper class of triangles $\xi = \xi(\Hom_{\T}(M,-))$ $($see {\rm Lemma \ref{exmproperclass}}$)$, add$M \subseteq \mathcal{P}(\xi)$. In fact, for any triangle $X\xra{} Y \xra{} Z \xra{} X[1]$ in $\xi$, one has an exact sequence $0 \xra{} \Hom_{\T}(M, X) \xra{} \Hom_{\T}(M, Y) \xra{} \Hom_{\T}(M,Z) \xra{} 0$ by definition. Thus $M \in \mathcal{P}(\xi)$.

\vskip5pt

{\rm (2)}  \ $(\mathcal P(\xi), \T)$ is always a cotorsion pair with respect to $\xi$; and moreover
the  cotorsion pair $(\mathcal P(\xi), \T)$ with respect to $\xi$ is complete if and only if $\mathcal T$ has enough $\xi$-projective objects.
If this is the case, then $(\mathcal P(\xi), \T)$ is a hereditary complete cotorsion pair with respect to $\xi$.

\vskip5pt
The  cotorsion pair $(\mathcal P(\xi), \T)$ with respect to $\xi$ is called {\it the $\xi$-projective cotorsion pair with respect to $\xi$}.

\end{rem}

\begin{prop}\label{thesame} \ Let $(\T, [1], \tri)$ be a triangulated category,  $\xi$ a proper class of triangles,
and $(\mathcal{X}, \mathcal{Y})$ a hereditary complete cotorsion pair with respect to $\xi$ such that $\omega=\mathcal{X}\cap \mathcal{Y}$ contravariantly finite in $\T$.
Then the $(\xi, \omega)$-model structure $(\CoFib_{\omega}^\xi, \ \Fib_{\omega}^\xi, \ \Weq_{\omega}^\xi)$ as defined in {\rm(\ref{Introconstruction})}
is $\xi$-triangulated if and only if $\mathcal T$ has enough $\xi$-projective objects and $\omega = \mathcal P(\xi)$.
\end{prop}
\begin{proof} \ Assume that the model structure $(\CoFib_{\omega}^\xi, \ \Fib_{\omega}^\xi, \ \Weq_{\omega}^\xi)$  is $\xi$-triangulated. By Theorem \ref{abelianhtpcat},
$(\mathcal X, \ \T, \ \Y)$ is a Hovey triple with respect to $\xi$, and hence $(\omega, \mathcal{T})$ is a complete cotorsion pair with respect to $\xi$.
It remains to show $\omega = \mathcal{P}(\xi)$ (and then $\mathcal T$ has enough $\xi$-projective objects, by the completeness).

\vskip5pt

If $P\in \omega$, then for any triangle $X \longrightarrow Y\longrightarrow Z \longrightarrow X[1]$ in $\xi$,
$X[-1] \longrightarrow Y[-1]\longrightarrow Z[-1] \longrightarrow X$ is also a triangle in $\xi$.
Thus one has exact sequence
$$\Hom_\mathcal T(P, X[-1])\longrightarrow \Hom_\mathcal T(P, Y[-1])\longrightarrow \Hom_\mathcal T(P, Z[-1])
\longrightarrow \Hom_\mathcal T(P, X)\longrightarrow \Hom_\mathcal T(P, Y).$$
On the other hand, by Theorem \ref{exactseq} one has exact sequences
$$\Hom_\mathcal T(P, X[-1])\longrightarrow \Hom_\mathcal T(P, Y[-1])\longrightarrow \Hom_\mathcal T(P, Z[-1])\longrightarrow \Ext_\xi(P, X[-1]) =0$$
$$\Hom_\mathcal T(P, X)\longrightarrow \Hom_\mathcal T(P, Y)\longrightarrow \Hom_\mathcal T(P, Z)\longrightarrow \Ext_\xi(P, X) = 0.$$
Putting together one gets exact sequence
$$0\longrightarrow \Hom_\mathcal T(P, X)\longrightarrow \Hom_\mathcal T(P, Y)\longrightarrow \Hom_\mathcal T(P, Z)\longrightarrow 0$$
i.e.,  $P$ is a $\xi$-projective object.

\vskip5pt

Let $P$ be a $\xi$-projective object. Since  $(\omega, \mathcal{T})$
is a complete cotorsion pair with respect to $\xi$,
there is a triangle $K \longrightarrow P'\stackrel v{\longrightarrow} P \longrightarrow K[1]$ in $\xi$ with $P'\in \omega$.
Since $P$ is $\xi$-projective, there is an exact sequences
 \ $0\longrightarrow \Hom_\mathcal T(P, K)\longrightarrow \Hom_\mathcal T(P, P')\stackrel{(P, v)}\longrightarrow \Hom_\mathcal T(P, P) \longrightarrow 0.$
Thus $v$ is a splitting epimorphism, and hence $P$ is a direct suammand of $P'\in \omega$. This proves
$P\in \omega$.

\vskip5pt

Conversely, assume that $\mathcal T$ has enough $\xi$-projective objects and $\omega = \mathcal P(\xi)$.  By Theorem \ref{mainthm}, the classes of cofibrant objects, of fibrant objects, and of trivial objects, of the model structure
$(\CoFib_{\omega}^\xi, \ \Fib_{\omega}^\xi, \ \Weq_{\omega}^\xi)$  are respectively $\mathcal X$,  $\mathcal T$ and $\Y$.
Comparing {\rm(\ref{Introconstruction})} with the definition of a $\xi$-triangulated model structure (cf. Definition \ref{defnxitriangulatedmds}) one sees that
the $(\xi, \omega)$-model structure $(\CoFib_{\omega}^\xi, \ \Fib_{\omega}^\xi, \ \Weq_{\omega}^\xi)$  is $\xi$-triangulated if and only if
$$\{\mbox{morphism} \ f  \mid  \Hom_\mathcal T(W, f) \ \mbox{is surjective}, \ \forall \ W\in\omega\} = \{\xi\mbox{-proper epimorphism}\}$$
and
$$\{\mbox{splitting monomorphism} \ f \ | \ \hoc f\in\omega\} = \{\xi\mbox{-proper monomorphism} \ f  \mid  \hoc f\in\omega\}.$$
\vskip5pt

For any morphism $f: A\longrightarrow B$ such that $\Hom_\mathcal T(W, f)$ is surjective for any $W\in \omega = P(\xi)$,
taking a $\xi$-proper epimorphism $v: P\longrightarrow B$ with  $P$ a $\xi$-projective object, then
$v = fh$ for some $h: P\longrightarrow A$. By Lemma \ref{monicepimorphism}, $f$ is a $\xi$-proper epimorphism.
On the other hand, for any $\xi$-proper epimorphism $f: A\longrightarrow B$ and for any object $W\in \omega = P(\xi)$, by the definition of a
$\xi$-projective object one knows that $\Hom_\mathcal T(W, f)$ is surjective. This proves the first equality. For the second equlaity, it suffices to show
the inclusion  $$\{\xi\mbox{-proper monomorphism} \ f  \mid  \hoc f\in\omega\}\subseteq \{\mbox{splitting monomorphism} \ f \ | \ \hoc f\in\omega\}.$$
This is clear  since $\omega = \mathcal P(\xi)$.
\end{proof}

\begin{cor}\label{gpctp} \ Let $(\T, [1], \tri)$ be a triangulated category,  $\xi$ a proper class of triangles,
and $\mathcal{X}$ and $\mathcal{Y}$ full additive subcategories of $\mathcal{T}$ which are closed under isomorphisms and direct summands. Put $\omega =\mathcal{X}\cap \mathcal{Y}$.
Then the following are equivalent.

\vskip5pt

$(1)$ \ The triple $(\CoFib_{\omega}^\xi, \ \Fib_{\omega}^\xi, \ \Weq_{\omega}^\xi)$ as defined in {\rm(\ref{Introconstruction})} is a $\xi$-triangulated model structure on $\T;$

\vskip5pt

$(2)$ \  $(\mathcal{X}, \mathcal{Y})$ is a hereditary complete cotorsion pair with respect to $\xi$, $\mathcal T$ has enough $\xi$-projective objects,  and $\omega = \mathcal P(\xi);$

\vskip5pt

$(3)$ \  $(\mathcal{X}, \mathcal{Y})$ is a generalized $\xi$-projective cotorsion pair with respect to $\xi$, and $\mathcal T$ has enough $\xi$-projective objects$;$

\vskip5pt

$(4)$ \ $(\CoFib_{\omega}^\xi, \ \Fib_{\omega}^\xi, \ \Weq_{\omega}^\xi)$ is a $\xi$-projective model structure.

\end{cor}
\begin{proof} \  $(1) \Longrightarrow (2)$: \ By Theorem \ref{onlyif},
$(\mathcal{X}, \mathcal{Y})$ is a hereditary complete cotorsion pair with respect to $\xi$ such that $\omega=\mathcal{X}\cap \mathcal{Y}$ contravariantly finite in $\T$.
By Proposition \ref{thesame}, $\mathcal T$ has enough $\xi$-projective objects and $\omega = \mathcal{P}(\xi)$.
\vskip5pt

\vskip5pt

$(2)\Longrightarrow (3)$: \ Since $\mathcal T$ has enough $\xi$-projective objects  and $\omega = \mathcal P(\xi)$, it follows that $\omega$ is contravariantly finite in $\T$.
By Proposition \ref{thesame},  $(\CoFib_{\omega}^\xi, \ \Fib_{\omega}^\xi, \ \Weq_{\omega}^\xi)$ is a $\xi$-triangulated model structure.
By Theorem \ref{mainthm}, $\Y$ is the class of trivial objects.
Thus $\Y$ is thick (cf. Theorem \ref{abelianhtpcat}). By definition $(\mathcal{X}, \mathcal{Y})$ is a generalized $\xi$-projective cotorsion pair with respect to $\xi$.

\vskip5pt

$(3) \Longrightarrow (4)$: \ By the assumption, $\mathcal T$ has enough $\xi$-projective objects, $\omega = \mathcal P(\xi)$, and $\Y$ is $\xi$-thick.
In particular,  $\Y$ is closed under the hocokernel of $\xi$-proper monomorphisms and $\omega$ is contravariantly finite in $\T$.
By Theorem \ref{heredity}, $(\mathcal{X}, \mathcal{Y})$ is a hereditary complete cotorsion pair with respect to $\xi$. By Proposition \ref{thesame},
$(\CoFib_{\omega}^\xi, \ \Fib_{\omega}^\xi, \ \Weq_{\omega}^\xi)$ is a $\xi$-triangulated model structure,
and moreover, each object is fibrant. By definition it is a $\xi$-projective model structure.

\vskip5pt

$(4) \Longrightarrow (1):$  \ This follows directly from the definition.\end{proof}

\subsection{$\xi$-triangulated model structures which are weakly $\xi$-projective}

\begin{prop} \label{xi-projective} \ {\rm (Compare [G, 4.6])} \
The intersection of the class of $\xi$-triangulated model structures and the class of weakly $\xi$-projective model structures
is the class of $\xi$-projective model structures.

\vskip5pt

More explicitly, let $(\T, [1], \tri)$ be a triangulated category,  $\xi$ a proper class of triangles,
and $(\CoFib, \Fib, \Weq)$ a $\xi$-triangulated model structure on $\mathcal T$, with {\rm Hovey} triple $(\mathcal C, \mathcal F, \mathcal W)$ with respect to $\xi$. Then the following are equivalent.

\vskip5pt

{\rm (1)} \ it is a weakly $\xi$-projective model structure$;$

\vskip5pt

{\rm (2)} \ it is a $\xi$-projective model structure$;$

\vskip5pt

{\rm (3)} \ $\mathcal T$ has enough $\xi$-projective objects
and the trivially cofibrant objects coincide with $\mathcal P(\xi);$

\vskip5pt

{\rm (4)} \ $(\mathcal{C}, \mathcal{W})$ is a hereditary complete cotorsion pair with respect to $\xi$,
$\mathcal T$ has enough $\xi$-projective objects, and $\mathcal{C}\cap \mathcal{W} = \mathcal P(\xi).$
\end{prop}
\begin{proof} The implication $(1) \Longrightarrow (2)$ follows directly from the definition.

\vskip5pt

$(2) \Longrightarrow (3):$ \ Since $(\CoFib, \Fib, \Weq)$ is a $\xi$-triangulated model structure, it satisfies the conditions in Proposition \ref{weaklyproj}(4), thus it is
a weakly $\xi$-projective model structure. Since weakly $\xi$-projective model structures are precisely the $(\xi, \omega)$-model structures (cf. Theorem  \ref{brcorrespondence}),
$(\CoFib, \Fib, \Weq)$ is both a $\xi$-triangulated model structure and a $(\xi, \omega)$-model structure.
By Proposition \ref{thesame}, $\mathcal T$ has enough $\xi$-projective objects
and the trivially cofibrant objects coincide with $\mathcal P(\xi)$.

\vskip5pt

 $(3) \Longrightarrow (4):$  \ Since $(\CoFib, \Fib, \Weq)$ is a $\xi$-triangulated model structure,
it follows from Theorem \ref{abelianhtpcat} that  $\mathcal W$ is $\xi$-thick,  $(\mathcal{C}\cap\mathcal W, \mathcal F) = (\mathcal{P}(\xi), \mathcal F)$
and $(\mathcal{C}, \mathcal F\cap\mathcal W)$
are complete cotorsion pairs with respect to $\xi$. We claim that $\mathcal{F} =  \mathcal{T}$.

\vskip5pt

In fact, for any object $P\in \mathcal P(\xi)$, by  the definition of a $\xi$-projective object, one knows that any triangle
 $K \longrightarrow M \longrightarrow P \longrightarrow K[1]$ in $\xi$ splits, and hence $\Ext_\xi(P, K) = 0, \ \forall \ K\in\mathcal T.$
Thus  $\Ext_\xi(\mathcal P(\xi), \T) = 0$, and hence $\T \subseteq \mathcal P(\xi)^\perp = \mathcal F$. Therefore $\mathcal{F} =  \mathcal{T}$.
This proves the claim.

\vskip5pt

Thus, $(\mathcal{C}, \mathcal W) = (\mathcal{C}, \mathcal F\cap\mathcal W)$ is a complete cotorsion pairs with respect to $\xi$.
In order to see the heredity, by Theorem \ref{heredity}, it suffices to see that
$\mathcal W$ is closed under the hocokernels of $\xi$-proper monomorphisms. This follows from the $\xi$-thickness of $\mathcal W$.

\vskip5pt

 $(4) \Longrightarrow (1):$  \  By definition, a $\xi$-triangulated model structure satisfies the conditions in Proposition \ref{weaklyproj}(4). Thus, it remains
to show each object is fibrant, i.e., $\mathcal{F} =  \mathcal{T}$.
For any object $P\in \mathcal P(\xi)$, by  the definition of a $\xi$-projective object, one knows that any triangle
 $K \longrightarrow M \longrightarrow P \longrightarrow K[1]$ in $\xi$ splits, and hence $\Ext_\xi(P, K) = 0, \ \forall \ K\in\mathcal T.$
Thus  $\Ext_\xi(\mathcal P(\xi), \T) = 0$. Since $\mathcal T$ has enough $\xi$-projective objects,
by Fact \ref{completectp}, $(\mathcal P(\xi), \mathcal{T})$
is a complete cotorsion pair with respect to $\xi$, i.e., $(\mathcal{C}\cap \mathcal W, \mathcal{T})$
is a complete cotorsion pair with respect to $\xi$. On the other hand,
by Theorem \ref{abelianhtpcat} $(\mathcal{C}\cap \mathcal W, \mathcal{F})$
is also a  cotorsion pairs with respect to $\xi$. It follows that $\mathcal{F} =  \mathcal{T}$. \end{proof}

\appendix
\section{\bf Fundamental Theorem of Model Categories}

\vskip5pt

As mentioned in Introduction, the proof of the following theorem is almost the same as the original one given by Quillen [Q1, Theorem 1'].
We only include the place where the condition that $\T$ is a model category can be replaced by the conditions
${\rm (i)}$, ${\rm (ii)}$,  ${\rm (iii)}$, ${\rm (iv)}$, ${\rm (v)}$ in {\rm Theorem \ref{appendixhtpthm}}, or more explicitly, the place where the
existence of pull-backs and push-outs can be replaced by the conditions (iv) and (v) in {\rm Theorem \ref{appendixhtpthm}}. For convenience, we will also explicitly mention
which parts in the original proof in [Q1] can be omitted.

\begin{thm}\label{appendixhtpthm} \ {\rm (Fundamental Theorem of Model Categories, \ [Q1, Theorem 1'])}

\vskip5pt  Let $\mathcal M$ be a category satisfying the following conditions$:$

\vskip5pt

${\rm (i)}$ \ \  $\mathcal M$ has the initial object and the final object$;$

${\rm (ii)}$ \  $\mathcal M$ has finite coproducts and finite products$;$

${\rm (iii)}$ \ \ There is a model structure $(\CoFib, \ \Fib,  \ \Weq)$ on $\mathcal M;$

${\rm (iv)}$ \ \ For any trivial cofibration $i: A\longrightarrow B$ and any morphism $u: A\longrightarrow C$,  there exists a weak push-out square
$$\xymatrix@R=0.5cm{
A\ar[r]^-u\ar[d]_-{i} & C \ar[d]^-{i'}\\
B\ar[r] & D
}$$
such that $i'$ is also a trivial cofibration.

\vskip5pt

${\rm (v)}$ \ \ For any trivial fibration $p: C\longrightarrow D$ and any morphism $v: B\longrightarrow D$,
there exists a weak pull-back square $$\xymatrix@R=0.5cm{
A\ar[r]\ar[d]_-{p'} & C \ar[d]^-{p}\\
B\ar[r]^-v & D
}$$ such that $p'$ is also a trivial fibration.

\vskip5pt

Then the composition of the embedding
 $\mathcal{M}_{cf}\hookrightarrow \mathcal{M}$ and the localization functor $\mathcal{M} \longrightarrow {\rm Ho}(\mathcal{M})$ induces an equivalence
 \ $\pi\mathcal{M}_{cf} \cong {\rm Ho}(\mathcal{M})$  of categories. \end{thm}

\begin{deflem}{\rm ([Q1, Definition 3, Definition 4, Lemma 1, Lemma 2])} \label{lhtp} \ Assume that $\mathcal M$ is a category
satisfying the conditions ${\rm (i)}$, ${\rm (ii)}$,  ${\rm (iii)}$, ${\rm (iv)}$, ${\rm (v)}$ in {\rm Theorem \ref{htpthm}}. Let $f, \ g \in \Hom_\mathcal M(A, B)$.

  \vskip5pt

  ${\rm (1)}$ \  There exists a commutative diagram
      $$\xymatrix@C=1.5cm{A\oplus A\ar[d]_-{(1,1)}\ar[rd]^-{(\partial_0,\partial_1)} \ar[r]^-{(f,g)} & B \\
      A & \widetilde{A}\ar[l]_-\sigma\ar[u]_-h }$$
      such that $\sigma$ is a weak equivalence if and only if there exists a commutative diagram
      $$\xymatrix@C=1.5cm{A\oplus A\ar[d]_-{(1,1)}\ar[rd]^-{(\partial_0,\partial_1)} \ar[r]^-{(f,g)} & B \\
      A & \widetilde{A}\ar[l]_-\sigma\ar[u]_-h }$$
      such that $\sigma$ is a weak equivalence and $(\partial_0, \partial_1)$ is a cofibration.

  \vskip10pt

If this is the case, one says that $f$ is {\it left homotopic} to $g$,  and denote it by $f\stackrel{l}\sim g$; and when $(\partial_0, \partial_1)$ is a cofibration, then $h:\widetilde{A}\longrightarrow B$ is called a left homotopy from $f$ to $g$.

  \vskip5pt

Moreover, if $A$ is a cofibrant object and if $(\partial_0, \partial_1)$ is a cofibration, then $\partial_0$ and $\partial_1$ are trivial cofibrations.

 \vskip10pt

 ${\rm (1')}$ \ There exists a commutative diagram
      $$\xymatrix@C=1.5cm{\widetilde{B}\ar[rd]^-{\binom{d_0}{d_1}} & B\ar[l]_-s \ar[d]^-{\binom{1}{1}} \\
      A\ar[u]^-k \ar[r]^-{\binom{f}{g}} & B\times B}$$
      such that $s$ is a weak equivalence if and only if there exists a commutative diagram
  $$\xymatrix@C=1.5cm{\widetilde{B}\ar[rd]^-{\binom{d_0}{d_1}} & B\ar[l]_-s \ar[d]^-{\binom{1}{1}} \\
      A\ar[u]^-k \ar[r]^-{\binom{f}{g}} & B\times B}$$
      such that $s$ is a weak equivalence and $\binom{d_0}{d_1}$ is a fibration.

  \vskip10pt

If this is the case, one says that $f$ is {\it right homotopic} to $g$ if $f$ and $g$, and denote it by $f\stackrel{r}\sim g$; an when
$\binom{d_0}{d_1}$ is a fibration, then $k: A\longrightarrow \widetilde{B}$ is called a right homotopy from $f$ to $g$.

\vskip5pt

Moreover, if $B$ is a fibrant object and if $\binom{d_0}{d_1}$ is a fibration, then $d_0$ and $d_1$ are trivial fibrations.

  \end{deflem}

 The proof of Definition - Lemma \ref{lhtp} is completely same as the original one in [Q1, Lemma 1, Lemma 2].
 Although in the proof of [Q1, Lemma 2], a push-out is used, but this push-out follows from the universal property of the coproduct.
 Namely, if $A$ is a cofibrant object, by the universal property of the coproduct one has the push-out
      $$\xymatrix{
          \emptyset\ar[d]\ar[r] & A\ar[d]^{\binom{1}{\emptyset}}\ar@/^25pt/[rdd]^{\partial_0} &  \\
          A\ar[r]^-{\binom{\emptyset}{1}}\ar@/_10pt/[rrd]^{\partial_1} & A\oplus A\ar@{-->}[rd]^-{(\partial_0,\partial_1)} & \\
            &  &  \widetilde{A} }$$
      Then by Fact \ref{elementpropmodel}(3),  $\binom{1}{\emptyset}: A \to A \oplus A$ is a cofibration. Thus $\partial_0 =(\partial_0, \partial_1)\binom{1}{\emptyset}$ is also a cofibration. By the following commutative diagram
      $$\xymatrix{ A\oplus A \ar[rd]^-{(\partial_0,\partial_1)}\ar[d]_-{(1,1)} \\
      A & \widetilde{A} \ar[l]_-\sigma }$$
      one sees that $\sigma \partial_0= \Id_A$. So by the two out of three axiom $\partial_0$ is also a weak equivalence, and hence $\partial_0$ is a trivial cofibration.
      Similarly $\partial_1$ is a trivial cofibration.

\vskip5pt

Also, the notion of cylinder object and path object in [Q1, Definition 4], Corollary and Lemma 3 in [Q1, p. 1.7] can be omitted.

 \begin{lem}{\rm ([Q1 Lemma 4])} \label{equivrel1} \ Assume that $\mathcal M$ is a category
satisfying the conditions ${\rm (i)}$, ${\rm (ii)}$,  ${\rm (iii)}$, ${\rm (iv)}$, ${\rm (v)}$ in {\rm Theorem \ref{htpthm}}.

\vskip5pt

  $(1)$ \ Let $A$ be a cofibrant object. Then $\overset{l}{\sim}$ is an equivalence relation on $\Hom_\mathcal M(A, B)$.

  \vskip5pt

  $(1')$ \ Let $B$ be a fibrant object. Then $\overset{r}{\sim}$ is an equivalence relation on $\Hom_\mathcal M(A, B)$.
  \end{lem}

  \begin{proof} \ (1) \ In fact, the proof of the transitivity is the unique place different from the original proof in [Q1]. Let $f_0\overset{l}{\sim}f_1$ and $f_1\overset{l}{\sim}f_2$. Then by Definition - Lemma \ref{lhtp} there exist the following commutative diagrams
    $$\xymatrix@C=1.5cm{A\oplus A\ar[d]_-{(1,1)}\ar[rd]^-{(\partial_0,\partial_1)} \ar[r]^-{(f_0,f_1)} & B && A\oplus A\ar[d]_-{(1,1)}\ar[rd]^-{(\partial_0',\partial_1')} \ar[r]^-{(f_1,f_2)} & B  \\
    A & A_0\ar[l]_-{\sigma_0}\ar[u]_-{h_0} && A & A_1\ar[l]_-{\sigma_1}\ar[u]_-{h_1}}$$
    such that $\sigma_0$ and $\sigma_1$ are weak equivalences, and $\partial_0, \partial_1,  \partial'_0, \partial'_1$ are trivial cofibrations.

    \vskip5pt

    For the trivial cofibration $\partial_0'$ and morphism $\partial_1$, by the condition (iv) in Theorem \ref{appendixhtpthm},  there exists a weak push-out
$$\xymatrix{A\ar[d]_-{\partial_0'}\ar[r]^{\partial_1} & A_0\ar[d]^{in_1}\\
        A_1\ar[r]^{in_2}& A_2}$$
such that $in_1$ is a trivial cofibration.

\vskip5pt

Since $\sigma_0 \partial_1 = {\rm Id}_A = \sigma_1 \partial_0'$, there exists a morphism $\sigma_2: A_2 \xra{} A$
such that $\sigma_2 \circ in_1 = \sigma_0$ and $\sigma_2 \circ in_2 = \sigma_1$.
Since $h_0 \partial_1 =f_1= h_1 \partial_0'$, there exists a morphism $h_2: A_2 \xra{} B$ such that $h_2 \circ in_1= h_0$ and $h_2\circ in_2= h_1$.
See commutative diagrams below.
    $$\xymatrix{A\ar[d]_-{\partial_0'}\ar[r]^{\partial_1} & A_0\ar[d]_-{in_1}\ar@/^10pt/[rdd]^{\sigma_0} &&&    A\ar[d]_-{\partial_0'}\ar[r]^{\partial_1} & A_0\ar[d]_-{in_1}\ar@/^10pt/[rdd]^{h_0} &  \\
        A_1\ar[r]^{in_2}\ar@/_10pt/[rrd]_-{\sigma_1} & A_2\ar@{-->}[rd]^{\sigma_2} &&&  A_1\ar[r]^{in_2}\ar@/_10pt/[rrd]_-{h_1} & A_2\ar@{-->}[rd]^{h_2} & \\
          &  &  A &&& &B }$$
    Thus one has the following commutative diagram
    $$\xymatrix@C=2cm{A\oplus A\ar[d]_-{(1,1)}\ar[rrd]^-{(in_1\circ\partial_0, \ in_2\circ\partial_1')} \ar[rr]^-{(f_0, \ f_2)} && B \\
    A && A_2\ar[ll]_-{\sigma_2}\ar[u]_-{h_2} }$$
    Since $\sigma_2 \circ in_1 = \sigma_0$ is a weak equivalence, by the two out of three axiom  $\sigma_2$ is also a weak equivalence. By definition $f_0\overset{l}{\sim}f_2$.

    \vskip5pt

    $(1')$ \ This can be dually proved, by using the condition (v) in Theorem \ref{appendixhtpthm}.  We omit the details. \end{proof}

From now on the proof of Theorem \ref{appendixhtpthm} is completely the same as Quillen's original proof, and we omit the details.


\vskip5pt

\begin{thebibliography}{99}	
\bibitem[AS]{AS}  J. Asadollahi, S. Salarian, Gorenstein objects in triangulated categories, J. Algebra 281(2004), 264-286.
\bibitem[BBD]{BBD}  A. A. Beilinson, J. Bernstein, P. Deligne, Faisceaux pervers, Ast\'erisgue 100. Soc. Math. France, Paris, 1982.

\bibitem[Bec]{Bec} H.Becker, Models for singularity categories, Adv. Math. 254(2014), 187-232.

\bibitem[B]{B} A. Beligiannis, Relative homological algebra and purity in triangulated categories, J. Algebra 227(1)(2000), 268-361.

\bibitem[BR]{BR}  A. Beligiannis, I. Reiten, Homological and homotopical aspects of torsion theories, Mem. Amer. Math. Soc. 188(883), 2007.


  \bibitem[Bon]{Bon}  M. V. Bondarko, Weight structures vs. $t$-structures; weight filtrations, spectral sequences, and complexes (for motives and in general), J. K-Theory 6(3)(2010), 387-504.
Preprint available at http://arxiv.org/abs/0704.4003v1

\bibitem[CLZ]{CLZ} J. Cui, X. S. Lu, P. Zhang, Model structure from one hereditary complete cotorsion pair, J. Pure Appl. Algebra (to appear),
arXiv: 2401.08078v2.


\bibitem[CRZ]{CRZ} J. Cui, S. Rong, P. Zhang, Cotorsion pairs and model structures on Morita rings, J. Algebra 661(2025), 1-81.



\bibitem[DS]{DS} W. G. Dwyer, J. Spalinski, Homotopy theories and model categories, in: Handbook of algebraic topology (Amsterdam), North-Holland, Amsterdam, 1995, pp. 73-126.

\bibitem[G1]{G1} J. Gillespie, Model structures on exact categories, J. Pure Appl. Algebra 215(12)(2011), 2892-2902.
\bibitem[G2]{G2} J. Gillespie, How to construct a Hovey triple from two cotorsion pairs, Fund. Math. 230(3)(2015), 281-289.
\bibitem[G3]{G3} J. Gillespie, Abelian model category theory, Cambridge Studies in Adv. Math. 215, Cambridge Univ. Press, 2025.
 \bibitem[Hap]{Hap} D. Happel, Triangulated categories in the representation theory of finite dimensional algebras, Lond. Math. Soc Lecture Note Ser. 119,
  Cambridge, New York, New Rochelle, Melbourne, Sydney, Cambridge Univ. Press, 1988.

\bibitem[Hov1]{Hov1} \ M. Hovey, Model categories,  Math. Surveys and Monographs 63, Amer. Math. Soc., Providence, 1999.
\bibitem[Hov2]{Hov2} \ M. Hovey, Cotorsion pairs, model category structures, and representation theory,  Math. Z. 241(3)(2002), 553-592.
\bibitem[Hir]{Hir} P. S. Hirschhorn, Model categories and their localizations, Math. Surveys and
Monographs 99, Amer. Math. Soc., Providence, 2003.

\bibitem[Hub]{Hub} A. Hubery, Notes on the octahedral axiom, Available at http://www.maths.leeds.ac.uk/ahubery/octahedral.pdf.

\bibitem[HZZ]{HZZ} J. S. Hu, D. D. Zhang, P. Y. Zhou, Model structure arising from one hereditary cotorsion pair on extriangulated categories, arXiv: 2406.14031v1.

\bibitem[M]{M} H. R. Margolis, Spectra and Steenrod Algebra, North Holland Math. Library, Vol. 29, North-Holland, Amsterdam, 1983.

\bibitem[NP]{NP} H. Nakaoka, Y. Palu, Extriangulated categories, Hovey twin cotorsion pairs and model structures, Cahiers
de Topologie et Geometrie Differentielle Categoriques, Volume LX-2 (2019), 117-193.

\bibitem[N1]{N1} A. Neeman, Some new axioms for triangulated categories, J. Algebra 139(1)(1991), 221-255.

\bibitem[N2]{N2} A. Neeman, Triangulated categories, Annals of Math. Studies 148, Princeton Univ. Press, Princeton, NJ, 2001.

\bibitem[P]{P} D. Pauksztello, Compact corigid objects in triangulated categories and co-t-structures, Cent. Eur. J. Math. 6(1) (2008) 25-42.

\bibitem[Q1]{Q1} \ D. Quillen, Homotopical algebra, Lecture Notes in Math. 43, Springer-Verlag, 1967.

\bibitem[Q2]{Q2} \ D. Quillen, Rational Homotopy Theory, Ann. Math. 90(2)(1969), 205-295.

\bibitem[\v{S}]{S} \ J. \v{S}t'ov\'{\i}\v{c}ek, Exact model categories, approximation theory, and cohomology of quasi-coherent sheaves, in: Advances in Representation Theory of Algebras, EMS Series of Congress Reports, European Math. Soc. Publishing House, 2014, pp. 297 - 367.

\bibitem[V]{V} J. L. Verdier, Cat\'egories d\'eriv\'ees, \'etat 0 in SGA $4\frac{1}{2}$//Lectures Notes in Math. 569, 262-311, Berlin: Springer-Verlag, 1977.

  \bibitem[Y]{Y} X. Y. Yang, Model structures on triangulated categories, Glasg. Math. J. 57(2)(2015), 263-284.

\end{thebibliography}
\end{document}